\documentclass[12pt]{amsart}
\usepackage{amsmath,amscd,amssymb,amsfonts}
\def\ssbull{\raise.2ex\hbox{${\scriptscriptstyle\bullet}$}}
\def\scirc{\,\raise.2ex\hbox{${\scriptstyle\circ}$}\,}
\def\mopls{\hbox{$\bigoplus$}}
\def\({{\rm (}}
\def\){{\rm )}}
\def\bA{{\mathbb A}}
\def\bC{{\mathbb C}}
\def\bD{{\mathbb D}}
\def\bP{{\mathbb P}}
\def\bQ{{\mathbb Q}}
\def\bR{{\mathbb R}}
\def\bZ{{\mathbb Z}}
\def\cA{{\mathcal A}}
\def\cC{{\mathcal C}}
\def\cD{{\mathcal D}}
\def\cF{{\mathcal F}}
\def\cH{{\mathcal H}}
\def\cM{{\mathcal M}}
\def\cO{{\mathcal O}}
\def\cV{{\mathcal V}}
\def\For{\hbox{{\rm For}}}
\def\real{\hbox{{\rm real}}}
\def\can{\hbox{{\rm can}}}
\def\sp{\hbox{{\rm sp}}}
\def\CH{\hbox{{\rm CH}}}
\def\IC{\hbox{{\rm IC}}}
\def\Im{\hbox{{\rm Im}}}
\def\End{\hbox{{\rm End}}}
\def\Ext{\hbox{{\rm Ext}}}
\def\Gr{\hbox{{\rm Gr}}}
\def\Gal{\hbox{{\rm Gal}}}
\def\Coker{\hbox{{\rm Coker}}}
\def\Ker{\hbox{{\rm Ker}}}
\def\Hom{\hbox{{\rm Hom}}}
\def\MF{\hbox{{\rm MF}}}
\def\MFW{\hbox{{\rm MFW}}}
\def\MHM{\text{\rm MHM}}
\def\Perv{\hbox{\rm Perv}}
\def\Sym{\text{\rm Sym}}
\def\Var{\text{\rm Var}}
\def\Spec{\text{{\rm Spec}}\,}
\def\supp{\hbox{{\rm supp}}\,}
\def\red{\text{\rm red}}
\def\gl{\text{\rm gl}}
\def\go{\text{\rm go}}
\def\simto{\buildrel{\sim\,\,}\over\to}

\begin{document}
\title{On the Formalism of Mixed Sheaves}
\author{Morihiko Saito}
\address{RIMS Kyoto University, Kyoto 606-8502 Japan}
\maketitle

\centerline{\bf Introduction}

\bigskip\noindent
It has been conjectured by Beilinson [1], etc. that there would exist the
category of {\it mixed motivic sheaves}
$ \cM\cM(X) $ for an algebraic variety
$ X $ defined over an algebraic number field
$ k, $ such that the motivic cohomology of
$ X $ with rational coefficients can be calculated by the higher extension
group of the constant sheaves on
$ X $ if
$ X $ is smooth proper.
For example,
$$
\CH^{p}(X)_{\bQ} = \Ext_{D^{b}\cM\cM(X)}^{2p}(\bQ_{X}^{\cM},
\bQ_{X}^{\cM}(p))
\leqno(0.1)
$$
should hold for
$ X $ smooth proper, where
$ \CH^{p}(X)_{\bQ} $ is the Chow group of
$ X $ with rational coefficients, and
$ \bQ_{X}^{\cM} $ denotes the constant object in
$ D^{b}\cM\cM(X), $ the bounded derived category of
$ \cM\cM(X) $.
In this paper we give a formalism of mixed sheaves which might be useful
for the construction of mixed motivic sheaves.

Let
$ k $ be a field embeddable in
$ \bC, $ and
$ \cV(k) $ the category of separated algebraic varieties over
$ k $.
We choose and fix an embedding
$ k \to \bC $.
For
$ X \in \cV(k), $ let
$ X(\bC) $ denote the set of closed points of
$ X_{\bC} \,(:= X\otimes_{k}\bC) $ with classical topology.
For a field
$ A $ of characteristic zero, let
$ {D}_{c}^{b}(X(\bC), A) $ denote the bounded derived categories of
complexes of
$ A $-Modules on
$ X(\bC) $ with constructible cohomologies, and
$ \Perv(X(\bC), A) $ the abelian category of
$ A $-perverse sheaves [3] on
$ X(\bC) $ which is a full subcategory of
$ {D}_{c}^{b}(X(\bC),A) $.

A {\it theory of}
$ A $-{\it mixed sheaves on}
$ \cV(k) $ consists of
$ A $-linear abelian categories
$ \cM(X) $ with
$ A $-linear forgetful functors
$ \For : \cM(X) \to \Perv(X(\bC), A) $ for
$ X \in \cV(k) $ such that the categories
$ \cM(X) $ are stable by the dual functor
$ \bD, $ the external product
$ \boxtimes $, the pull-backs
$ j^{*} $ by open embeddings
$ j, $ and the (cohomological) direct images
$ H^{i}f_{*} $ by affine morphisms
$ f, $ in a compatible way with the functor
$ \For $, and
$ \cM(X) $ contains a constant object
$ A_{X}^{\cM}[n] $ if
$ X $ is smooth of pure dimension
$ n $.
Moreover they should satisfy natural properties associated with these
functors.
Then we show that the bounded derived categories
$ D^{b}\cM(X) $ are stable by the standard functors
$ f_{*}, f_{!}, f^{*}, f^{!}, \psi_{g}, \varphi_{g,1}, \bD,\boxtimes,
\otimes , \cH{om} $ in a compatible way with the corresponding
functors on the underlying complexes of
$ A $-Modules by the forgetful functor, and they satisfy natural
properties associated with these functors.
Using these, most of the arguments in [28] [30] [31] are valid in this
setting.
Exception comes mainly from the Hodge theoretic description of the Picard
variety of a smooth proper variety over
$ \bC, $ which becomes quite nontrivial in other situations
(e.g., for a variety over a number field).

As a corollary, the notion of {\it geometric origin} can be defined as in
[3]; the full subcategory
$ \cM(X)^{\go} $ of
$ \cM(X) $ is obtained by iterating the cohomological standard functors
as above to the constant object
$ \bQ^{\cM} $ on
$ \Spec k $ and also by taking subquotients in
$ \cM(X) $.
In the case of varieties over a number field
$ k, $ it is interesting whether we can get the category of mixed motivic
sheaves in this way.
For example, for
$ A = \bQ, $ we can define
$ \cM(X) $ by the category consisting of objects
$ ((M, F), K_{\sigma}, K_{l}; W) $ where
$ (M,F) $ is a filtered regular holonomic
$ \cD $-Module on
$ X, K_{\sigma} $ is a
$ \bQ $-perverse sheaf on
$ X_{\sigma}(\bC) $ with
$ X_{\sigma} := X\otimes_{k,\sigma}\bC $ for each embedding
$ \sigma : k \to \bC, K_{l} $ is an \'etale
$ \bQ_{l} $-perverse sheaf [3] on
$ \overline{X} := X\otimes_{k}\overline{k} $ with an action of
$ \Gal(\overline{k}/k) $ (compatible with the natural action on
$ \overline{X}) $ for each prime number
$ l, $ and
$ W $ is a finite filtration on
$ M, K_{\sigma}, K_{l} $,
such that they have comparison isomorphisms as in [11] [12] [14] [20]
(using [3]), and
$ ((M\otimes_{k,\sigma}\bC, F), K_{\sigma}; W) $ is a mixed Hodge Module
on
$ X_{\sigma} $ for any
$ \sigma $.
See (1.8) below.
Here we assume also that each graded piece
$ \Gr_{k}^{W}((M, F), K_{\sigma}, K_{l}) $ has a pairing which induces a
polarization of the mixed Hodge Module on
$ X_{\sigma} $.
(This condition is satisfied for the objects of geometric origin.)
Then they satisfy the axiom of mixed sheaves, and the full subcategory
consisting of the objects of geometric origin might be close to the
category of mixed motivic sheaves.
In this example, the Hodge conjecture is replaced with the conjecture in
[20] that an absolute Hodge cycle [11] [14] is algebraic, which will be
called the {\it Hodge type conjecture}.

Note that, if the motivic sheaves exist so that (0.1) holds,
then the Hodge type conjecture should follow from (0.1).
In fact, we have a natural surjective morphism
$$
\Ext_{D^{b}\cM\cM(X)}^{2p}(\bQ_{X}^{\cM}, \bQ_{X}^{\cM}(p)) \to
\Hom_{\cM\cM(\Spec k)}(\bQ^{\cM}, H^{2p}(X, \bQ(p))^{\cM})
\leqno(0.2)
$$
by the adjunction for
$ a_{X} : X \to \Spec k $ and the decomposition theorem for
$ (a_{X})_{*}\bQ_{X}^{\cM}(p), $ where
$ H^{2p}(X, \bQ(p))^{\cM} := H^{2p}(a_{X})_{*}\bQ_{X}^{\cM}(p) $ and the
target of (0.2) should be the group of absolute Hodge cycles by the second
hypothesis.

In this note we show a partial converse.
We construct a cycle map
$$
\CH^{p}(X)_{\bQ} \to \Ext_{D^{b}M(X{)}^{\go}}^{2p}(\bQ_{X}^{\cM},
\bQ_{X}^{\cM}(p))
\leqno(0.3)
$$
and reduce its surjectivity to the Hodge type conjecture for any smooth
projective varieties over
$ k $ (see (8.9)).
The injectivity of (0.3) seems to be related with the surjectivity of the
cycle map of Bloch's higher Chow group
$ \CH^{p}(X,1)_{\bQ} $ to
$ \Ext^{2p-1}(\bQ_{X}^{\cM},\bQ_{X}^{\cM}(p)) $ (see (8.13)).
Note that the injectivity can be reduced to the case
$ k = \bC $ and
$ \cM(X) = \MHM(X) $ (the category of mixed Hodge Modules on
$ X), $ because the natural map
$ \CH^{p}(X)_{\bQ} \to \CH^{p}(X_{\bC})_{\bQ} $ is injective (see (8.5.2)).
The injectivity of the cycle map (0.3) would imply Bloch's conjecture [5]
(see [30]).
Its bijectivity would be related with Deligne's remark [19, 4.16] and
Murre's results [26] (cf. [30, II, (3.4)]).

As a corollary of the above construction, we get also, for example, a
natural
$ k $-structure on the intersection cohomology of an irreducible variety
defined over
$ k, $ and furthermore, a complex of
$ \cM(\Spec k) $ whose cohomology gives the intersection cohomology, where
$ \cM(\Spec k) $ is a full subcategory of the category consisting of
the families of realizations [12] [14] [20].
This holds also for the intersection cohomology (or cohomology with, or
without, compact supports) with coefficients in a local system of
geometric origin.

Some of the arguments in this paper are similar to [28] [30] [31].
But we repeated some arguments because of the changes coming from the
difference of varieties over
$ k $ and
$ \bC $ (for example,
$ X_{\bC} $ may be reducible even if
$ X $ is irreducible, and this makes the calculation of
$ \End(\bQ_{X}^{\cM}) $ more complicated).

In \S 1, we formulate the axiom of mixed sheaves.
In \S 2--5, we construct the direct image, pull-backs, dual, and nearby
and vanishing cycle functors.
In \S 6, we discuss about the weight.
In \S 7, we introduce the notion of geometric origin.
In \S 8, we apply these to the study of cycles maps.

\bigskip\bigskip
\centerline{\bf \S 1. Axiom of Mixed Sheaves}

\bigskip
\noindent
{\bf 1.1.}
Let
$ k $ be a field embeddable in
$ \bC $.
Let
$ \cV(k) $ be the category of separated algebraic varieties over
$ k, $ and
$ \cV $ a full subcategory of
$ \cV(k) $ such that:

\noindent
(i) the objects of
$ \cV $ are stable by cartesian products in
$ \cV(k) $,

\noindent
(ii) for an open embedding
$ j : X \to Y $ of
$ \cV(k) $,
$ X \in \cV $ if
$ Y \in \cV $,

\noindent
(iii)
$ \bP^{1} $ and
$ \Spec k $ belong to
$ \cV $.

\noindent
Actually, we may assume
$ \cV $ is
$ \cV(k) $ or the full subcategory consisting of smooth varieties
of
$ \cV(k) $.

We choose and fix an embedding
$ k \to \bC $.
Let
$ A $ be a field of characteristic zero.
For
$ X \in \cV(k), $ let
$ X(\bC), {D}_{c}^{b}(X(\bC), A), $ and
$ \Perv(X(\bC),A) $ be as in the introduction.
We say that an abelian category is an
$ A $-linear abelian category, if
$ \Hom(M,N) $ has a structure of
$ A $-module for objects
$ M, N $ and the composition of morphisms is compatible with the structure
of
$ A $-module.
A functor of an
$ A $-linear abelian categories is called an
$ A $-linear functor if the map of morphisms is
$ A $-linear.
In this paper, a functor of
$ A $-linear abelian categories is assumed to be an
$ A $-linear functor.

A theory of
$ A $-mixed sheaves on
$ \cV $ consists of
$ A $-abelian categories
$ \cM(X) $ with
$ A $-linear functors
$ \For : \cM(X) \to \Perv(X(\bC), A) $ for
$ X \in \cV $ such that:

\begin{itemize}
\item[(1.1.1)]
The functor
$ \For $ is faithful
and exact.

\smallskip
\item[(1.1.2)]
$ \For(M) $ is quasi-unipotent in the sense of [21], and
$ \For(M) $ has a stratification defined
over
$ k $ such that the restrictions of the cohomology sheaves of
$ \For(M) $ to each stratum are local systems.

\smallskip
\item[(1.1.3)]
Each
$ M \in \cM(X) $ has a finite increasing
filtration
$ W $, called the {\it weight filtration}, which
is strictly compatible with any morphism of
$ \cM(X) $ (i.e.,
$ M \to \Gr_{i}^{W}M $ is an exact functor).

\smallskip
\item[(1.1.4)]
The graded pieces
$ \Gr_{i}^{W}M $ are semisimple for
$ M \in \cM(X) $.

\end{itemize}

Moreover they should satisfy the properties in (1.2--6) below.

\medskip\noindent
{\it Remark.} We will denote also by
$$
\For : D^{b}\cM(X) \to {D}_{c}^{b}(X(\bC), A)
\leqno(1.1.5)
$$
the composition of
$ \For : D^{b}\cM(X) \to D^{b} \Perv(X(\bC), A) $ with
$ \real : D^{b} \Perv(X(\bC), A) $
$\to {D}_{c}^{b}(X(\bC), A) $ in [3], and
$ \For(M) $ is called the underlying
$ A $-complex (or perverse sheaf) of
$ M \in D^{b}\cM(X) $ (or
$ \cM(X)), $ where
$ D^{b} $ denotes the bounded derived category (associated with an abelian
category).
See [35].
Then
$$
\For\scirc H^{k} ={}^{p}\cH^{k}\scirc \For
\leqno(1.1.6)
$$
by [3].
In particular, we get

\begin{itemize}
\item[(1.1.7)]
A morphism of
$ D^{b}\cM(X) $ is an isomorphism
if and only if so is its image by
$ \For $,
\end{itemize}
because we have by (1.1.1):
\begin{itemize}
\item[(1.1.8)]
$ \For : \cM(X) \to \Perv (X(\bC), A) $ commutes with
$ \Ker, \Im, \Coker $.
\end{itemize}

\medskip\noindent
{\bf 1.2.}~{\it Dual.} There exists a contravariant functor
$$
\bD : \cM(X) \to \cM(X)
\leqno(1.2.1)
$$
with a functorial isomorphism
$$
\For\scirc \bD = \bD\scirc \For ,
\leqno(1.2.2)
$$
such that
$$
W_{i}(\bD M) = \bD(M/W_{-i-1}M)
\leqno(1.2.3)
$$
for
$ M \in \cM(X) $.
There exists a functorial isomorphism
$$
\bD^{2} = id
\leqno(1.2.4)
$$
which corresponds to the natural isomorphism on the underlying perverse
sheaves by the functor
$ \For $.

\medskip\noindent
{\bf 1.3.}~{\it Open Pull-back.} For an open embedding
$ j : X \to Y $ of varieties of
$ \cV, $ there exists a functor
$$
j^{*} : \cM(Y) \to \cM(X)
\leqno(1.3.1)
$$
with a functorial isomorphism
$$
\For\scirc j^{*} = j^{*}\scirc \For ,
\leqno(1.3.2)
$$
such that
$$
W_{i}(j^{*}M) = j^{*}(W_{i}M)
\leqno(1.3.3)
$$
for
$ M \in \cM(X) $.
We denote
$ j^{*} $ also by
$ j^{!} $.
Then there exists a functorial isomorphism
$$
\bD\scirc j^{!} = j^{*}\scirc \bD
\leqno(1.3.4)
$$
which corresponds to the natural isomorphism on the underlying perverse
sheaves by the functor
$ \For $.

If
$ j = id, $ there exists a functorial isomorphism
$$
j^{*} = id
\leqno(1.3.5)
$$
which corresponds to the natural isomorphism on the underlying perverse
sheaves by the functor
$ \For $.

For open embeddings
$ j : X \to Y, j' : Y \to Z, $ there exists a functorial
isomorphism
$$
j^{*}\scirc j'^{*} = (j'j)^{*}
\leqno(1.3.6)
$$
which corresponds to the natural isomorphism on the underlying perverse
sheaves by the functor
$ \For $.

\medskip\noindent
{\it Remark.} If
$ j : X \to Y $ is an isomorphism of varieties of
$ \cV $,
(1.3.1) is an equivalence of categories by (1.3.5--6).

\medskip\noindent
{\bf 1.4.}~{\it Direct Image.} For an affine morphism
$ f : X \to Y $ of varieties of
$ \cV, $ there exists a cohomological functor
$$
H^{k}f_{*} : \cM(X) \to \cM(Y)
\leqno(1.4.1)
$$
with a functorial isomorphism
$$
\For\scirc H^{k}f_{*} ={}^{p}\cH^{k}f_{*}\scirc \For ,
\leqno(1.4.2)
$$
such that
$$
W_{i+k}(H^{k}f_{*}M) = 0\quad \text{if}\quad W_{i}M = 0
\leqno(1.4.3)
$$
for
$ M \in \cM(X) $.
For an open embedding
$ j : Y' \to Y, $ let
$ X' = f^{-1}(Y') $ so that we have a cartesian diagram:
$$\CD
X @>{j'}>> X' \\
@V{f}VV @VV{f'}V \\
Y @>>{j}> Y'
\endCD
$$
Then there exists a functorial isomorphism
$$
j^{!}\scirc H^{k}f_{*} = H^{k}f'_{*}\scirc j'^{!}
\leqno(1.4.4)
$$
which corresponds to the natural isomorphism on the underlying perverse
sheaves by the functor
$ \For $.

If
$ f $ is a closed embedding or an affine open embedding, we define
$ f_{*} $ by
$ H^{0}f_{*}, $ because
$ H^{k}f_{*} = 0 $ for
$ k \ne 0 $ by (1.4.2).
Then, for a closed embedding
$ i : X \to Y $ of varieties of
$ \cV, $ there exists a functorial isomorphism
$$
\bD\scirc i_{*} = i_{*}\scirc \bD
\leqno(1.4.5)
$$
which corresponds to the natural isomorphism on the underlying perverse
sheaves by the functor
$ \For $, and
$$
i_{*} : \cM(X) \to \cM_{X}(Y)
\leqno(1.4.6)
$$
is an equivalence of categories, where
$ \cM_{X}(Y) $ is the full subcategory of
$ \cM(Y) $ consisting of the objects
$ M $ such that
$ \supp M \subset i(X) $ (i.e.,
$ j^{*}M = 0 $ for
$ j : Y \setminus i(X) \to Y) $.

If
$ j : X \to Y $ is an affine open embedding of varieties of
$ \cV, $ there exist functorial isomorphism and morphism:
$$
j^{*}j_{*} = id\quad \text{and}\quad id \to j_{*}j^{*},
\leqno(1.4.7)
$$
which correspond to the natural isomorphism and morphism on the underlying
perverse sheaves by the functor
$ \For $.

Let
$ f : X \to Y $ and
$ g : Y \to Z $ be affine morphisms of varieties of
$ \cV $.
Then, for
$ M \in \cM(X) $ such that
$ H^{i}f_{*}M = 0 $ for
$ i \ne 0, $ there exists a functorial isomorphism
$$
H^{i}g_{*}H^{0}f_{*}M = H^{i}(gf)_{*}M,
\leqno(1.4.8)
$$
which corresponds to the natural isomorphism on the underlying perverse
sheaves by the functor
$ \For $.

\medskip\noindent
{\it Remarks.} (i) By (1.4.6), we have
$$
\cM(X_{\red}) = \cM(X).
\leqno(1.4.9)
$$

\medskip

(ii) If
$ j = id, (1.3.5) $ and (1.4.7) imply a canonical isomorphism
$$
j_{*} = id.
\leqno(1.4.10)
$$
If
$ f $ is an isomorphism,
$ f_{*} $ is an equivalence of categories and
$ f^{*} $ is a quasi-inverse of
$ f_{*} $ by (1.4.7--8) (1.4.10).

\medskip\noindent
{\bf 1.5.}~{\it External Product.} There exists a bifunctor
$$
\boxtimes : \cM(X) \times \cM(Y) \to \cM (X\times Y)
\leqno(1.5.1)
$$
with a functorial isomorphism
$$
\For(M\,\boxtimes\,N) = \For(M)\,\boxtimes\, \For(N)
\leqno(1.5.2)
$$
for
$ M \in \cM(X), N \in \cM(Y), $ such that
$$
W_{n}(M\,\boxtimes\,N) = \sum_{i+j=n} W_{i}M\,\boxtimes\,W_{j}N,
\leqno(1.5.3)
$$
where
$ W_{i}M\,\boxtimes\,W_{j}N $ is a subobject of
$ M\,\boxtimes\,N, $ because
$ \boxtimes $ is exact for both factors by (1.5.2) and (1.1.1).
There exists a functorial isomorphism
$$
(M\,\boxtimes\,N)\,\boxtimes\,L = M\,\boxtimes\,(N\,\boxtimes\,L)
\leqno(1.5.4)
$$
for
$ M \in \cM(X), N \in \cM(Y), L \in \cM(Z), $ which corresponds to
the natural isomorphism on the underlying perverse
sheaves.

Let
$ \iota : X\times Y \to Y\times X $ be an isomorphism defined
by
$ \iota (x,y) = (y,x) $.
Then there exists a functorial isomorphism
$$
\iota_{*}(M\,\boxtimes\,N) = N\,\boxtimes\,M
\leqno(1.5.5)
$$
for
$ M \in \cM(X), N \in \cM(Y), $ which corresponds to the
natural isomorphism on the underlying perverse sheaves.

For an open embedding
$ j : X' \to X $ and an affine morphism
$ f : X \to Z $ of varieties of
$ \cV, $ there exist functorial isomorphisms
$$
\bD(M\,\boxtimes\,N) = \bD M\,\boxtimes\,\bD N,
\leqno(1.5.6)
$$
$$
(j\times id)^{*}(M\,\boxtimes\,N) = j^{*}M\,\boxtimes\,N,
\leqno(1.5.7)
$$
$$
H^{k}(f\times id)_{*}(M\,\boxtimes\,N) = H^{k}f_{*}M\,\boxtimes\,N
\leqno(1.5.8)
$$
for
$ M \in \cM(X), N \in \cM(Y), $ which correspond to the natural
isomorphisms on the underlying perverse sheaves.

\medskip\noindent
{\bf 1.6.}~{\it Constant Object.} There exists an object
$ A^{\cM} $ of
$ \cM(\Spec k) $ with an isomorphism
$$
\For(A^{\cM}) = A,
\leqno(1.6.1)
$$
such that
$$
\Gr_{i}^{W}A^{\cM} = 0\quad \text{for}\quad i \ne 0.
\leqno(1.6.2)
$$
There exists an isomorphism
$$
\bD(A^{\cM}) = A^{\cM}
\leqno(1.6.3)
$$
which corresponds to the natural isomorphism on the underlying vector
spaces by the functor
$ \For $.

There exists a functorial isomorphism
$$
A^{\cM}\,\boxtimes\,M = M
\leqno(1.6.4)
$$
for
$ M \in \cM(X), $ which corresponds to the natural isomorphism
on the underlying perverse sheaves by the functor
$ \For $, where
$ \Spec k \times X $ is naturally identified with
$ X, $ and
$ \cM(\Spec k \times X) $ with
$ \cM(X) $ using Remarks after (1.3) and (1.4).

For a smooth variety
$ X \in \cV $ of pure dimension
$ n, $ there exists
$ A_{X}^{\cM}[n] \in \cM(X) $ with an isomorphism
$$
\For(A_{X}^{\cM}[n]) = A_{X(\bC)}[n],
\leqno(1.6.5)
$$
and a morphism
$$
A^{\cM} \to H^{-n}(a_{X})_{*}(A_{X}^{\cM}[n])
\leqno(1.6.6)
$$
which corresponds to the natural morphism
$ A \to H^{0}(a_{X(\bC)})_{*}A_{X(\bC)} $ by the functor
$ \For $, where
$ a_{X} : X \to \Spec k $ and
$ a_{X(\bC)} : X(\bC) \to \Spec \bC $ are natural
morphisms.

Let
$ \delta : X \to X\times X $ denote the diagonal morphism for
a smooth pure dimensional variety
$ X \in \cV $.
Then, for
$ M \in \cM(X), $ there exists a morphism
$$
\bD M\,\boxtimes\,M \to \delta_{*}\bD A_{X}^{\cM}
\leqno(1.6.7)
$$
in the bounded derived category
$ D^{b}\cM(X), $ which corresponds to the natural pairing of
the underlying perverse sheaf
$ \bD \For(M)\otimes \For(M) \to \bD A_{X(\bC)} $ by the functor
$ \For $ and the adjunction for
$ \delta , $ where
$ A_{X}^{\cM} = (A_{X}^{\cM}[n])[-n] \in D^{b}\cM(X), $ and
$ \delta_{*}, \bD $ are naturally extended to the derived categories
(because they are exact functors).
Moreover, the morphism
$$
\Hom_{\cM(X)}(N,\bD M) \to \Hom_{D}b_{\cM(X\times X)}
(N\,\boxtimes\,M, \delta_{*}\bD A_{X}^{\cM})
\leqno(1.6.8)
$$
obtained by the composition with (1.6.7) is bijective for
$ M, N \in \cM(X) $.

\medskip\noindent
{\bf 1.7.}~{\it Remark.} Since the functor
$ \For : \cM(X) \to \Perv(X(\bC), A) $ is faithful, a morphism of
$ \cM(X) $ is uniquely determined by its underlying morphism of perverse
sheaves.
So the compatibility of some morphisms of
$ \cM(X) $ can be reduced to that for the underlying perverse sheaves.

\medskip\noindent
{\bf 1.8.}~{\it Examples.}
We are interested in the case
$ \cV = \cV(k) $.
By (6.15) below, we may assume
$ \cV $ is the full subcategory of smooth varieties.

\medskip

(i) A basic example is:
$ \cM(X) = \MHM(X_{\bC},A) $ the category of mixed Hodge Modules on
$ X_{\bC} $ with
$ A $-coefficients [27] [28], where
$ A $ is a subfield of
$ \bR $.

\medskip

(ii) Let
$ \MF_{rh}(\cD_{X}) $ be the category of regular holonomic
$ \cD_{X} $-Modules
$ M $ with a good filtration
$ F $ (i.e.,
$ M\otimes_{k}\bC $ is regular holonomic on
$ X_{\bC} $ (cf. Remark below) and
$ \Gr^{F}M $ is coherent over
$ \Gr^{F}\cD_{X}), $ and
$ \MFW_{rh}(\cD_{X}) $ the category of
$ (M,F,W) $ where
$ (M,F) \in \MF_{rh}(\cD_{X}) $ and
$ W $ is a finite filtration of
$ M $.
Here we assume varieties are smooth using (6.15).
Then we have a natural functor
$$
\MFW_{rh}(\cD_{X}) \to \MFW_{rh}(\cD_{X_{\bC}})
\leqno(1.8.1)
$$
If
$ A $ is a subfield of
$ \bR, $ we can define
$ \cM(X) $ using the fiber product of
$ \MHM(X_{\bC},A) $ and
$ \MFW_{rh}(\cD_{X}) $ over
$ \MFW_{rh}(\cD_{X_{\bC}}) $ (i.e., the underlying
$ \cD_{X_{\bC}} $-Module of a mixed Hodge Module is defined over
$ X $ as a
$ \cD_{X} $-Module with filtrations
$ F $ and
$ W) $.
Here we assume that a polarization on the graduation by the weight
filtration
$ W $ of the mixed Hodge Module is defined also on the underlying filtered
$ \cD_{X} $-Module
$ \Gr_{i}^{W}(M,F) $ as the self duality
$ \bD(\Gr_{i}^{W}(M,F)) \simeq (\Gr_{i}^{W}(M,F))(i) $, where the dual
functor
$ \bD $ and Tate twist
$ (i) $ are defined as in [27, \S 2].
Then we can check the conditions of mixed sheaves.
For example, (1.6.8) for filtered
$ \cD_{X} $-Modules follows from the same argument as in [34], using the
resolution by canonical induced filtered
$ \cD $-Modules in [27, \S 2].

\medskip

(iii) Let
$ \overline{X} = X\otimes_{k}\overline{k} $ with
$ \overline{k} $ the algebraic closure of
$ k $ in
$ \bC $.
For a prime number
$ l $, let
$ \Perv_{G}(\overline{X}, \bQ_{l}) $ be the category of \'etale perverse
sheaves on
$ \overline{X} $ with
$ \bQ_{l} $-coefficients [3] endowed with action of the Galois group
$ G = \Gal(\overline{k}/k) $ compatible with its action on
$ \overline{X} $ (i.e., for
$ \cF\in \Perv_{G}(\overline{X}, \bQ_{l}) $ and
$ g \in G $, there is an isomorphism
$ \cF\simto (\alpha_{g})^{*}\cF $ in
$ \Perv(\overline{X}, \bQ_{l}) $ compatible with multiplication of
$ G $ (cf. [13]), where
$ \alpha_{g} $ denotes the action of
$ g $ on
$ \overline{X}) $.
Since we will restrict to the full subcategory of the objects of geometric
origin, it would not be necessary to assume further condition on the action
in this stage.
We have a functor
$$
\Perv_{G}(\overline{X}, \bQ_{l})\to \Perv (X(\bC), \bQ_{l})
\leqno(1.8.2)
$$
by [3].
This functor would be compatible with direct image, dual and external
product by [loc. cit.].
Then, for
$ A = \bQ $, we can define
$ \cM(X) $ using the fiber product of
$ \MHM(X_{\bC},A) $ and
$ \Perv_{G}(\overline{X}, \bQ_{l}) $ over
$ \Perv(X(\bC), \bQ_{l}) $, where we assume that the weight filtration
$ W $ of the mixed Hodge Module
$ M $ is also defined on the corresponding etale perverse sheaf,
and a polarization of the Hodge Module
$ \Gr_{i}^{W}M $ is defined also on the corresponding etale perverse
sheaf
$ \Gr_{i}^{W}K $ as the self duality
$ \bD(\Gr_{i}^{W}K)\simeq (\Gr_{i}^{W})(i) $.
Then the conditions of mixed sheaves should be satisfied.
For example, (1.6.8) for
$ \Perv_{G}(\overline{X}, \bQ_{l}) $ is reduced to
$$
\Hom(A, \bR\cH{om}(B,C)) = \Hom(A\otimes^{L}B, C)
\leqno(1.8.3)
$$
which would be compatible with (1.8.2).

\medskip

(iv) We can also combine (ii) and (iii), i.e., use the fiber product of
$ \MHM(X_{\bC},A) $,
$ M_{rh}(\cD_{X}) $ and
$ \Perv_{G}(\overline{X}, \bQ_{l}) $, where
$ A = \bQ $.
An object consists of
$ ((M, F), K, K_{l}; W) $ where
$ (M,F) \in \MF_{rh}(\cD_{X}), K \in \Perv(X(\bC), \bQ), K_{l}\in \Perv_{G}
(\overline{X}, \bQ_{l}) $, and
$ W $ is a finite filtration on
$ M, K, K_{l} $.
Moreover, they have comparison isomorphisms compatible with
$ W $ as above, such that
$ ((M\otimes_{k}\bC, F), K; W) $ is a mixed Hodge Module on
$ X_{\bC} $ and the graduation of
$ W $ has a self-duality inducing a polarization of Hodge Module as above.

\medskip

(v) Since the above example depends on the choice of the prime
$ l $ and the embedding
$ k \to \bC $, it would be more natural to consider
all the prime numbers
$ l $ and all the embeddings
$ \sigma : k \to \bC, \overline{\sigma} : \overline{k}
\to \bC $ (extending
$ \sigma $) where
$ \overline{k} $ is an algebraic closure of
$ k $.
This means that an object consists of
$ ((M, F), K_{\sigma}, K_{l}; W) $ where
$ M \in M_{rh}(\cD_{X}), F $ is a good
filtration on
$ M, K_{\sigma} \in \Perv (X_{\sigma}(\bC), \bQ) $ with
$ X_{\sigma} = X\otimes_{k,\sigma}\bC $ for each embedding
$ \sigma : k \to \bC,
K_{l}\in \Perv_{G}(\overline{X}, \bQ_{l}) $ for each prime number
$ l $,and
$ W $ is a finite filtration on
$
M, K_{\sigma}, K_{l} $, such that they have comparison isomorphisms
for any
$ \sigma $ and
$ \overline{\sigma}, l $ (see also [11] [12] [14] [20]), and
$ ((M\otimes_{k,\sigma}\bC,
F), K_{\sigma}; W) $ is a mixed Hodge Module on
$ X_{\sigma} $ for any
$ \sigma . $ Moreover, the graduation of
$ W $ has a self-duality inducing a polarization of Hodge Module.
In this way, we would get a natural extension of [11] [12] [14] [20].

\medskip\noindent
{\it Remark.} Let
$ X $ be a smooth variety over a subfield
$ k $ of
$ \bC $, and
$ X_{\bC} = X\otimes_{k}\bC $.
We have the ring of differential operators
$ \cD_{X} $ on
$ X $ by Grothendieck [16] (cf. also [8] [34, (1.20.2)]).
A
$ \cD_{X} $-Module
$ M $ is called holonomic or regular holonomic if so is
$ M\otimes_{k}\bC $,
where
$ M\otimes_{k}\bC $ is defined by the scalar extension of the global
sections on each affine open subset.
Then a holonomic
$ \cD_{X} $-Module
$ M $ has the filtration
$ V $ of Kashiwara [22] and Malgrange [24] on its direct image by the
embedding by graph of a function on
$ X $, because it admits locally the
$ b $-function (i.e., Bernstein polynomial) as in [36] (this can be
checked also using the Weyl algebra argument as in [8]).
Let
$ (M,F) $ be a regular holonomic
$ \cD_{X} $-Module with a good filtration, and assume
$ (M,F)\otimes_{k}\bC $ underlies a Hodge Module.
Since the scalar extension by
$ k \to \bC $ is exact and faithful, the intersection of
$ F $ and
$ V $ commutes with the scalar extension, and the condition of Hodge
Module on
$ F $ and
$ V $ is satisfied over
$ k $.
By the same argument as in [28, (2.8)] we can define the direct image of
$ (M,F) $ by an open embedding whose complement is a locally principal
divisor, using the filtration
$ V $.
So the functors in (1.8) can be defined for the underlying filtered
$ \cD_{X} $-Modules
$ (M,F) $.

\medskip\noindent
{\bf 1.9.}~{\it Remark.} Let
$ K $ be a subfield of
$ \bC $ containing
$ k $ such that
$ K $ is finitely generated over
$ k $.
If a theory of mixed sheaves
$ \cM(X) $ over
$ \cV(k) $ is given, we can define a theory of mixed sheaves
$ \cM_{K}(X) $ over
$ \cV(K) $ as follows.
Let
$ S = \Spec R $ for a subring
$ R $ of
$ K $, which is finitely generated as
$ k $-algebra, and generates
$ K $.
For
$ X \in \cV(K) $, there is an open subvariety
$ U $ of
$ S $ and
$ X_{U} \in \cV(k) $ with a morphism
$ X_{U} \to U $ whose fiber at the generic point is
$ X $.
(Here
$ X_{U} $ is essentially unique by restricting to sufficiently small open
subvarieties.)
We define
$ \cM_{K}(X) $ by the inductive limit of
$ \cM(X_{U}) $ where
$ U $ runs over open subvarieties of
$ S $.
The inclusion
$ K \to \bC $ determines a closed point
$ x $ of
$ S_{\bC} $ which is not contained in any divisor defined over
$ k $, and
$ \For $ is defined by restricting the underlying
$ A $-perverse sheaf to the fiber at
$ x $.

\newpage
\centerline{\bf \S 2. Direct Images}

\bigskip
\noindent
{\bf 2.1.}
Let
$ f : X \to Y $ be a morphism of varieties of
$ \cV $.
Define a cohomological functor
$ H^{k}f_{!} $ by
$$
H^{k}f_{!} = \bD\scirc H^{-k}f_{*}\scirc \bD
\leqno(2.1.1)
$$
so that
$$
\For\scirc H^{k}f_{!} ={}^{p}\cH^{k}f_{!}\scirc \For ,
\leqno(2.1.2)
$$
$$
W_{i+k}(H^{k}f_{!}M) = H^{k}f_{!}M\quad \text{if}\quad W_{i}M = M
\leqno(2.1.3)
$$
for
$ M \in \cM(X) $.
If
$ f $ is a closed embedding or an affine open embedding, we define
$ f_{!} $ by
$ H^{0}f_{!} $ as in (1.4) so that
$$
f_{!}\scirc \bD = \bD\scirc f_{*}.
\leqno(2.1.4)
$$
For a closed embedding
$ i $, we have by (1.4.5):
$$
i_{!} = i_{*}.
\leqno(2.1.5)
$$

\medskip\noindent{{\bf 2.2.~Lemma.}} {\it
Let
$ \{U_{i}\} $ be a finite affine open covering of
$ X $, and
$ U_{I} = \bigcap_{i\in I} U_{i} $ with natural inclusion
$ j_{I} : U_{I} \to X $.
Then, for
$ M \in \cM(X) $, there is a natural quasi-isomorphism of
$ M $ to the
Cech complex whose
$ p^{th} $ component is
$$
\mopls_{|I|=p+1} (j_{I})_{*}(j_{I})^{*}M\,\,\, \text{for }p \ge 0\,\,\,
\text{and }0\,\,\, \text{otherwise } .
\leqno(2.2.1)
$$
}

\medskip\noindent
{\it Proof.} By (1.4.7), this is reduced to the isomorphism
$$
(j_{1})_{*}(j_{1})^{*}(j_{2})_{*}(j_{2})^{*} = (j_{\{1,2\}})_{*}
(j_{\{1,2\}})^{*}
\leqno(2.2.2)
$$
which follows from (1.3.6) (1.4.4) (1.4.8).

\medskip\noindent
{\it Remark.} If
$ X \in V $ is a disjoint union of open subvarieties
$ U_{1}, U_{2} $ with natural inclusions
$ j_{1} : U_{1} \to X, j_{2} : U_{2} \to X $.
Then we have an equivalence of categories
$$
((j_{1})^{*}, (j_{2})^{*}) : \cM(X) \simto \cM(U_{1})
\times \cM(U_{2})
\leqno(2.2.3)
$$
with quasi-inverse defined by
$ (M_{1},M_{2}) \to (j_{1})_{*}M_{1} \oplus (j_{2})_{*}M_{2} $ as a
special case of the above Lemma.

\medskip\noindent{{\bf 2.3.~Corollary}} {\it
The objects and morphisms of
$ \cM(X) $ are defined locally.
}

\medskip\noindent{{\bf 2.4.~Proposition.}} {\it
Let
$ f : X \to Y $ be a morphism of varieties of
$ \cV $.
Then we have exact functors of triangulated categories:
$$
f_{*}, f_{!} : D^{b}\cM(X) \to D^{b}\cM(Y)
\leqno(2.4.1)
$$
such that
$$
f_{!}\scirc \bD = \bD\scirc f_{*},
\leqno(2.4.2)
$$
and they are compatible with the corresponding functors and isomorphism
on the underlying
$ A $-complexes by the functor
$ \For $ in {\rm (1.1.5)}.
If
$ f $ is affine, the induced cohomological functors coincide with
{\rm (1.4.1) (2.1.1)} \(i.e., the isomorphism of the underlying perverse
sheaves is lifted to an isomorphism of
$ \cM(Y)) $.
}

\medskip\noindent{\it Proof.}
By definition, it is enough to show the assertion for
$ f_{*} $.
It is checked using an affine covering
$ \{U_{i}\} $ of
$ X $ and a sheaf theoretic resolution in [2] (because
$ |k| = \infty ) $.
In fact, for a finite number of objects
$ M_{k} \in \cM(X) $, there are surjective morphisms
$ M'_{k} \to M_{k} $ compatible with any morphisms between
$ M_{k} $ (i.e., the morphisms are naturally lifted to morphisms between
$ M'_{k}) $ such that
$ H^{k}(f_{I})_{*}(j_{I})^{*}M'_{k} = 0 $
$ (k \ne 0) $ for any
$ f $ (i.e.,
$ M'_{k} $ can be taken independently of
$ f) $, where
$ f_{I} = fj_{I} $ with
$ j_{I} $ as in (2.2).
(In fact,
$ M'_{k} $ is a finite direct sum of
$ (j_{U})_{!}(j_{U})^{!}M $ for affine open subvarieties
$ U $ with
$ j_{U} : U \to X $, where
$ U $ are the complements of generic hyperplane sections of affine open
subvarieties which cover
$ X $.
See [loc. cit.].)
So, for a bounded complex
$ M $, we have a resolution
$ M' \to M $ such that
$ M' $ is bounded above and
$$
H^{k}(f_{I})_{*}(j_{I})^{*}M'^{j} = 0\quad \text{for}\quad k \ne 0.
\leqno(2.4.3)
$$
The direct image
$ f_{*}M $ is defined by the
Cech double complex whose
$ (p,q) $-component is
$$
\mopls_{|I|=p+1} H^{0}(f_{I})_{*}(j_{I})^{*}M'^{q}\,\,\, \text{for }
p \ge 0\,\,\, \text{and }0\,\,\, \text{otherwise } .
\leqno(2.4.4)
$$
Here we use also the truncation
$ \tau $ in [35], because
$ H^{k}f_{*}M = 0 $ for
$ k << 0 $.
We can check the independentness of the choice of resolution and open
covering.

For the last assertion, we may assume
$ X, Y $ affine, because the assertion is local.
Then we can use the trivial covering of
$ X $ for the definition of direct image, and the assertion follows, because
the fact that
$ H^{k}f_{*} = 0 $ for
$ i > 0 $ and
$ H^{k}f_{*} $ is the cohomological left derived functor of
$ H^{0}f_{*} $.

\medskip\noindent
{\it Remark.} If
$ f $ is an open embedding, it is not necessary to take a resolution of
$ M $, because the direct image by an affine open embedding is exact.

\medskip\noindent{{\bf 2.5.~Corollary}} {\it
Let
$ U $ be an affine open subvariety of
$ X $, and
$ M \in \cM(U) $ such that
$ H^{k}(fj)_{*}M = 0 $ (resp.
$ H^{k}(fj)_{!}M = 0) $ for
$ k \ne 0 $, where
$ j : U \to X $.
Then we have a natural isomorphism in
$ D^{b}\cM(Y) : $
}
$$
f_{*}(j_{*}M) = H^{0}(fj)_{*}M\,\,\, \text{(resp. } f_{!}(j_{!}M) =
H^{0}(fj)_{!}M).
\leqno(2.5.1)
$$

\medskip\noindent
{\it Proof.} This follows from the proof of (2.4).
It is enough to show the assertion for
$ f_{*} $.
Take an affine open covering of
$ X $ such that
$ U $ is a member of the covering.
Then we get a morphism
$ f_{*}(j_{*}M) \to (fj)_{*}M $ by definition of direct image, because
$ j^{*}(j_{*}M) = M $.
Then the assertion follows from the last assertion of (2.4).

\medskip\noindent{{\bf 2.6.~Corollary}} {\it
For morphisms
$ f : X \to Y, g : Y \to Z $ of varieties of
$ \cV $, we have functorial isomorphisms
$$
g_{*}f_{*} = (gf)_{*},\,\,\,g_{!}f_{!} = (gf)_{!}
\leqno(2.6.1)
$$
which correspond to the natural isomorphisms on the underlying
$ A $-complexes by the functor
$ \For $ in {\rm (1.1.5)}.
}

\medskip\noindent{\it Proof.}
Take affine open coverings
$ \{U_{i}\}, \{V_{i}\} $ of
$ X, Y $ such that
$ f(U_{i}) \subset V_{i} $.
Then the assertion is reduced to (1.4.8) using (2.5), because the
resolution
$ M' $ in (2.4) is independent of
$ f $.

\medskip\noindent{{\bf 2.7.~Proposition.}} {\it
Let
$ f : X \to Y $ be a morphism of varieties of
$ \cV $, and
$ Z \in \cV $.
Then we have a functorial isomorphism
$$
(f\times id)_{*}(M\,\boxtimes\,N) = f_{*}M\,\boxtimes\,N,\,\,\,
(f\times id)_{!}(M\,\boxtimes\,N) = f_{!}M\,\boxtimes\,N
\leqno(2.7.1)
$$
for
$ M \in D^{b}\cM(X), N \in D^{b}\cM(Z) $.
}

\medskip\noindent{\it Proof.}
This follows from the definition of direct images and (1.5.8).

\medskip\noindent
{\it Remarks.} (i) By (1.5.2), the external product is exact for both
factors, and it is naturally extended to
$$
\boxtimes : D^{b}\cM(X) \times D^{b}\cM(Y) \to
D^{b}\cM(X\times Y).
\leqno(2.7.2)
$$
We have a natural isomorphism
$$
H^{k}(M\,\boxtimes\,N) = \mopls_{i+j=k} H^{i}M\,\boxtimes\,H^{j}N
\leqno(2.7.3)
$$
induced by the inclusion
$ \Ker\,d\,\boxtimes\,\Ker\,d \to \Ker(d\,\boxtimes\,id \pm
id\,\boxtimes\,d) $.

(ii) We can deduce
$$
(id\times f)_{*}(N\,\boxtimes\,M) = N\,\boxtimes\,f_{*}M,\,\,\,
(id\times f)_{!}(N\,\boxtimes\,M) = N\,\boxtimes\,f_{!}M
\leqno(2.7.4)
$$
from (2.7.1), using (1.5.5).
Combined with (2.6), we have
$$
(f\times g)_{*}(M\,\boxtimes\,N) = f_{*}M\,\boxtimes\,g_{*}N,\,\,\,
(f\times g)_{!}(M\,\boxtimes\,N) = f_{!}M\,\boxtimes\,g_{!}N
\leqno(2.7.5)
$$
for
$ f : X \to Y, g : Z \to W $ and
$ M \in D^{b}\cM(X), N \in D^{b}\cM(Z) $, using the decomposition
$ f\times g = (f\times id)\scirc (id\times g) $.

\bigskip\bigskip
\centerline{\bf \S 3. Pull-Backs}

\bigskip
\noindent
{\bf 3.1.}
Let
$ f : X \to Y $ be a morphism of varieties of
$ \cV $.
We define functors
$$
f^{*}, f^{!} : D^{b}\cM(Y) \to D^{b}\cM(X)
\leqno(3.1.1)
$$
by the left and right adjoint functors of
$ f_{*}, f_{!} $ respectively, i.e., we have functorial isomorphisms
(called the adjoint relations):
$$
\Hom(f^{*}N,M) = \Hom(N,f_{*}M),\,\,\,\Hom(M,f^{!}N) = \Hom(f_{!}M,N)
\leqno(3.1.2)
$$
for
$ M \in D^{b}\cM(X), N \in D^{b}\cM(Y) $.
So these functors are unique if they exist.
By (2.4.2) and (2.6), we have
$$
f^{!}\scirc \bD = \bD\scirc f^{*}
\leqno(3.1.3)
$$
$$
f^{*}g^{*} = (gf)^{*},\,\,\,f^{!}g^{!} = (gf)^{!}
\leqno(3.1.4)
$$
for
$ f : X \to Y, g : Y \to Z $.

The first isomorphism of (3.1.2) is equivalent to two morphisms
$ \alpha : N \to f_{*}f^{*}N $ and
$ \beta : f^{*}f_{*}M \to M $ such that
$$
\beta\scirc f^{*}\alpha : f^{*}N \to f^{*}f_{*}f^{*}N \to
f^{*}N
\leqno(3.1.5)
$$
$$
f_{*}\beta\scirc \alpha : f_{*}M \to f_{*}f^{*}f_{*}M \to
f_{*}M
\leqno(3.1.6)
$$
are isomorphisms and one of them is the identity (which implies that the
two morphisms are inverse of each other).
This is checked using a commutative diagram
$$
\CD
f^{*}N @>>> f^{*}f_{*}f^{*}N @>>> f^{*}N
\\
@. @VVV @VVV
\\
@. f^{*}f_{*}M @>>> M
\endCD
\leqno(3.1.7)
$$
and a similar diagram.
The argument is similar for the second isomorphism of (3.1.2).

The morphisms
$$
N \to f_{*}f^{*}N\quad \text{and}\quad f_{!}f^{!}M \to M
\leqno(3.1.8)
$$
associated with (3.1.2) are called the {\it restriction} and {\it Gysin}
morphisms respectively.

We will show the existence of
$ f^{*}, f^{!} $ by factorizing
$ f $ into a composition of closed embedding and projection.
See (3.7) and (3.8--9).
The arguments show also that they are exact functors of triangulated
categories,
which are compatible with the corresponding functors on the underlying
$ A $-complexes by
$ \For $, and the adjoint relation (3.1.2) is also compatible with
$ \For $.
Moreover, we will show the commutativity with the external products, e.g.,
$$
(f\times id)^{*}(N\,\boxtimes\,M) = f^{*}N\,\boxtimes\,M,\,\,\,
(f\times id)^{!}(N\,\boxtimes\,M) = f^{!}N\,\boxtimes\,M
\leqno(3.1.9)
$$
for
$ M \in D^{b}\cM(Z), N \in D^{b}\cM(Y) $ (similarly for
$ id\times f) $.

We first check the compatibility with the previous definition of open
pull-backs:

\medskip\noindent{{\bf 3.2.~Proposition.}} {\it
Let
$ j : X \to Y $ be an open embedding of varieties of
$ \cV $.
Then we have canonical functorial morphisms and isomorphisms:
$$
j_{!}j^{!} \to id \to j_{*}j^{*}\,\,\, \text{and }
j^{*}j_{*} = id = j^{!}j_{!},
\leqno(3.2.1)
$$
which correspond to the natural morphisms and isomorphisms on the underlying
$ A $-complexes by the functor
$ \For $ in {\rm (1.1.5)}.
Moreover, {\rm (3.2.1)} induces the adjoint relations
$$
\Hom(j^{*}N,M) = \Hom(N,j_{*}M),\,\,\,\Hom(M,j^{!}N) = \Hom(j_{!}M,N)
\leqno(3.2.2)
$$
for
$ M \in D^{b}\cM(X), N \in D^{b}\cM(Y) $.
}

\medskip\noindent{\it Proof.}
The first assertion follows from the definition of direct image (cf. Remark
after (2.4)) and (1.4.7), using (2.2).
For (3.2.2), we show the first isomorphism, because the argument is similar
for the second.
The morphism of (3.2.2) is induced by
$ \alpha : N \to j_{*}j^{*}N $, and its inverse by
$ \beta : j^{*}j_{*}M \to M $ in (3.2.1).
We have to show (3.1.5--6).
In this case,
$ \beta $ is an isomorphism, and
$ \beta^{-1} : j^{*}N \to j^{*}j_{*}j^{*}N $ coincides with
$ j^{*}\alpha : j^{*}N \to j^{*}j_{*}j^{*}N $ by definition.
So (3.1.5) is the identity.
For (3.1.6), the assertion follows from (1.1.7).

\medskip\noindent{{\bf 3.3.~Proposition.}} {\it
Let
$ i : X \to Y $ be a closed embedding of varieties of
$ \cV(k) $ such that
$ Y \in \cV $.
Let
$ j : Y \setminus i(X) \to Y $ be the natural morphism.
Then we have exact functors of triangulated categories:
$$
i_{!}i^{!}, i_{*}i^{*} : D^{b}\cM(Y) \to {D}_{X}^{b}\cM(Y)
\leqno(3.3.1)
$$
and canonical distinguished triangles of functors
$$
\to i_{!}i^{!} \to id \to j_{*}j^{*} \to ,
\,\,\,\to j_{!}j^{!} \to id \to i_{*}i^{*} \to,
\leqno(3.3.2)
$$
which are compatible with the corresponding functors and triangles on the
underlying
$ A $-complexes by the functor
$ \For $ in {\rm (1.1.5)}, where
$ {D}_{X}^{b}\cM(Y) $ is the full subcategory of
$ D^{b}\cM(Y) $ consisting of the objects
$ M $ such that
$ \supp H^{k}M \subset i(X) $ \(cf. {\rm (1.4.6))}.
Furthermore, these functors commute with external products as in
{\rm (3.1.9)}.
}

\medskip\noindent{\it Proof.}
It is enough to show the assertion for
$ i_{!}i^{!} $ by duality.
Let
$ \{U_{i}\} $ be a finite affine open covering of
$ Y \setminus i(X) $.
Then, for a complex
$ M $,
$ i_{!}i^{!}M $ is defined by a
Cech double complex whose
$ (p,q) $-component is
$$
\mopls_{|I|=p} (j_{I})_{*}(j_{I})^{*}M^{q}\,\,\, \text{for }p \ge
0\,\,\, \text{and }0\,\,\, \text{otherwise } ,
\leqno(3.3.3)
$$
and
$ i_{!}i^{!}M \to M $ is defined by the natural projection, where
$ j_{I} $ is as in (2.2) and
$ j_{\emptyset} = id $.
Moreover, the mapping cone of
$ i_{!}i^{!}M \to M $ is naturally quasi-isomorphic to
$ j_{*}j^{*}M $ by Remark after (2.4).
We can check the independentness of the covering, using a covering whose
members are the union of the member of two coverings, so that a morphism
of complexes in the derived category is constructed.
Here it is enough to construct a morphism, because it is a
quasi-isomorphism by (1.1.1) and (1.1.6).
The last assertion follows from (2.7.1).

\medskip\noindent{{\bf 3.4.~Corollary}} {\it
With the above notation, assume
$ X \in \cV $.
Then we have canonical cohomological functors
$$
H^{k}i^{!}, H^{k}i^{*} : D^{b}\cM(Y) \to \cM(X)
\leqno(3.4.1)
$$
and canonical long exact sequences of cohomological functors
$$
\aligned
&\to i_{!}H^{k}i^{!} \to H^{k} \to H^{k}j_{*}j^{*} \to
i_{!}H^{k+1}i^{!} \to ,
\\
&\to H^{k}j_{!}j^{!} \to H^{k} \to i_{*}H^{k}i^{*} \to
H^{k+1}j_{!}j^{!} \to ,
\endaligned
\leqno(3.4.2)
$$
which are compatible with the corresponding cohomological functors
$^{p}\cH^{k}i^{!},^{p}\cH^{k}i^{*} $ and long exact sequences on
the underlying perverse sheaves by the functor
$ \For $.
Furthermore, these functors commute with
$ \boxtimes\,M $ for
$ M \in \cM(Z) $ as in {\rm (3.1.9)}.
}

\medskip\noindent{\it Proof.}
This follows from (3.3) and (1.4.6), because
$ i_{!} = i_{*} $ and
$ j^{*} = j^{!} $ are exact functors and commute with external products,
and
$ \boxtimes\,M $ is an exact functor.

\medskip\noindent{{\bf 3.5.~Proposition.}} {\it
With the notation of {\rm (3.3)}, the morphisms
$ i_{!}i^{!} \to id $ and
$ id \to i_{*}i^{*} $ in {\rm (3.3.2)} induce the adjoint relations
$$
\Hom(M,i_{!}i^{!}N) = \Hom(M,N),\,\,\,\Hom(i_{*}i^{*}N,M) = \Hom(N,M)
\leqno(3.5.1)
$$
for
$ M \in {D}_{X}^{b}\cM(Y), N \in D^{b}\cM(Y) $, where
$ {D}_{X}^{b}\cM(Y) $ is as in {\rm (3.3)}.
}

\medskip\noindent{\it Proof.}
This follows from the triangles in (3.3.2) and (3.2), because
$ j^{*}M = 0 $ (which is reduced to the case
$ M \in \cM_{X}(Y)) $.

\medskip\noindent
{\it Remark.} By (3.5), the existence of the pull-back functors in the
closed embedding case (cf. (3.7)) is reduced to the equivalence of
categories:
$$
i_{!} = i_{*} : D^{b}\cM(X) \simto {D}_{X}^{b}\cM(Y),
\leqno(3.5.2)
$$
which will be proved in (5.6).
So we will consider the projection case.

\medskip\noindent
{\bf 3.6.}
Let
$ X, Y \in \cV $, and
$ p : X\times Y \to Y $ be the second projection.
Assume there exists
$ A_{X}^{\cM} \in D^{b}\cM(X) $ with an isomorphism
$ \For(A_{X}^{\cM}) = A_{X(\bC)} $ and a morphism
$$
A^{\cM} \to H^{0}(a_{X})_{*}A_{X}^{\cM}
\leqno(3.6.1)
$$
whose image by
$ \For $ is the natural morphism
$ A \to H^{0}(a_{X(\bC)})_{*}A_{X(\bC)} $.
(If
$ X $ is smooth,
$ A_{X}^{\cM} $ exists by (1.6).)
Note that (3.6.1) induces naturally a morphism in
$ D^{b}\cM(\Spec k) $:
$$
A^{\cM} \to (a_{X})_{*}A_{X}^{\cM},
\leqno(3.6.2)
$$
because
$ H^{0}(a_{X})_{*}A_{X}^{\cM} = 0 $ for
$ k < 0 $.
Let
$$
p^{*}N = A_{X}^{\cM}\,\boxtimes\,N,\,\,\,p^{!}N = \bD A_{X}^{\cM}
\,\boxtimes\,N,
\leqno(3.6.3)
$$
and define the morphism
$ \alpha : N \to p_{*}p^{*}N $,
$ \alpha : p_{!}p^{!}N \to N $ in (3.1) by
$$
N = A^{\cM}\,\boxtimes\,N \to (a_{X})_{*}A_{X}^{\cM}\,\boxtimes\,N =
p_{*}p^{*}N
\leqno(3.6.4)
$$
and its dual, where the first and last isomorphisms are induced by (1.6.4)
combined with (1.5.5), and (2.7.1) respectively, and the middle morphism
by (3.6.2).
The construction of
$ \beta $ will be given in (3.8).
We will use the fact that the composition
$$
N \to p_{*}p^{*}N \to p_{*}(i_{*}i^{*})p^{*}N
\leqno(3.6.5)
$$
is an isomorphism for a section
$ i : Y \to X\times Y $ of
$ p $, where the morphisms are induced by (3.6.4) and (3.3.2).
This follows from (1.1.7), because the underlying morphism of
$ A $-complexes is an isomorphism.

\medskip

In the later part of this section we will use (3.5.2).
Note that the proof of (3.5.2) in \S 5 does not use the later part of this
section and \S 4.
We first get the following Theorem by (3.5) and (3.5.2).

\medskip\noindent{{\bf 3.7.~Theorem.}} {\it
Let
$ i, j $ be as in {\rm (3.3)}.
Then we have exact functors of triangulated categories:
$$
i^{!}, i^{*} : D^{b}\cM(Y) \to D^{b}\cM(X)
\leqno(3.7.1)
$$
compatible with the corresponding functors on the underlying
$ A $-complexes by the functor
$ \For $ in {\rm (1.1.5)}, so that their compositions with
$ i_{!} = i_{*} $ coincides with {\rm (3.3.1)}, and
$ i^{!}, i^{*} $ commute with external product as in {\rm (3.1.9)}.
Moreover, the morphisms
$ i_{!}i^{!} \to id $ and
$ id \to i_{*}i^{*} $ in {\rm (3.3.2)} induce the adjoint relations
$$
\Hom(M,i^{!}N) = \Hom(i_{!}M,N),\,\,\,\Hom(i^{*}N,M) = \Hom(N,i_{*}M)
\leqno(3.7.2)
$$
for
$ M \in D^{b}\cM(X), N \in D^{b}\cM(Y) $.
}

\medskip\noindent
{\it Remark.} We have a functorial isomorphism
$$
i^{!}i_{!} = i^{*}i_{*} = id,
\leqno(3.7.3)
$$
because
$ i_{!}i^{!}M = i_{*}i^{*}M = M $ for
$ M \in {D}_{X}^{b}\cM(Y) $ by definition.

\medskip\noindent{{\bf 3.8.~Theorem.}} {\it
Let
$ p : X\times Y \to Y $ be a projection.
Assume
$ A_{X}^{\cM} \in D^{b}\cM(X) $ in {\rm (3.6)} exists.
Then the pull-back functors
$ p^{*}, p^{!} $ are defined by {\rm (3.6.3)}, so that they commute
with external product as in {\rm (3.1.9)}.
Furthermore,
$ A_{X}^{\cM} $ in {\rm (3.6)} is unique and we have a canonical
isomorphism
$ A_{X}^{\cM} = {a}_{X}^{*}A^{\cM} $ for
$ a_{X} : X \to \Spec k $.
}

\medskip\noindent{\it Proof.}
We show the assertion for
$ p^{*} $, because the dual argument holds for
$ p^{!} $.
The morphism
$ \alpha : N \to f_{*}f^{*}N $ in (3.1) is defined in (3.6.4).
For the construction of
$ \beta $, we use a commutative diagram due to Kashiwara (see also
[28] [33]):
$$
\CD
X \times Y @<{{q}_{2}}<< X \times X \times Y Y @<{i}<< X
\times Y
\\
@V{p_{1}}VV @V{q_{1}}VV
\\
Y Y @<{{p}_{2}}<< X \times Y
\endCD
\leqno(3.8.1)
$$
where
$ q_{1}, q_{2} $ are induced by the projections of
$ X\times X $ to the first and second factors,
$ f = p_{1} = p_{2} $ and
$ q_{1}i = q_{2}i = id $.
We define
$ \beta : f^{*}f_{*}M \to M $ by
$$
{p}_{2}^{*}(p_{1})_{*}M = (q_{1})_{*}{q}_{2}^{*}M \to
(q_{1})_{*}i_{*}i^{*}{q}_{2}^{*}M = i^{*}{q}_{2}^{*}M = M
\leqno(3.8.2)
$$
where the first three morphisms are induced by (2.7.1), (3.3.2), (2.6.1),
and the last by (2.6.1) and (3.6.5) (applied to
$ q_{2} $ and
$ i) $.
Then
$ (p_{1})_{*} $ of the first morphism of (3.6.5) (applied to
$ q_{2} $ and
$ i) $ coincides with
$$
\alpha : (p_{1})_{*}M \to (p_{2})_{*}{p}_{2}^{*}((p_{1})_{*}M)
$$
by (2.6.1), (2.7.1).
So (3.1.6) is the identity.
Since (3.1.5) is an isomorphism by (1.1.7), we get the first assertion.
The last assertion follows from (1.6.4) and the uniqueness of the adjoint
functor for
$ p = a_{X} $.

\medskip\noindent{{\bf 3.9.~Proposition.}} {\it
In the notation of {\rm (3.8)},
$ A_{X}^{\cM} $ \(and hence
$ p_{*}, p_{!}) $ exists.
}

\medskip\noindent{\it Proof.}
It is enough to show the existence of
$ A_{X}^{\cM} $.
The case
$ X $ smooth is clear by (1.6).
If there is a closed embedding
$ i : X \to Y $ with
$ Y $ smooth (e.g.,
$ Y = \bA^{n}) $, we have
$ i^{*}A_{Y}^{\cM} \in D^{b}\cM(X) $ using (5.6), and it satisfies the
condition of
$ A_{X}^{\cM} $ by an argument similar to (3.10).
Moreover, it is unique by (3.8).
In general we may assume that
$ X $ is covered by two open subvarieties
$ U_{1}, U_{2} $ such that
$ A_{{U}_{i}}^{\cM} $ exists for
$ i = 1, 2 $.
Let
$ U_{3} = U_{1} \cap U_{2} $ with natural morphisms
$ j_{i} : U_{i} \to X, j'_{i} : U_{3} \to U_{i} $.
Then we have a canonical isomorphism
$ (j'_{i})^{*}A_{{U}_{i}}^{\cM} = A_{{U}_{3}}^{\cM} $ with a morphism
$ A_{{U}_{i}}^{\cM} \to (j'_{i})_{*}A_{{U}_{3}}^{\cM} $ for
$ i = 1, 2 $ by uniqueness of
$ A_{{U}_{3}}^{\cM} $.
So
$ A_{X}^{\cM} $ is defined by the mapping cone:
$$
C(\mopls_{i=1,2} (j_{i})_{*}A_{{U}_{i}}^{\cM} \to (j_{3})_{*}
A_{{U}_{3}}^{\cM})[-1]
$$
whose underlying
$ A $-complex is
$$
C(\mopls_{i=1,2} (j_{i})_{*}A_{U_{i^{(\bC)}}} \to
(j_{3})_{*}A_{U_{3^{(\bC)}}})[-1] = A_{X(\bC)}.
$$

\medskip\noindent
{\it Remark.} For a morphism
$ f : X \to Y $ of varieties of
$ \cV $, we have
$$
f^{*}A_{Y}^{\cM} = A_{X}^{\cM}
\leqno(3.9.1)
$$
by (3.8) and (3.1.4).
So (3.1.8) induces the morphisms
$$
A_{Y}^{\cM} \to f_{*}A_{X}^{\cM},\,\,\,f_{!}\bD A_{X}^{\cM} \to \bD
A_{Y}^{\cM}
\leqno(3.9.2)
$$
which are called the {\it restriction} and {\it Gysin} (or {\it trace})
morphisms respectively.
They are compatible with the composition of morphisms of varieties,
because the adjoint isomorphism (3.1.2) is compatible with the composition
by the proof of (3.1.4).
So they are compatible with the restriction and Gysin morphisms on the
$ A $-complexes by the functor
$ \For $, because it is clear by definition in the case of closed embedding
and projection.

\medskip\noindent{{\bf 3.10.~Proposition.}} {\it
For a cartesian diagram in
$ \cV(k) : $
$$
\CD
X @<{g'}<< X'
\\
@V{f}VV @V{f'}VV
\\
Y @<{g}<< Y'
\endCD
\leqno(3.10.1)
$$
such that
$ X, Y, X', Y' \in \cV $, we have canonical functorial isomorphisms
}
$$
g^{!}f_{*} = f'_{*}g'^{!},\,\,\,g^{*}f_{!} = f'_{!}g'^{*}.
\leqno(3.10.2)
$$

\medskip\noindent
{\it Proof.} If
$ g $ is a projection, the assertion follows from (2.7).
If
$ g $ is a closed embedding, it follows from the definition of direct image
and pull-backs, using open coverings
$ \{U_{i}\}_{i\in I}, \{V_{i}\}_{i\in I} $ of
$ X, Y $ such that
$ \{U_{i}\}_{i\in J}, \{V_{i}\}_{i\in J} $ are coverings of
$ X \setminus X', Y \setminus Y' $ respectively.
See [28, (4.4.3)] for detail.

\bigskip\bigskip
\centerline{\bf \S 4. Duality}

\bigskip
\noindent
{\bf 4.1.~Definition.} Let
$ X \in \cV $, and
$ \delta : X \to X\times X $ the diagonal morphism.
For
$ M, N \in D^{b}\cM(X) $, let
$$
N\otimes M = \delta^{*}(N\,\boxtimes\,M),\,\,\,\cH{om}
(N,M) = \delta^{!}(\bD N\,\boxtimes\,M).
\leqno(4.1.1)
$$
Then
$$
\For(N\otimes M) = \For(N)\otimes\For(M),\,\,\,
\For(\cH{om}(N,M)) = \cH{om} (\For(N),\For(M)),
\leqno(4.1.2)
$$
cf. [8] [32].
By (3.1.9) (3.1.4) (1.5.4) (1.5.5), we have
$$
(N\otimes M)\otimes L = N\otimes (M\otimes L),
\leqno(4.1.3)
$$
$$
N\otimes M = M\otimes N,
\leqno(4.1.4)
$$
because
$ \iota^{*} = (\iota^{-1})_{*} $ by
$ (\iota^{-1})_{*}\iota_{*} = id, \iota_{*}(\iota^{-1})_{*} = id $.
We have also
$$
A_{X}^{\cM}\otimes M = M\,\,\, \text{for }M \in \cM(X)
\leqno(4.1.5)
$$
by (3.6.5).
The image of these isomorphisms by the functor
$ \For $ is the natural isomorphisms.

\medskip\noindent
{\bf 4.2.~Definition.} Let
$ X $
$ = \bP^{1} $.
Define
$$
A^{\cM}(-1) = H^{2}(a_{X})_{!},A_{X}^{\cM},\,\,\,A^{\cM}(1) =
\bD(A^{\cM}(-1)).
\leqno(4.2.1)
$$
and
$$
A^{\cM}(n) = A^{\cM}(n-1)\otimes A^{\cM}(1), A^{\cM}(-n) =
A^{\cM}(1-n)\otimes A^{\cM}(-1)
\leqno(4.2.2)
$$
by induction on
$ n > 0 $.
Here
$$
M\otimes N = M\,\boxtimes\,N\quad \text{for}\quad M, N \in \cM
(\Spec k)
\leqno(4.2.3)
$$
using the identification
$ \Spec k \times \Spec k = \Spec k $.
We have a canonical isomorphism
$$
\For(A^{\cM}(n)) = A(n)
\leqno(4.2.4)
$$
where the right hand side is defined by
$$
A(n) = A\otimes_{\bZ} \bZ(n)\quad \text{with}\quad
\bZ(n) = (2\pi i)^{n}\bZ \subset \bC,
\leqno(4.2.5)
$$
cf. also [27, 2.5.7].
We have canonical isomorphisms
$$
A^{\cM}(n)\otimes A^{\cM}(m) = A^{\cM}(n+m)
\leqno(4.2.6)
$$
for
$ m, n \in \bZ $, whose image by
$ \For $ is the natural isomorphisms, because the case
$ nm = -1 $ follows from (1.6.7).

We define the {\it Tate twist}
$ (n) $ by
$$
M(n) = M\,\boxtimes\,A^{\cM}(n)\quad \text{for}\quad M \in \cM
(X)
\leqno(4.2.7)
$$
using the identification
$ X \times \Spec k = X $.
It is an exact functor and is naturally extended to
$ D^{b}\cM(X) $.
We have
$$
(M(n))(m) = M(n+m)
\leqno(4.2.8)
$$
$$
\bD(M(n))= (\bD M)(-n)
\leqno(4.2.9)
$$
by (4.2.6) (1.5.6), and the Tate twist commutes with direct image and
pull-backs.

\medskip\noindent{{\bf 4.3.~Proposition.}} {\it
Let
$ X $ be a purely
$ n $-dimensional smooth variety of
$ \cV $.
Then we have a canonical isomorphism
$$
\bD A_{X}^{\cM} = A_{X}^{\cM}(n)[2n]
\leqno(4.3.1)
$$
whose image by the functor
$ \For $ is the canonical isomorphism.
In particular, we have the trace morphism
$$
H^{2n}(a_{X})_{!}A_{X}^{\cM}(n) \to A^{\cM}
\leqno(4.3.2)
$$
whose image by
$ \For $ is the natural trace morphism.
}

\medskip\noindent{\it Proof.}
By uniqueness of
$ A_{X}^{\cM} $ (cf. (3.8)), it is enough to show the last assertion,
because the dual of (4.3.2) is
$$
A^{\cM} \to H^{0}(a_{X})_{*}(\bD A_{X}^{\cM})(-n)[-2n].
$$
Here we may restrict
$ X $ to any dense open subvariety
$ U $ of
$ X $, because the inclusion
$ j : U \to X $ induces an isomorphism
$$
H^{2n}(a_{U})_{!}A_{U}^{\cM} \to H^{2n}(a_{X})_{!}A_{X}^{\cM}.
\leqno(4.3.3)
$$
If
$ X = \bA^{1}, (4.3.2) $ is clear by (4.2.1) and (4.3.3).
If
$ X = \bA^{n} $, we have a decomposition
$ X = X'\times X'' $ with
$ X' = \bA^{1}, X'' = \bA^{n-1} $, and
$$
(a_{X})_{!}A_{X}^{\cM} = (a_{X'})_{!}A_{X'}^{\cM}\,\boxtimes\,
(a_{X''})_{!}A_{X''}^{\cM}
\leqno(4.3.4)
$$
by (2.7.5), because
$ A_{X}^{\cM} = A_{X'}^{\cM}\,\boxtimes\,A_{X''}^{\cM} $ by (3.8).
So we get the assertion.

In general, we may assume that there is a proper \'etale morphism
$ f : X \to Y $ such that
$ Y $ is an open subvariety of
$ \bA^{n} $, by shrinking
$ X $.
By Lemma (4.4) below, we have a canonical isomorphism
$ f^{!}A_{Y}^{\cM} = f^{*}A_{Y}^{\cM} = A_{X}^{\cM} $ which induces
a morphism
$ f_{!}A_{X}^{\cM} \to A_{Y}^{\cM} $ by (3.1.8).
So we get a morphism
$$
H^{2n}(a_{X})_{!}A_{X}^{\cM} \to H^{2n}(a_{Y})_{!}A_{Y}^{\cM},
\leqno(4.3.5)
$$
and (4.3.2) is defined by composition with the trace morphism of
$ Y $.

\medskip\noindent{{\bf 4.4.~Lemma.}} {\it
Let
$ f : X \to Y $ be an \'etale morphism of smooth varieties of
$ \cV $.
Then we have a canonical isomorphism
$$
f^{!}M = f^{*}M\quad \text{for}\quad M \in \cM(Y)
\leqno(4.4.1)
$$
whose image by the functor
$ \For $ is the natural isomorphism.
}

\medskip\noindent{\it Proof.}
Let
$ X' = X\times_{Y}X $ with the natural morphisms
$ i : X' \to X\times X, f' : X' \to Y $.
Taking the pull-back of the pairing (1.6.7), we get
$$
f^{!}\bD M\,\boxtimes\,f^{!}M \to (f\times f)^{!}
\delta_{*}\bD A_{Y}^{\cM}
\leqno(4.4.2)
$$
by (3.1.9), where
$ \delta $ denotes the diagonal morphisms.
We have
$$
(f\times f)^{!}\delta_{*}\bD A_{Y}^{\cM} = i_{*}f'^{!}\bD A_{Y}^{\cM} =
i_{*}\bD
A_{Y'}^{\cM}
\leqno(4.4.3)
$$
by (3.10) (3.1.4).
Since
$ f $ is \'etale, the diagonal
$ X $ is an open and closed subvariety of
$ Y' $, and we have a natural morphism
$ i_{*}\bD A_{Y'}^{\cM} \to \delta_{*}\bD A_{X}^{\cM} $.
So we get a morphism
$$
f^{!}\bD M \to \bD f^{!}M
\leqno(4.4.4)
$$
by (1.6.8).
This is an isomorphism, because so is the underlying morphism of perverse
sheaves.
So we get the assertion by (3.1.3).

\medskip\noindent
{\bf 4.5.}
Let
$ X $ be a purely
$ n $-dimensional smooth variety of
$ \cV $.
We say that
$ L \in D^{b}\cM(X) $ is {\it smooth}, if
$ \For(L) $ is a local system (in particular,
$ L[n] \in \cM(X)) $.
For smooth
$ L \in D^{b}\cM(X) $, we define
$$
L^{*} = (\bD L)(-n)[-2n]
\leqno(4.5.1)
$$
so that
$$
\For(L^{*}) = \For(L)^{*} := \cH{om}_{A} (\For(L),A)
\leqno(4.5.2)
$$
$$
(A_{X}^{\cM})^{*} = A_{X}^{\cM}
\leqno(4.5.3)
$$
by (4.3.1).
In particular,
$ L^{*} $ is smooth.
Then (1.6.7) (4.3.1) induce a pairing
$$
L^{*}\otimes L \to A_{X}^{\cM},
\leqno(4.5.4)
$$
whose image by the functor
$ \For $ is the natural pairing
$$
\For(L)^{*}\otimes \For(L) \to A_{X(\bC)}.
\leqno(4.5.5)
$$
It induces further
$$
(L\otimes L^{*})\otimes (L^{*}\otimes L) \to A_{X}^{\cM},
\leqno(4.5.6)
$$
using (4.1.3--5), so that
$$
(L^{*}\otimes L)^{*} = L\otimes L^{*}
\leqno(4.5.7)
$$
Taking the dual of (4.5.4), we get
$$
A_{X}^{\cM} \to L\otimes L^{*}.
\leqno(4.5.8)
$$

\medskip\noindent{{\bf 4.6.~Proposition.}} {\it
Let
$ X $ be a smooth variety in
$ \cV $, and
$ L $ a smooth object of
$ \cM(X) $.
Then,
$ M\otimes L \in \cM(X) $ for
$ M \in \cM(X) $, and the functor
$ \otimes L $ is exact.
Moreover, we have a canonical isomorphism
$$
\Hom(N, M\otimes L^{*}) = \Hom(N\otimes L,M)\quad \text{for }M, N
\in D^{b}\cM(X)
\leqno(4.6.1)
$$
induced by {\rm (4.5.4) (4.1.5)}.
}

\medskip\noindent{\it Proof.}
The first two assertions are clear by the compatibility with the functor
$ \For $.
We define the inverse of (4.6.1) using (4.5.8) (4.1.5).
Then the assertion is reduced to the case
$ M, N \in \cM(X) $ by the exactness of
$ \otimes L, \otimes L^{*} $, and we can check the assertion using the
underlying perverse sheaves.

\medskip\noindent{{\bf 4.7.~Lemma.}} {\it
For
$ M, N \in \cM(X) $, we have
$$
\Ext^{k}(N\otimes M, \bD A_{X}^{\cM}) = 0\quad \text{for}\quad k < 0
\leqno(4.7.1)
$$
and
$$
\Ext^{0}(N\otimes M, \bD A_{X}^{\cM}) \to
\Ext^{0}(\For(N)\otimes \For(M), \bD A_{X(\bC)})
\leqno(4.7.2)
$$
is injective, where
$ \Ext^{0} $ on the left hand side is taken in
$ D^{b}\cM(X) $, and
$ \Ext^{0} $ on the right in
$ {D}_{c}^{b}(X(\bC),A) $.
}

\medskip\noindent{\it Proof.}
By adjunction for
$ a_{X} : X \to \Spec k $, we have
$$
\Ext^{k}(N\otimes M, \bD A_{X}^{\cM}) = \Ext^{k}((a_{X})_{!}(N\otimes M),
A^{\cM})
\leqno(4.7.3)
$$
in a compatible way with the functor
$ \For $, because
$$
\bD A_{X}^{\cM} = {a}_{X}^{!}A^{\cM}
\leqno(4.7.4)
$$
by (3.8) (3.1.3) (1.6.3).
So the assertion is reduced to
$$
H^{k}(a_{X})_{!}(\For(N)\otimes \For
(M)) = 0\quad \text{for}\quad k > 0,
\leqno(4.7.5)
$$
using (1.1.1) (applied to
$ \Hom(H^{0}(a_{X})_{!}(N\otimes M), A^{\cM})) $ and the compatibility
with
$ \For $.
But (4.7.5) follows from the vanishing of
$ \Ext^{k}(\For(N), \bD\For(M))) = 0 $ for
$ k < 0 $ (cf. [3]) by the isomorphism
$$
\Ext^{k}(\For(N)\otimes \For(M), \bD A_{X(\bC)}) =
\Ext^{k}(\For(N), \bD \For(M))).
\leqno(4.7.6)
$$

\medskip\noindent{{\bf 4.8.~Proposition.}} {\it
Let
$ \{U_{i}\} $ be an open covering of
$ X $, and
$ j_{I} : U_{I} := \bigcap_{i\in I} U_{i} \to X $ the natural inclusion.
For
$ M, N \in \cM(X) $, we have a spectral sequence
$$
{E}_{1}^{p,q} = \mopls_{|I|=p+1}\Ext^{q}({j}_{I}^{!}N\otimes
{j}_{I}^{!}M,\bD A_{{U}_{I}}^{\cM}) \Rightarrow
\Ext^{p+q}(N\otimes M, \bD A_{X}^{\cM}),
\leqno(4.8.1)
$$
such that
$ {E}_{1}^{p,q} = 0 $ for
$ p < 0 $.
In particular, a morphism
$ N\otimes M \to \bD A_{X}^{\cM} $ in
$ D^{b}\cM(X) $ is defined locally on
$ X $ for
$ M, N \in \cM(X) $.
}

\medskip\noindent{\it Proof.}
It is enough to show (4.8.1), because
$ {E}_{1}^{p,q} = 0 $ for
$ q < 0 $ by (4.7.1).
It is induced by the filtration
$ \sigma [9] $ on the dual of the
Cech complex (2.2.1) (applied to
$ N) $ by the next lemma, using the adjunction for
$ j_{I} $.

\medskip\noindent{{\bf 4.9.~Lemma.}} {\it
Let
$ i $ and
$ j $ be respectively a closed embedding and an open embedding of
varieties of
$ \cV $.
Then, for
$ M, N \in D^{b}\cM(X) $, we have canonical isomorphisms
}
$$
(j_{!}j^{!}N)\otimes M = j_{!}(j^{!}N\otimes j^{!}M) = j_{!}j^{!}
(N\otimes M)
\leqno(4.9.1)
$$
$$
(i_{*}i^{*}N)\otimes M = i_{*}(i^{*}N\otimes i^{*}M) = i_{*}i^{*}
(N\otimes M)
\leqno(4.9.2)
$$

\medskip\noindent
{\it Proof.} Let
$ j : X \to Y $.
The first assertion (4.9.1) follows from (3.10) (2.7) (3.1.4) (3.1.9)
using the diagram:
$$
\CD
X @= X @>>> Y
\\
@VVV @VVV @VVV
\\
X\times X @>>> X\times Y @>>> Y\times Y
\endCD
$$
The argument is similar for (4.9.2).

\medskip\noindent{{\bf 4.10.~Proposition.}} {\it
For
$ M \in D^{b}\cM(X) $, we have a canonical morphism
$$
\bD M\otimes M \to \bD A_{X}^{\cM}
\leqno(4.10.1)
$$
in the bounded derived category
$ D^{b}\cM(X) $, such that its image by the functor
$ \For $ is the natural pairing of
$ \bD \For(M) $ and
$ \For(M) $, which corresponds to the identity on
$ \bD \For(M) $ by the natural isomorphism
$$
\Ext^{k}(K, \bD L) = \Ext^{k}(K\otimes L, \bD A_{X(\bC)})
\leqno(4.10.2)
$$
for
$ K, L \in {D}_{c}^{b}(X(\bC),A) $.
}

\medskip\noindent{\it Proof.}
Assume first
$ M \in \cM(X) $.
Then (4.10.1) is unique by the injectivity of (4.7.2).
So the assertion is local by (4.8), and we may assume that
$ X $ has a closed embedding
$ i $ into smooth
$ Y \in \cV (e.g., Y = \bA^{n}) $.
Then we have the canonical morphism
$$
\bD i_{!}M\otimes i_{!}M \to \bD A_{X}^{\cM}
\leqno(4.10.3)
$$
by (1.6.7), and
$$
\bD i_{!}M\otimes i_{!}M = i_{!}(\bD M\otimes M)
\leqno(4.10.4)
$$
by an argument similar to (4.9).
So the assertion follows from the adjunction for
$ i $.

Now consider the general case.
We have the spectral sequence
$$
{E}_{1}^{p,q} = \Ext^{q}((\bD M\,\boxtimes\,M)^{-p}, \delta_{*}\bD
A_{X}^{\cM}) \Rightarrow \Ext^{p+q}(\bD M\otimes M, \bD A_{X}^{\cM}),
\leqno(4.10.5)
$$
by the filtration
$ \sigma [9] $ on
$ \bD M\,\boxtimes\,M $, using also the adjunction for the diagonal
morphism
$ \delta : X \to X\times X $.
We have a canonical element
$ e \in {E}_{1}^{0,0} $ by (4.10.1) (applied to each component
$ M^{k} $ of
$ M) $.
Since
$ {E}_{1}^{p,q} = 0 $ for
$ q < 0 $, it is enough to show that
$ e $ belongs to
$ \Ker(d_{1} : {E}_{1}^{0,0} \to {E}_{1}^{1,0}) $.
By the injectivity of (4.7.2), the assertion is reduced to that for the
underlying perverse sheaves, and we can check it using (4.10.2).
The compatibility of (4.10.1) with the functor
$ \For $ is checked using the spectral sequence (4.10.5) for the underlying
$ A $-complexes, which is defined by the convolution of filtrations on
$ \bD \For(M) $ and
$ \For(M) $ in [3].

\medskip\noindent{{\bf 4.11.~Theorem.}} {\it
For
$ M, N \in D^{b}\cM(X), we $ have a canonical isomorphism
$$
\Ext^{k}(N, \bD M) = \Ext^{k}(N\otimes M, \bD A_{X}^{\cM})
\leqno(4.11.1)
$$
induced by the composition with {\rm (4.10.1)}, and it is compatible
with {\rm (4.10.2)} by the functor
$ \For $.
}

\medskip\noindent{\it Proof.}
The last assertion follows from (4.10), because (4.10.2) is also induced
by the composition with the natural pairing.
Since (4.11.1) is functorial, we may assume
$ M, N \in \cM(X) $.

We proceed by induction on
$ \dim X $.
Let
$ U $ be an open subvariety of
$ X, Y = X \setminus U $ with the natural inclusions
$ i : Y \to X, j : U \to X $.
We have
$$
\Ext^{k}(i_{*}i^{*}N, \bD M) = \Ext^{k}(i^{*}N, \bD i^{*}M)
$$
$$
\Ext^{k}(j_{!}j^{!}N, \bD M) = \Ext^{k}(j^{!}N, \bD j^{!}M)
$$
by adjunction for
$ i, j $, and
$$
\Ext^{k}((i_{*}i^{*}N)\otimes M, \bD A_{X}^{\cM}) =
\Ext^{k}(i^{*}N\otimes i^{*}M, \bD A_{Y}^{\cM})
$$
$$
\Ext^{k}((j_{!}j^{!}N)\otimes M, \bD A_{X}^{\cM}) =
\Ext^{k}(j^{!}N\otimes j^{!}M, \bD A_{U}^{\cM})
$$
by (4.9) (4.7.4) and adjunction.
By inductive hypothesis, the assertion is reduced to that for
$ U $ using the distinguished triangle (3.3.2).
So we may assume
$ X $ smooth and pure dimensional, and
$ (H^{k}M)[- \dim X] $ smooth, by replacing
$ X $.
Then we may assume further
$ M $ smooth, and the assertion follows from (4.6.1), where
$ L, N, M $ in (4.6.1) are
$ M, N, \bD A_{X}^{\cM} $ respectively in this situation, so that (4.6.1)
is equivalent to (4.11.1).

\medskip\noindent{{\bf 4.12.~Theorem.}} {\it
Let
$ f : X \to Y $ be a morphism of varieties in
$ \cV $.
Then we have a canonical functorial morphism
$$
f_{!} \to f_{*}
\leqno(4.12.1)
$$
which corresponds to the natural morphism
$ f_{!} \to f_{*} $ by the functor
$ \For $.
In particular, {\rm (4.12.1)} is an isomorphism if
$ f $ is proper.
}

\medskip\noindent{\it Proof.}
Taking the direct image of the pairing (4.10.1), we have
$$
f_{!}\bD M\,\boxtimes\,f_{!}M \to \delta_{*}f_{!}\bD
A_{X}^{\cM} \to \delta_{*}\bD A_{Y}^{\cM}
\leqno(4.12.2)
$$
by (2.6--7), where the last morphism is induced by (3.9.2).
Here
$ \delta $ denotes the diagonal morphisms.
By (4.11.1), we get
$$
f_{!}\bD M \to \bD f_{!}M,
\leqno(4.12.3)
$$
by (1.6.8), and the assertion follows by replacing
$ M $ with
$ \bD M $.

\medskip\noindent
{\it Remark.} Using (4.10--11) instead of (1.6.7--8), Lemma (4.4) is
valid without assumption on the smoothness of varieties.

\newpage
\centerline{\bf \S 5. Vanishing Cycles}

\bigskip

We first give the proof of (3.5.2).

\medskip\noindent{{\bf 5.1.~Lemma.}} {\it
Let
$ S^{*} = \Spec k[t,t^{-1}] $.
Then we have
$ L_{i} \in D^{b}\cM(S^{*}) $ with morphisms
$ u_{i} : L_{i} \to L_{i+1} (i\ge 0) $ such that
$ \For(L_{i}) $ is an indecomposable local system of rank
$ i + 1 $ with unipotent monodromy,
$ \For(u_{i}) $ is injective, and
$ L_{0} = A_{{S}^{*}}^{\cM} $.
}

\medskip\noindent{\it Proof.}
Let
$ X = S^{*}\times S^{*} $ with
$ p : X \to S^{*} $ the second projection.
Let
$ i_{1} : S^{*} \to X $ be the diagonal morphism, and
$ i_{2} : S^{*} \to X $ the section of
$ \pi $ whose image is
$ \{t=1\}\times S^{*} $.
We define
$ \pi^{*}N = A_{{S}^{*}}^{\cM}\boxtimes N $ and
$ N \to \pi_{*}\pi^{*}N $ as in (3.6) for
$ N \in D^{b}\cM(S^{*}) $.
The composition
$$
N \to \pi_{*}\pi^{*}N \to \pi_{*}((i_{a})_{*}(i_{a})^{*})\pi^{*}N
\leqno(5.1.1)
$$
is an isomorphism by (3.6.5), where the last morphism is induced by
(3.3.2).
So we get morphisms
$$
u_{a} : \pi_{*}\pi^{*}N \to \pi_{*}((i_{a})_{*}(i_{a})^{*})
\pi^{*}N \simeq N\quad \text{for }a = 1, 2.
\leqno(5.1.2)
$$
Let
$ L = C(\mopls_{a=1,2} u_{a} : \pi_{*}\pi^{*}A_{{S}^{*}}^{H} \to
\mopls_{a=1,2}A_{{S}^{*}}^{\cM}) $ for
$ N = A_{{S}^{*}}^{\cM} $ so that we have a distinguished triangle
$$
\to \pi_{*}\pi^{*}A_{{S}^{*}}^{\cM} \to
\mopls_{a=1,2} A_{{S}^{*}}^{\cM} \to L \to
\leqno(5.1.3)
$$
Then
$ \For(L) $ is an indecomposable local system of rank
$ 2 $, and we have a morphism
$$
A_{{S}^{*}}^{\cM} \to L
\leqno(5.1.4)
$$
by the long exact sequence associated with the triangle.

We define
$ \otimes^{i} L $ by
$ (H^{1}\delta^{*}(L\,\boxtimes\,\cdots\,\boxtimes\,L))[-1], $ where
$ \delta : S^{*} \to $
$ S^{*}\times{\cdots} \times S^{*} $ is the diagonal morphism.
Its underlying
$ A $-complex is the
$ i $-times tensor of
$ \For(L) $.
We have an action of the symmetric group
$ S_{n} $ on
$ L\,\boxtimes\,\cdots\,\boxtimes\,L $ (and hence on
$ \otimes^{i} L) $, which is induced by (1.5.5) (and is well-defined
because it is uniquely determined by the underlying endomorphism on
$ \For(L)\,\boxtimes\,\cdots\,\boxtimes\, \For(L)) $.
We denote the action of
$ \sigma \in S_{n} $ by
$ \sigma_{*} $.
Then the symmetric tensor
$ \Sym^{i}L $ of
$ L $ is defined by the quotient of
$ \otimes^{i} L $ divided by
$ \sum_{\sigma \in S_{n}} \Im(\sigma_{*} - id) $ in the abelian category
$ \cM(X)[-1] $.
Let
$ L_{i} = \Sym^{i}L $.
Then
$ \For(L_{i}) $ is an indecomposable local system of rank
$ i + 1 $, and (5.1.4) induces
$ L_{i} \to L_{i+1} $ for
$ i \ge 0 $, where
$ L_{0} = A_{{S}^{*}}^{\cM} $.

\medskip\noindent
{\it Remark.} The above construction of
$ L $ is inspired by [4].

\medskip\noindent{{\bf 5.2.~Proposition.}} {\it
Let
$ X \in \cV, S = \Spec k[t] $, and
$ S^{*} $ as above.
Let
$ Y = X\times S, Y^{*} = X\times S^{*} $ with natural morphisms
$ p : Y \to S, j : Y^{*} \to Y, i : X\times \{0\} \to
Y $.
Then we have the nearby and vanishing cycle functors with unipotent
monodromy:
$$
\psi_{p,1},\,\,\,\varphi_{p,1} : \cM(Y) \to \cM
(X)
\leqno(5.2.1)
$$
together with functorial morphisms:
$$
\sp : H^{-1}i^{*}M \to \psi_{p,1}M,\,\,\,
\can : \psi_{p,1}M \to \varphi_{p,1}M
\leqno(5.2.2)
$$
for
$ M \in \cM(Y) $, which are compatible with the corresponding
functors
$^{p}\psi_{p,1},^{p}\varphi_{p,1} $ and morphisms
$ \sp $,
$ \can $ on the underlying perverse sheaves {\rm [10]} by the functor
$ \For $, where
$^{p}\psi_{p,1} := \psi_{p,1}[-1],^{p}\varphi_{p,1} :=
\varphi_{p,1}[-1] $ on the underlying perverse sheaves.
Moreover,
$ \psi_{p,1}M $ depends only on
$ j^{*}M $, and
$ \psi_{p,1}j_{*} : \cM(Y^{*}) \to \cM(X) $ is also denoted
by
$ \psi_{p,1} $.
}

\medskip\noindent{\it Proof.}
For
$ M' \in \cM(Y^{*}), we $ define
$$
M'_{i} = H^{0}\delta^{*}(M'\,\boxtimes\,L_{i}),
\leqno(5.2.3)
$$
where
$ \delta : X\times S^{*} \to X\times S^{*}\times S^{*} $ is induced
by the diagonal morphism of
$ S^{*} $.
For
$ k \ne 0 $, we have
$ H^{k}\delta^{*}(M'\,\boxtimes\,L_{i}) = 0 $ by the compatibility with
$ \For $ (cf. (3.4)).
So we get
$$
\delta_{*}M'_{i}= (\delta_{*}\delta^{*})(M'\,\boxtimes\,L_{i})\,\,\,
\text{in }D^{b}\cM(X\times S^{*}\times S^{*}).
\leqno(5.2.4)
$$
In particular, we have a canonical isomorphism
$$
M'_{0} = M'
\leqno(5.2.5)
$$
by (4.1.5).
If
$ M' = j^{*}M $ for
$ M \in \cM(Y) $, it induces a canonical morphism
$$
u_{i} : M \to j_{*}M'_{i}.
\leqno(5.2.6)
$$
Applying
$ j_{!}j^{!} $, we get
$$
v_{i} : j_{!}j^{!}M \to j_{!}M'_{i}
\leqno(5.2.7)
$$
with a canonical morphism of the mapping cones
$ C(v_{i}) \to C(u_{i}) $.
Using (1.4.6), we define
$ \psi_{p,1}M', \varphi_{p,1}M $ for
$ M' \in \cM(Y^{*}), M \in \cM(Y) $ so that
$$
i_{*}\psi_{p,1}M' = H^{-1}(C(j_{!}M'_{i} \to j_{*}M'_{i}))
\leqno(5.2.8)
$$
$$
i_{*}\varphi_{p,1}M = H^{-1}(C(C(v_{i}) \to C(u_{i})))
\leqno(5.2.9)
$$
for
$ i $ sufficiently large, and (5.2.2) is induced by the distinguished
triangle
$$
\to C(j_{!}j^{!}M \to M) \to C(j_{!}M'_{i} \to j_{*}M'_{i}) \to
C(C(v_{i}) \to C(u_{i})) \to .
\leqno(5.2.10)
$$
The independentness of
$ i $ and the compatibility with
$ \For $ are checked as in [32].

\medskip\noindent{{\bf 5.3.~Lemma.}} {\it
With the above notation, we have canonical functorial isomorphisms
$$
\psi_{p,1}(N\,\boxtimes\,A_{S}^{\cM}[1]) = N
\leqno(5.3.1)
$$
$$
\psi_{p,1}M = 0,\,\,\,i_{*}\varphi_{p,1}M = M
\leqno(5.3.2)
$$
for
$ N \in \cM(X), M \in \cM(Y) $ such that
$ \supp M \subset X\times \{0\} $.
}

\medskip\noindent{\it Proof.}
The last assertion is clear, because
$ j^{*}M = M'_{i} = 0 $ if
$ \supp M \subset X\times \{0\} $.
For (5.3.1), let
$ \pi : Y \to X $ be the projection, and
$ M = N\,\boxtimes\,A_{S}^{\cM} $.
It is enough to show
$$
H^{0}i^{*}M = N,
\leqno(5.3.3)
$$
because
$ sp $ in (5.2.2) is an isomorphism in this case by the compatibility
with the underlying morphism on the underlying perverse sheaves.
But (5.3.3) follows from (4.1.5), because
$$
(i_{*}i^{*})M = i_{*}H^{0}i^{*}M
\leqno(5.3.4)
$$
by definition (cf. (3.4)), where
$ H^{k}i^{*}M = 0 $ for
$ k \ne 0 $ is reduced to
$^{p}\cH^{k}i^{*} \For(M) = 0 $ for
$ k \ne 0 $.

\medskip\noindent
{\bf 5.4.~Definition.}
Let
$ X \in \cV $, and
$ g : X \to S := \Spec k[t] $ a function.
Let
$ i_{g} : X \to X\times S $ be the embedding by graph, and
$ j : X^{*} := X \setminus g^{-1}(0) \to X $ the natural inclusion.
We define the functors
$$
\psi_{g,1} : \cM(X) \to \cM(X), \varphi_{g,1}: \cM(X) \to \cM(X)
$$
by
$$
\psi_{g,1} = \psi_{p,1}(i_{g})_{*},\,\,\,\varphi_{g,1} = \varphi_{p,1}
(i_{g})_{*}.
\leqno(5.4.1)
$$
Here
$ \psi_{g,1}M $ depends only on
$ j^{*}M $ and
$ \psi_{g,1}j_{*} : \cM(X^{*}) \to \cM(X) $ is also denoted by
$ \psi_{g,1} $.
By (5.2.2), we have a canonical morphism
$$
\can : \psi_{g,1} \to \varphi_{g,1}.
\leqno(5.4.2)
$$
If
$ \supp M \subset g^{-1}(0) $, we have
$$
\psi_{g,1}M = 0,\,\,\,\varphi_{g,1}M = M
\leqno(5.4.3)
$$
by (5.3.2) using also (1.4.6) (1.4.8).
Here
$ \cM(Y) $ can be taken to be the target of the functors if
$ Y := g^{-1}(0) \in \cV $.
The compatibility with the above definition for
$ g = p $ is reduced to the compatibility with the direct image by
$ i_{g} $, and will be checked later.

\medskip\noindent{{\bf 5.5.~Proposition.}} {\it
With the above notation, we have a functor
$$
\xi_{g} : \cM(X) \to \cM(X)
$$
with functorial canonical exact sequences
$$
0 \to \psi_{g,1}M \to \xi_{g}M \to
M \to 0
\leqno(5.5.1)
$$
$$
0 \to j_{!}j^{!}M \to \xi_{g}M \to
\varphi_{g,1}M \to 0
\leqno(5.5.2)
$$
for
$ M \in \cM(X) $, where
$ j : U := X \setminus g^{-1}(0) \to X $ is the natural inclusion.
}

\medskip\noindent{\it Proof.}
Let
$ S = \Spec k[t], Y = X\times S $, and
$ i_{g} : X \to Y $ the embedding by graph of
$ g $.
Let
$ j_{g} : Y \setminus \Im \,i_{g} \to Y $ be the natural inclusion.
Then
$ \xi_{g} $ is defined by
$$
\xi_{g}M = \psi_{p,1}(j_{g})_{!}(j_{g})^{!}(M\,\boxtimes\,A_{S}^{\cM}
[1]).
\leqno(5.5.3)
$$
By (3.3.2), we have an exact sequence
$$
0 \to (i_{g})_{*}M \to (j_{g})_{!}(j_{g})^{!}(M\,\boxtimes\,A_{S}^{\cM}
[1]) \to M\,\boxtimes\,A_{S}^{\cM}[1] \to 0
\leqno(5.5.4)
$$
using the same argument as (5.3.3--4) (e.g.,
$ H^{k}(i_{g})^{*}(M\,\boxtimes\,A_{S}^{\cM}) = M $ if
$ k = 0 $ and
$ 0 $ otherwise).
Taking
$ \psi_{p,1} $, this induces (5.5.1) by (5.3.1).

The proof of (5.3.2) is same as in [28, 2.22] using the compatibility with
the functor
$ \For $.

\medskip\noindent
{\it Remark.} By (5.5.1--2), we have canonical isomorphisms
$$
\xi_{g}M = \varphi_{g,1}M = M
\leqno(5.5.5)
$$
for
$ M \in \cM(X) $ supported in
$ g^{-1}(0) $.

\medskip\noindent{{\bf 5.6.~Theorem.}} {\it
Let
$ i : X \to Y $ be a closed embedding of varieties of
$ \cV $.
Then
$$
i_{!} = i_{*} : D^{b}\cM(X) \to {D}_{X}^{b}\cM(Y),
\leqno(5.6.1)
$$
is an equivalence of categories.
}

\medskip\noindent{\it Proof.}
By (1.4.6), it is enough to show
$$
\Ext^{k}(M,N) = \Ext^{k}(i_{*}M,i_{*}N)
\leqno(5.6.2)
$$
for
$ M, N \in \cM(X) $, where the left (resp. right) hand side means
the extension group in
$ \cM(X) $ (resp.
$ \cM(Y)) $.
It is equivalent to that the functor
$$
N \to \Ext^{k}(i_{*}M,i_{*}N)
$$
is effaceable for
$ k > 0 [3] [2] $ (see also [32, 2.3]).
By (3.2), it is reduced to the case
$ X, Y $ affine, using an open covering
$ \{U_{\alpha}\} $ of
$ Y $ and an injective morphism
$ N \to \mopls_{\alpha} (j'_{\alpha})_{*}(j'_{\alpha})^{*}N $, where
$ j'_{\alpha} : X \cap U_{\alpha} \to X $,
and (1.4.4) (1.4.8) are also used.
We may further assume
$ X = g^{-1}(0) $ for a function
$ g $, factorizing the embedding
$ i $.
Then
$ \varphi_{g,1} $ is a quasi-inverse of
$$
D^{b}\cM_{X}(Y) \to {D}_{X}^{b}\cM(Y)
\leqno(5.6.3)
$$
using (5.4.3) and the quasi-isomorphisms
$ M \leftarrow \xi_{g}M \to \varphi_{g,1}M $ in (5.5.1--2) for
$ M \in {D}_{X}^{b}\cM(Y) $.

\medskip\noindent
{\it Remark.} Due to (5.6.1), the direct image by a closed embedding will
be sometimes omitted to simplify the notation.

\medskip

{\it Form now on we will use the results in} \S 3 {\it and} \S 4,
because (3.5.2) is proved above.

\medskip\noindent{{\bf 5.7.~Proposition.}} {\it
For a function
$ g : X \to S $ as in {\rm (5.4)}, let
$ X^{*} = X \setminus g^{-1}(0) $.
Then we have the nearby cycle functors
$$
\psi_{g}, \psi_{g,\ne 1} : \cM(X^{*}) \to \cM
(X)
\leqno(5.7.1)
$$
with the decomposition
$$
\psi_{g}M = \psi_{g,1}M \oplus \,\,\,\psi_{g,\ne 1}M\,\,\, \text{for }
M \in \cM(X^{*})
\leqno(5.7.2)
$$
They have a functorial endomorphism
$$
N : \psi_{g} \to \psi_{g}(-1) \,\text{\(same for }
\psi_{g,\ne 1}).
\leqno(5.7.3)
$$
These functors and endomorphism correspond to
$^{p}\psi_{g},^{p}\psi_{g,\ne 1} $ and
$ N $ on the underlying perverse sheaves by the functor
$ \For $, where
$ N = (2\pi i)^{-1}\otimes \log T_{u} $ with
$ T = T_{s}T_{u} $ the Jordan decomposition of the monodromy on the
underlying perverse sheaves, and
$^{p}\psi_{g,\ne 1} $ is the non-unipotent monodromy part of
$^{p}\psi_{g} $.
Here
$ (-1) $ is the Tate twist \(cf. {\rm (4.2))}.
}

\medskip\noindent{\it Proof.}
Let
$ S^{*} = \Spec k[t, t^{-1}] $, and
$ \pi : S^{*} \to S^{*} $ the
$ n $-fold covering such that
$ \pi^{*}t = t^{n} $.
We have a canonical morphism
$$
A_{{S}^{*}}^{\cM} \simto \pi_{*}A_{{S}^{*}}^{\cM}
\leqno(5.7.4)
$$
by (3.9.2).
Dualizing this, we get
$$
\pi_{*}A_{{S}^{*}}^{\cM} \simto A_{{S}^{*}}^{\cM}
\leqno(5.7.5)
$$
by (4.3) (4.12), whose composition with (5.7.4) is the identity on
$ A_{{S}^{*}}^{\cM} $.
In the notation of (5.2), let
$$
\psi_{p}M = \psi_{p,1}(M\otimes p'^{*}\pi_{*}A_{{S}^{*}}^{\cM}))\quad
\text{for}\quad M \in \cM(X^{*})
\leqno(5.7.6)
$$
for
$ n $ sufficiently divisible, where
$ p' $ is the restriction of
$ p $ to
$ X^{*} $.
Then we define
$ \psi_{p}M $ as in (5.4.1).
The independentness of
$ n $ (sufficiently divisible) is checked using the direct image
of (5.7.4), because the monodromy of
$^{p}\psi_{g}(\For(M)) $ is quasi-unipotent as a consequence of (1.1.2)
(this is checked by reducing to the normal crossing case [28, \S 3]).
We have
$ \psi_{g,\ne 1} $ with the decomposition(5.7.2) by (5.7.4--5).

Now it remains to construct
$ N $.
It is enough to do it on
$ \psi_{p,1} $ by (5.7.6).
With the notation of (4.2), we have an isomorphism
$$
H^{1}(a_{S^{*}})_{*}A_{{S}^{*}}^{\cM} = A^{\cM}(-1)
\leqno(5.7.7)
$$

\noindent
(canonical up to a sign).
In fact, we have natural inclusions
$ j : S^{*} \to X = \bP^{1}, i_{a} : \{P_{a}\} \to X
(a = 1, 2) $ with
$ \{P_{1}, P_{2}\} = X \setminus S^{*} $, and a distinguished triangle
$$
\to \mopls_{a=1,2} (i_{a})_{!}(i_{a})^{!}A_{X}^{\cM} \to
A_{X}^{\cM} \to j_{*}A_{{S}^{*}}^{\cM} \to
$$
by (3.3.2).
Then (5.7.7) follows from the associated long exact sequence, because
(4.3) implies
$ (i_{a})^{!}A_{X}^{\cM} = A^{\cM}(-1)[-2] $ using (3.1.3).
By construction of
$ L $ in (5.1), (5.7.7) implies a surjective morphism
$$
L \to A_{{S}^{*}}^{\cM}(-1)
\leqno(5.7.8)
$$
whose kernel is
$ A_{{S}^{*}}^{\cM} $.
This is extended to
$$
L_{i} \to L_{i-1}(-1)
\leqno(5.7.9)
$$
which is induced by the sum of the morphisms
$$
L\,\boxtimes\,\cdots\,\boxtimes\,L \to L\,\boxtimes\,\cdots\,\boxtimes\,
A_{{S}^{*}}^{\cM}(-1)\,\boxtimes\,\cdots\,\boxtimes\,L.
$$
We can check that (5.7.9) induces the morphism
$ N $.

\medskip\noindent{{\bf 5.8.~Proposition.}} {\it
We have a canonical functorial morphism
$$
\Var : \varphi_{g,1} \to \psi_{g,1}(-1)
\leqno(5.8.1)
$$
whose image by the functor
$ \For $ coincides with
$ \Var $ in {\rm [27] [28]}.
In particular,
$ \Var\scirc\can = N $,
$ \can\scirc\Var = N $.
}

\medskip\noindent{\it Proof.}
With the notation of (5.2), we have a morphism
$ C(u_{i}) \to M'_{i-1}(-1) $ by (5.7.9), and this induces (5.8.1).

\medskip\noindent{{\bf 5.9.~Proposition.}} {\it
Let
$ i : X \to Y $ be as in {\rm (3.3)}, and assume
$ X = g^{-1}(0) $ for a function
$ g $.
Then
$ we $ have canonical functorial isomorphisms
$$
i_{*}i^{*} = C(\can : \psi_{g,1} \to \varphi_{g,1}),\,\,\,
i_{!}i^{!} = C(-\Var : \varphi_{g,1} \to \psi_{g,1}(-1))[-1].
\leqno(5.9.1)
$$
whose image by the functor
$ \For $ is the natural isomorphism on the underlying
$ A $-complexes in {\rm [28] [32]}.
}

\medskip\noindent{\it Proof.}
We show the assertion for
$ i_{*}i^{*} $, because a similar argument holds for
$ i_{!}i^{!} $.
Since
$ i_{*}i^{*} = C(j_{!}j^{!} \to id) $ and we have
$ id \leftarrow \xi_{g} \to \varphi_{g,1}, $ it is enough to show an
isomorphism
$$
\varphi_{g,1}C(j_{!}j^{!}\to id) = C(\can : \psi_{g,1}\to\varphi_{g,1}).
\leqno(5.9.2)
$$
But this follows from the isomorphism
$$
\can : \psi_{g,1}j_{!}j^{!} \to \varphi_{g,1}j_{!}j^{!}.
\leqno(5.9.3)
$$
The compatibility with
$ \For $ can be checked using (5.5.5).

\medskip\noindent{{\bf 5.10.~Proposition.}} {\it
Let
$ f : X \to Y $ be a proper morphism of varieties of
$ \cV, g $ a function on
$ Y $, and
$ h = gf $.
Then we have canonical isomorphisms
$$
f_{*}\scirc\psi_{h} = \psi_{g}\scirc f_{*},\,\,\,
f_{*}\scirc\varphi_{h,1} = \varphi_{g,1}\scirc f_{*},
\leqno(5.10.1)
$$
whose images by the functor
$ \For $ coincide with the natural isomorphisms of the underlying perverse
sheaves.
}

\medskip\noindent{\it Proof.}
With the notation of (5.2), let
$ X' = X\times S, Y' = Y\times S, f' = f\times id $ with projections
$ p : X' \to S, q : Y' \to S $.
Then it is enough to show
$$
f'_{*}\scirc\psi_{p} = \psi_{q}\scirc f'_{*},\,\,\,
f'_{*}\scirc\varphi_{p,1} = \varphi_{q,1}\scirc f'_{*}.
\leqno(5.10.2)
$$
Since the external product commutes with direct image (2.7), the assertion
follows from the commutativity of pull-back with direct image (3.10).

\medskip\noindent
{\it Remark.} We can check
$$
\psi_{g,1}M' = i_{*}H^{-1}i^{*}j_{*}(M'\otimes g'^{*}L_{i}),\,\,\,
\varphi_{g,1}M = i_{*}H^{-1}i^{*}(C(M \to j_{*}(M'\otimes g'^{*}L_{i}))
,
\leqno(5.10.3)
$$
using the same argument as above (applied to the direct image by the
embedding by graph of
$ g) $.

\medskip\noindent{{\bf 5.11.~Proposition.}} {\it
We have canonical functorial isomorphisms
$$
\psi_{g}\scirc \bD = (\bD\scirc \psi_{g})(1),\,\,\,
\varphi_{g,1}\scirc \bD = \bD\scirc \varphi_{g,1},
\leqno(5.11.1)
$$
whose image by the functor
$ \For $ is the natural isomorphism on the underlying perverse sheaves
in {\rm [29]}.
}

\medskip\noindent{\it Proof.}
Let
$ \bD M\otimes M \to \bD A_{X}^{\cM} $ be the canonical pairing
in (4.10).
We can construct natural pairings
$$
\psi_{g}\bD M\otimes \psi_{g}M \to \bD
A_{X}^{\cM}(1),\,\,\,\varphi_{g,1}\bD M\otimes \varphi_{g,1}M \to
\bD A_{X}^{\cM}
\leqno(5.11.2)
$$
by the same argument as in [32], using the natural morphism
$ L_{i}\otimes L_{j} \to L_{i+j} $.
Then the assertion follows from (4.11), because their underlying pairings
are perfect.

\medskip\noindent{{\bf 5.12.~Proposition.}} {\it
Let
$ i : X \to Y $ be a closed embedding of varieties of
$ \cV $ such that
$ X $ is a locally principal divisor on
$ Y $.
Let
$ M \in \cM(Y) $ such that
$ \varphi_{g}M = 0 $ for a local defining equation
$ g $ of
$ X $.
Then we have a natural isomorphism
$$
i^{!}M = i^{*}M(-1)[-2]
\leqno(5.12.1)
$$
compatible with the natural isomorphism on the underlying
$ A $-complexes.
}

\medskip\noindent{\it Proof.}
Replacing
$ Y $ by the line bundle associated with the divisor
$ X $, and
$ X $ by the zero section of the line bundle, we may assume
$$
i^{!}A_{Y}^{\cM} = A_{X}^{\cM}(-1)[-2],
\leqno(5.12.2)
$$
because it holds for the underlying
$ A $-complex and we have (5.12.2) restricting to a dense open subvariety so
that the line bundle becomes trivial.
We have a canonical morphism
$$
i_{*}i^{!}A_{Y}^{\cM}\otimes M \to A_{Y}^{\cM}\otimes M = M
\leqno(5.12.3)
$$
by (3.3.2) (4.1.5).
By the same argument as (4.9.2), we have
$$
i_{*}i^{!}A_{Y}^{\cM}\otimes M = i_{*}(i^{!}A_{Y}^{\cM}\otimes i^{*}M).
\leqno(5.12.4)
$$
So we get a morphism
$$
i^{*}M(-1)[-2] = i^{!}A_{Y}^{\cM}\otimes i^{*}M \to i^{!}M
\leqno(5.12.5)
$$
by (5.12.2) (4.1.5) and adjunction.
This is an isomorphism, because its underlying morphism is.
In fact, we can check it using (5.9) and the construction in [27, 5.2].

\bigskip\bigskip
\centerline{\bf \S 6. Weight}

\bigskip
\noindent
{\bf 6.1.~Definition.} We say that
$ M \in \cM(X) $ has {\it weights}
$ \le n $ (resp.
$ \ge n) $ if
$ \Gr_{k}^{W}M = 0 $ for
$ k > n $ (resp.
$ k < n) $, and
$ M $ is pure of weight
$ n $ if
$ \Gr_{k}^{W}M = 0 $ for
$ k \ne n $.
We denote by
$ \cM(X,n) $ the full subcategory of
$ \cM(X) $ consisting of pure objects with weight
$ n $.
By (1.1.3--4),
$ \cM(X,n) $ is a semi-simple abelian full subcategory of
$ \cM(X) $.

We say that
$ M \in D^{b}\cM(X) $ has {\it weights}
$ \le n $ (resp.
$ \ge n) $ if
$ H^{i}M $ has weights
$ \le n + i $ (resp.
$ \ge n + i) $, and
$ M $ is pure of weight
$ n $ if
$ H^{i}M $ is pure of weight
$ n + i $.

Let
$ M \in \cM(X) $.
We say
$ \supp M = Z $ if
$ X \setminus Z $ is the maximal open subvariety of
$ X $ such that
$ j^{*}M = 0 $ where
$ j : X \setminus Z \to X $.
We say that
$ M $ has {\it strict support}
$ Z $ if
$ \supp M = Z $ and
$ M $ has neither nontrivial sub or quotient object with strictly smaller
support.
Let
$ \cM_{Z}(X,n) $ denote the full subcategory of
$ \cM(X,n) $ consisting of the objects with strict support
$ Z $ or
$ \emptyset $.
(Note that
$ \cM_{Z}(X,n) \ne \cM(X,n) \cap \cM_{Z}(X) $ in the notation
of (1.4.6).)
By semi-simplicity of
$ \cM(X,n), \cM_{Z}(X,n) $ is a semi-simple abelian full subcategory
of
$ \cM(X) $.

\medskip\noindent{{\bf 6.2.~Lemma.}} {\it
For
$ M \in \cM_{Z}(X,m), N \in \cM_{Z'}(X,n) $, we have
$$
\Hom(M,N) = 0\,\,\, \text{unless }Z = Z', m = n.
\leqno(6.2.1)
$$
For
$ M \in \cM(X,n) $, we have a unique decomposition:
$$
M = \mopls_{Z} M_{Z}
\leqno(6.2.2)
$$
such that
$ M_{Z} \in \cM_{Z}(X,n) $ \(called the decomposition by strict
support\,\), where
$ Z $ runs over irreducible closed subvarieties of
$ X $.
}

\medskip\noindent{\it Proof.}
The first assertion is clear by definition.
By semisimplicity, we get the second except for the irreducibility of
$ Z $.
But it follows from the next lemma, combined with (2.2.3) applied to an
open subset of
$ Z $.

\medskip\noindent{{\bf 6.3.~Lemma.}} {\it
Let
$ f : X \to Y $ be a locally closed embedding of varieties of
$ \cV $, i.e.,
$ f $ is the composition of a closed embedding
$ i : X \to U $ with an open embedding
$ j : U \to Y $.
Let
$ Z $ be the closure of
$ f(X) $.
Then we have a functor \(called the intermediate direct image\):
$$
f_{!*} : \cM(X) \to \cM_{Z}(Y)
\leqno(6.3.1)
$$
with a functorial isomorphism
$$
f^{*}f_{!*} = id,
\leqno(6.3.2)
$$
such that
$ f_{!*}(M) $ has neither nontrivial sub or quotient object supported in
$ Y \setminus f(X) $.
Moreover, if
$ M \in \cM_{Z}(Y) $ has neither nontrivial sub or quotient object
supported in
$ Y \setminus f(X) $, we have a canonical isomorphism
}
$$
f_{!*}f^{*}M = M.
\leqno(6.3.3)
$$

\medskip\noindent
{\it Proof.} We define
$$
f_{!*}(M) = \Im(H^{0}f_{!}M \to H^{0}f_{*}M).
\leqno(6.3.4)
$$
using (4.12.1).
By adjunction (3.1.2),
$ H^{0}f_{*}M $ has no nontrivial subobject support in
$ Y \setminus f(X) $, because
$ H^{k}f_{*}M = 0 $ for
$ k < 0 $.
So
$ f_{!*}(M) $ has no nontrivial subobject support in
$ Y \setminus f(X) $.
The dual argument shows the assertion on the quotient.
We have (6.3.2), because
$ j^{*}f_{!*}M = i_{*}M $ for
$ M \in \cM(X) $ and
$ i^{*}i_{*} = id $ (cf. (3.7.3)).
For
$ M \in \cM_{Z}(Y) $, we have
$ f^{!}M, f^{*}M \in \cM(X) $, and the composition of the natural
morphisms
$ H^{0}f_{!}f^{!}M \to M, M \to H^{0}f_{*}f^{*}M $ coincides
with the morphism in (6.3.4), and
$ M $ is identified with the image of the morphism.
So we get the last assertion

\medskip\noindent{{\bf 6.4.~Proposition.}} {\it
Let
$ X \in \cV $, and
$ Z $ be an irreducible closed subvariety of
$ X $ in
$ \cV(k) $.
Then we have an equivalence of categories
$$
\cM_{Z}(X,n) = \varinjlim_{U}\cM_{U}(X\setminus (Z\setminus U),n)_{s}
\leqno(6.4.1)
$$
where
$ U $ runs over nonempty smooth open subvarieties of
$ Z $, and the morphisms of the inductive system are induced by open
pull-backs.
Here
$ \cM_{U}(X\setminus (Z\setminus U),n)_{s} $ is the full subcategory
of
$ \cM_{U}(X\setminus (Z\setminus U),n) $ consisting of objects
$ M $ such that
$ \For(M)[- \dim Z] $ is a local system on
$ U $.
}

\medskip\noindent{\it Proof.}
Let
$ j : X\setminus (Z\setminus U) \to X $ denote the natural
inclusion.
Then (6.4.1) is induced by
$ j^{*} $, and its quasi-inverse by
$ j_{!*} $ using (6.3).

\medskip\noindent
{\it Remark.} If
$ Z $ contains an open subvariety
$ U $ belonging to
$ \cV, (6.4.1) $ becomes
$$
\cM_{Z}(X,n) = \varinjlim_{U}\cM
(U,n)_{s}
\leqno(6.4.2)
$$
by (1.4.6), where
$ U $ runs over nonempty smooth open subvarieties of
$ Z $ belonging to
$ \cV $.

\medskip\noindent{{\bf 6.5.~Lemma.}} {\it
Let
$ M \in D^{b}\cM(X) $, and
$ N \in \cM(Y) $ such that
$ N $ is nonzero and pure of weight
$ k $.
Then
$ M $ has weight
$ \le n $ \(resp.
$ \ge n) $ if and only if one of the following equivalent conditions
is satisfied:

\noindent
{\rm (i)}
$ M[j] $ has weight
$ \le n + j $ \(resp.
$ \ge n + j) $,

\noindent
{\rm (ii)}
$ M (i) $ has weight
$ \le n - 2i $ \(resp.
$ \ge n - 2i) $,

\noindent
{\rm (iii)}
$ \bD M $ has weight
$ \ge -n $ \(resp.
$ \le -n) $,

\noindent
{\rm (iv)}
$ M\,\boxtimes\,N $ has weight
$ \le n + k $ \(resp.
$ \ge n + k) $.
}

\medskip\noindent{\it Proof.}
The assertion for (i) (iii) (iv) is clear by definition and (1.2.3), (1.5.3).
As for (ii), it follows from the last assertion of the next lemma.

\medskip\noindent{{\bf 6.6.~Proposition.}} {\it
If
$ X \in \cV $ is smooth
$ A_{X}^{\cM} $ is pure of weight
$ 0 $, and
$ A^{\cM}(-1) $ is pure of weight
$ 2 $.
}

\medskip\noindent{\it Proof.}
We assume first
$ X_{\bC} $ irreducible so that
$ A_{X(\bC)}[\dim X] $ and hence
$ A_{X}^{\cM}[\dim X] \in \cM(X) $ are simple.
Then the weight filtration on
$ A_{X}^{\cM}[\dim X] $ is trivial, and
$ A_{X}^{\cM} $ is pure.
We have a natural nonzero morphism
$$
A^{\cM} \to H^{0}(a_{X})_{*}A_{X}^{\cM}
$$
by (3.9.2).
Using the last assertion of the next theorem, this implies that
$ A_{X}^{\cM} $ for
$ X = \bP^{1} $ is pure of weight
$ 0 $ (because
$ A^{\cM} $ is pure of weight
$ 0 $ by (1.6.2)), and the last assertion also follows by definition of
$ A^{\cM}(-1) $ (cf. (4.2.1)).
Then the weight of
$ A_{X}^{\cM} $ is determined by (4.3.1) using (iii) of (6.5) and the
last assertion.

For general
$ X $, we may restrict to any dense open subvariety, considering the
weight filtration on the local system
$ \For(A_{X}^{\cM}) $.
Then we have an \'etale morphism to an affine space, and the assertion
follows from the next theorem on the pull-back and (4.4.1).

\medskip\noindent{{\bf 6.7.~Theorem.}} {\it
Let
$ f : X \to Y $ be a morphism of varieties of
$ \cV $.
Then
$ f_{*}M, f^{!}N $ \(resp.
$ f_{!}M, f^{*}N) $ have weights
$ \ge n $ \(resp.
$ \le n) $ if so are
$ M \in D^{b}\cM(X), N \in D^{b}\cM(Y) $.
In particular, if
$ f $ is proper,
$ f_{!}M = f_{*}M $ is pure of weight
$ n $ if so is
$ M $.
}

\medskip\noindent{\it Proof.}
The assertion for the direct image by an affine morphism is clear by (1.4.3),
using (iii) of (6.5) for
$ f_{!} $.
This implies
$$
W_{k}i_{*}M = i_{*}W_{k}M
\leqno(6.7.1)
$$
for a closed embedding
$ i $.

We first prove the assertion on the pull-back by a closed embedding
$ i $.
By (6.7.1), it is enough to show the assertion for
$ i_{*}i^{*}, i_{!}i^{!} $, where
$ X $ may not necessarily belong to
$ \cV $ (cf. (3.3)).
We show the assertion for
$ i_{*}i^{*} $, because the dual argument holds for
$ i_{!}i^{!} $.
Since the assertion is local, we may assume
$ X $ is an intersection of principal divisors defined by
$ g_{k} $.
Using the open covering
$ \{U_{k}\} $ of
$ Y $ such that
$ U_{k} = Y \setminus {g}_{k}^{-1}(0),i_{*}i^{*} $ is isomorphic to the
composition of
$ (i_{k})_{*}(i_{k})^{*} $ by (2.2.2), where
$ i_{k} : {g}_{k}^{-1}(0) \to Y $.
So we may assume
$ X = g^{-1}(0) $ for a function
$ g $.
We may also assume
$ N \in \cM(Y) $ and
$ N $ is pure of weight
$ n $.
By (6.7.1), it is enough to show
$$
\Ker(j_{!}j^{!}N \to N)\,\,\, \text{has weights } \le
n -1,
\leqno(6.7.2)
$$
where
$ j : Y \setminus X \to Y $.
By (6.2.2), we may assume
$ N $ has strict support
$ Z $ not contained in
$ g^{-1}(0) $.
Since
$ j_{!}j^{!}N $ has weights
$ \le n, (6.7.2) $ is equivalent to
$$
\Gr_{n}^{W}j_{!}j^{!}N = N.
\leqno(6.7.3)
$$
By semisimplicity, it is enough to show
$$
\Hom(j_{!}j^{!}N,M) = 0\quad \text{if}\quad \supp M \subset X,
\leqno(6.7.4)
$$
and follows from the adjunction.

Now we show the assertion on the direct image.
Let
$ U $ be an affine open subvariety of
$ X $, and
$ Z = X \setminus U $ with natural inclusions
$ i : Z \to X, j : U \to X $.
Then we have a distinguished triangle
$$
\to i_{!}i^{!}M \to M \to j_{*}j^{*}M
\to ,
\leqno(6.7.5)
$$
and
$ i_{!}i^{!}M $ has weights
$ \ge n $ by the above argument.
Since the assertion is true for
$ (fj)_{*} $, we can proceed by induction on dimension of support.
(This completes the proof of (6.6) in the case
$ X_{\bC} $ irreducible.)

It remains to show the assertion for the pull-back by a smooth projection
whose fiber is an affine space, taking a factorization of
$ f $, because the assertion is local.
Using (iv) of (6.5), it follows from the assertion of (6.6) in the case
$ X_{\bC} $ irreducible.

\medskip\noindent{{\bf 6.8.~Proposition.}} {\it
Let
$ M \in D^{b}\cM(X) $, and assume
$ M $ is pure.
Then we have a non canonical isomorphism
}
$$
M = \mopls_{i} (H^{i}M)[-i]\,\,\, \text{in }D^{b}\cM(X).
\leqno(6.8.1)
$$

\medskip\noindent
{\it Proof.} This follows from the next lemma.

\medskip\noindent{{\bf 6.9.~Lemma.}} {\it
If
$ M, N \in D^{b}\cM(X) $ are pure of weight
$ m, n $.
Then
}
$$
\Ext^{i}(M,N) = 0\,\,\, \text{\it for }m < n + i.
\leqno(6.9.1)
$$

\medskip\noindent
{\it Proof.} This follows from the semisimplicity of pure objects of
$ D^{b}\cM(X) $ (cf. for example, [30, II, (4.5)]).

\medskip\noindent{{\bf 6.10.~Proposition.}} {\it
Let
$ f : X \to Y $ be a proper morphism of varieties of
$ \cV $.
If
$ M \in D^{b}\cM(X) $ is pure, we have a non canonical isomorphism
}
$$
f_{*}M = \mopls_{i} (H^{i}f_{*}M)[-i]\,\,\, \text{in }D^{b}\cM(Y).
\leqno(6.10.1)
$$

\medskip\noindent
{\it Proof.} This follows from (6.7--8).

\medskip\noindent{{\bf 6.11.~Proposition.}} {\it
Let
$ g : X \to S $ be as in {\rm (5.4)}, and
$ M \in \cM(X) $.
Let
$ L_{i}\psi_{g}M = \psi_{g}W_{i+1}M, L_{i}\varphi_{g,1}M = $
$ \varphi_{g,1}W_{i}M $.
Then the weight filtration
$ W $ on
$ \psi_{g}M, \varphi_{g,1}M $ is the monodromy filtration relative to
$ L $ {\rm [13]}, i.e.,
}
$$
N^{k} : \Gr_{i+k}^{W}\Gr_{i}^{L}\psi_{g}M \simto (\Gr_{i-k}^{W}
\Gr_{i}^{L}\psi_{g}M)(-k)\,\,\, \text{\it \(same for }\varphi_{g,1}).
\leqno(6.11.1)
$$

\medskip\noindent
{\it Remark.} In the case of
$ l $-adic perverse sheaves over a variety of characteristic
$ p $ [3], the above result was obtained by Gabber (unpublished).
The following proof is inspired by a discussion with Deligne, and should be
essentially the same as Gabber's proof.

\medskip\noindent{\it Proof of} (6.11).
By definition of relative monodromy filtration, we may assume
$ M $ is pure of weight
$ n $, and also
$ M $ has strict support
$ Z $ by (6.2.2).
If
$ Z \subset Y := g^{-1}(0) $, the assertion is clear.
So we may assume
$ Z \not\subset Y $.
Then it is enough to show the assertion for
$ \psi_{g}M $ by [27, 5.1.12], because
$ \can $ is surjective in this case. By duality (5.11.1), it is enough
to show
the injectivity of
$$
N^{k} : \Gr_{n-1+k}^{W}\psi_{g}M \to (\Gr_{n-1-k}^{W}\psi_{g}M)(-k).
\leqno(6.11.2)
$$
By definition, we may replace
$ g : X \to S $ by the projection
$ p : X\times S \to S $ in (5.2).
By (5.7.6), we can replace
$ \psi_{p}M $ by
$ \psi_{p,1}M $, because
$ \pi_{*}A_{{S}^{*}}^{\cM} $ is pure of weight zero, and
$ \delta^{*}(M\,\boxtimes\,\pi_{*}A_{{S}^{*}}^{\cM}) $ is pure by (6.14)
below, where we have
$$
M\otimes p'^{*}\pi_{*}A_{{S}^{*}}^{\cM} = \delta^{*}(M\,\boxtimes\,
\pi_{*}A_{{S}^{*}}^{\cM})
\leqno(6.11.3)
$$
by (3.1.4) (3.6.5).

Let
$ i : X\times \{0\} \to X\times S $ be the natural embedding.
Then we have an isomorphism
$$
H^{-1}i^{*}M = \Ker\,N \subset \psi_{p,1}M
\leqno(6.11.4)
$$
by (5.9.1), because
$ \Var $ is injective in this case (cf. [27]).
So
$ \Ker\,N $ has weights
$ \le n - 1 $ by (6.7), and
$$
\Ker(N : \Gr_{j}^{W}\psi_{p,1}M \to (\Gr_{j-2}^{W}\psi_{p,1}M)(-1)) = 0
\,\,\, \text{for }j > n - 1.
\leqno(6.11.5)
$$
If (6.11.2) is not injective for some
$ k > 0 $, let
$ j' $ be the minimal number such that
$$
N^{j} : \Gr_{n-1+k}^{W}\psi_{p,1}M \to (\Gr_{n-1+k-2j}^{W}
\psi_{p,1}M)(-j)
$$
is not injective, and
$ M' $ be the subobject of
$ \Gr^{W}\psi_{p,1}M $ defined by
$$
\sum_{i\ge 0} N^{i}(\Ker(N^{j'} : \Gr_{n-1+k}^{W}\psi_{p,1}M \to
(\Gr_{n-1+k-2j'}^{W}\psi_{p,1}M)(-j'))(i) .
$$
Let
$ K_{j} $ be the degree
$ j $ part of
$$
\Ker(N\,\boxtimes\,id + id\,\boxtimes\,N : M'\,\boxtimes\,M' \to
M'\,\boxtimes\,M'(-1)).
$$
Then
$ K_{j} \ne 0 $ for some
$ j > 2n - 2 $.
This contradicts to (6.11.5) applied to
$ H^{-1}\delta^{*}(M\,\boxtimes\,M) $, combined with the next lemma.

\medskip\noindent{{\bf 6.12.~Lemma.}} {\it
Let
$ X, Y \in \cV, S = \Spec k[t] $, and
$ X' = X\times S, Y' = Y\times S $ with projections
$ p : X' \to S, q : Y' \to S $.
Let
$ Z = X\times Y\times S $ with projection
$ h : Z \to S $.
Let
$ \delta : Z \to X'\times Y' $ be the morphism induced by the diagonal
morphism of
$ S $.
Let
$ M \in \cM(X',m), M' \in \cM(Y',n) $.
Then there is an open subvariety
$ U $ of
$ Z $ containing
$ X\times Y\times \{0\} $ such that
$ j^{*}H^{i}\delta^{*}(M\,\boxtimes\,M') = 0 $ for
$ i \ne -1 $ and
$ j^{*}H^{-1}\delta^{*}(M\,\boxtimes\,M') $ is pure of weight
$ n + m - 1 $, where
$ j : U \setminus X\times Y\times \{0\} \to Z $.
Moreover, we have a natural isomorphism
$$
(\psi_{p,1}M,N)\,\boxtimes\,(\psi_{q,1}M',N) =
(\psi_{h,1}H^{-1}\delta^{*}(M\,\boxtimes\,M'),N)
\leqno(6.12.1)
$$
where
$ N $ on the right hand side is defined by
$ N\,\boxtimes\,id + id\,\boxtimes\,N $.
}

\medskip\noindent{\it Proof.}
Let
$ S\times S \to S $ be the morphism defined by
$ (t_{1},t_{2}) \to t_{1} - t_{2} $, and
$ g : X\times Y \to S $ its composition with
$ p\times q $ so that
$ g^{-1}(0) = Z $.
Let
$$
U = Z \setminus (\supp \psi_{g,1}(M\,\boxtimes\,M') \setminus
X\times Y\times \{0\}).
$$
Here
$ \supp \psi_{g,1}(M\,\boxtimes\,M') \setminus X\times Y\times \{0\} $ is
a closed subvariety of
$ Z $, because it is true after taking the base change by
$ k \to \bC $ (using Milnor fibration which follows
from [18]).
Then the first assertion follows from (5.9.1) and the lemma below.
Since
$ i^{*}j_{*} = C(j_{!} \to j_{*}) $ in the notation of (3.3.2), we get the
morphism (6.12.1) by the morphism
$ id \to \delta_{*}\delta^{*} $ in (3.3.2) and the natural morphism
$ L_{i}\otimes L_{j} \to L_{i+j} $.
Then the morphism is an isomorphism, because the underlying morphism of
$ A $-complex is an isomorphism using Milnor fibration (cf. [loc. cit.]).

\medskip\noindent{{\bf 6.13.~Lemma.}} {\it
With the notation and assumption of {\rm (5.12)}, assume further
$ M \in \cM(Y) $ is pure of weight
$ n $.
Then
$ i^{!}M = i^{*}M(-1)[-2] $ is also pure of weight
$ n $.
}

\medskip\noindent{\it Proof.}
This follows from (6.7).

\medskip\noindent{{\bf 6.14.~Proposition.}} {\it
Let
$ M \in D^{b}\cM(X) $.
Then
$ M $ has weights
$ \le n $ \(resp.
$ \ge n) $ if and only if
$ (i_{x})^{*}M $ \(resp.
$ (i_{x})^{!}M) $ has weights
$ \le n $ \(resp.
$ \ge n) $ for any closed points
$ x $ of
$ X $, where
$ i_{x} : \{x\} = \Spec k(x) \to X $ is the natural morphism.
}

\medskip\noindent{\it Proof.}
By duality, it is enough to show the assertion for weights
$ \le n $.
If
$ M $ has weights
$ \le n, (i_{x})^{*}M $ has weights
$ \le n $ by (6.7).

Assume
$ (i_{x})^{*}M $ has weights
$ \le n $ for any closed point
$ x $ of
$ X $.
Let
$ U $ be a pure dimensional smooth open subvariety of
$ X $ with natural inclusion
$ j : U \to X $ such that
$ (j^{!}H^{k}M)[- \dim U] $ are smooth (cf. (4.5)).
Let
$ Z = X \setminus U $ with
$ i : Z \to X $.
By induction on
$ \dim X $, it is enough to show that
$ j^{!}M $ has weights
$ \le n $ using the triangle (3.3.2) and (6.7).
So we may assume
$ (H^{k}M)[- \dim U] $ are smooth.
If
$ M $ is smooth,
$ H^{k}(i_{x})^{*}M $
$ = 0 $ for
$ k \ne 0 $.
So we may assume further
$ M $ smooth.
Then the assertion follows, applying (6.13) inductively.

\medskip\noindent{{\bf 6.15.~Proposition.}} {\it
Let
$ \cV' $ be the full subcategory of
$ \cV(k) $ such that a variety
$ X $ belongs to
$ \cV' $ if and only if there is an open covering
$ \{U_{i}\} $ of
$ X $ with closed embeddings
$ U_{i} \to V_{i} $ such that
$ V_{i} \in \cV $.
Then a theory of mixed sheaves on
$ \cV $ is naturally extended to a theory of mixed sheaves on
$ \cV' $.
}

\medskip\noindent{\it Proof.}
For
$ X \in \cV' $, we consider the category
$ \cC(X) $ whose objects are closed embeddings
$ U \to V $ where
$ U $ is an open subvariety of
$ X $ and
$ V \in \cV $, and whose morphisms are pairs of morphisms of
$ U $ and of
$ V $, commutative with the closed embeddings.
Here the morphisms of
$ U $ are assumed to be natural morphisms (i.e., compatible with the natural
embedding into
$ X) $.
We say that
$ \{W_{i}\} $ is a covering family of
$ \cC(X) $ if
$ \{U_{i}\} $ is a covering of
$ X $ for
$ W_{i} = \{U_{i} \to V_{i}\} $.
We define
$ \cM(X) $ by the category whose object is
$ \{M_{W}\}_{W\in \cC(X)} $ with morphisms
$ u_{f} : f_{!}M_{W} \to M_{W'} $ for
$ f : W \to W' $ such that:

\noindent
(i)
$ M_{W} \in \cM_{U}(V) $ for
$ W = \{U \to V\} $,

\medskip\noindent
(ii) the restriction of
$ u_{f} $ to the complement of
$ U' \setminus U $ is an isomorphism, where
$ W = \{U \to V\}, W' = \{U' \to V'\} $,

\medskip\noindent
(iii)
$ u_{g}\scirc g_{!}u_{f} = u_{gf}\quad $ for
$ f : W \to W', g : W' \to W'' $.

\medskip\noindent
See also [33, 1.5].
For a morphism
$ f : X \to Y $ of
$ \cV $ inducing an isomorphism of their closed subvarieties
$ Z \simto Z' $, we have an equivalence of categories
$$
f_{!} : M_{Z}(X) \simto M_{Z'}(Y).
\leqno(6.15.1)
$$
This follows from (1.4.6) in the closed embedding case, and the general
case is checked by factorizing
$ f $ into the composition of a closed embedding and a smooth projection
locally, where the smooth projection has locally a section.
By a similar argument, we can check
$$
\bD\scirc f_{!} = f_{!}\scirc \bD,
\leqno(6.15.2)
$$
$$
W_{i}(f_{!}M_{W}) = f_{!}(W_{i}M_{W}).
\leqno(6.15.3)
$$
For
$ f : W \to W' $ and
$ \{M_{W}\} $ as above, we can show that
$ f_{!}M_{W} $ is independent of
$ f $, applying (6.15.1) to the projection of
$ V\times V' $ to
$ V $ (because the embedding by graph is its section).
By (6.15.1) and (2.2), we can also define
$ \cM(X) $ using a covering family
$ \{W_{i}\} $ of
$ \cC(X) $ (cf. [33, 1.6]).
In particular,
$ \cM(X) $ coincides with the previous
$ \cM(X) $ for
$ X \in \cV $.

The functor
$ \For $ is defined using [3].
The weight filtration is defined by
$ \{W_{i}M_{W}\} $ using (6.15.3), and the dual functor
$ \bD $ by
$ \{\bD M_{W}\} $ using (3.13.2).
We define the full subcategory
$ \cM(X,n) $ of pure objects with weight
$ n $ using
$ W $ as in (6.1).
Then we have the strict support decomposition (6.2.2), because it holds
locally, and the decomposition is unique so that the direct factors glue
globally by definition.
Then the semi-simplicity of pure objects is reduced to (6.4) in this
setting, where we use open direct images and pull-backs which will be
proved below.

The open pull-backs are defined naturally.
For an open embedding
$ j : X \to Y $ such that
$ j(X) = Y \setminus Z $ for a locally principal divisor
$ Z $, the functors
$ j_{!}, j_{*} $ are defined by extending the locally principal divisor to
$ V $ for
$ \{U \to V\} \in \cC(X) $, where the assertion is local
on
$ Y $ and we may assume
$ V $ affine and
$ Z \cap U $ is a principal divisor.
We can check the independentness of the choice of the extension of the
divisor to
$ V $, considering the union of the divisors and using the functor
$ \For $.
Then we can define the direct image as in (2.4), and the estimate of
weights is proved as in (6.7) using (6.7.5), where
$ i_{!}i^{!}M $ is defined as in (3.3) and the estimate of weights of
$ i_{!}i^{!}M $ follows from (6.7), because the assertion is local.

The external product is naturally defined using the commutativity with
direct image, because
$ W\times W' \in \cC(X\times Y) $ for
$ W \in \cC(X), W' \in \cC(Y) $.
We have (1.6.7--8), because the assertions are local by the same argument
as (4.8).
Then we can check the other conditions of \S 1.

\medskip\noindent
{\it Remark.} If
$ \cV $ is the full subcategory of smooth varieties, we have
$ \cV' = \cV(k) $.
So it is enough to define
$ \cM(X) $ satisfying the conditions in \S 1 for smooth varieties,
if one wants to construct
$ \cM(X) $ for
$ X \in \cV(k) $.

\bigskip\bigskip
\centerline{\bf \S 7. Geometric Origin}

\bigskip

Form now on, we assume
$ \cV = \cV(k) $.

\medskip\noindent
{\bf 7.1.}
Let
$ \{\cM(X)\} $ be a theory of
$ A $-mixed sheaves on
$ \cV(k) $ as above.
We define the categories
$ \cM(X)^{\go} $ by the smallest full subcategory of
$ \cM(X) $ satisfying the following conditions:

\noindent
(i)
$ \cM(X)^{\go} $ are stable by the cohomological functors
$ H^{i}f_{*}, H^{i}f_{!}, H^{i}f^{*}, H^{i}f^{!} $ for morphisms
$ f $,

\noindent
(ii)
$ \cM(X)^{\go} $ is stable by subquotients in
$ \cM(X) $,

\noindent
(iii)
$ A^{\cM} \in \cM(\Spec k)^{\go} $,

\noindent
where we assume also that an object of
$ \cM(X) $ isomorphic to an object of
$ \cM(X)^{\go} $ belongs to
$ \cM(X)^{\go} $.
We say that an object of
$ \cM(X)^{\go} $ is {\it of geometric origin} (cf. [3]).
By definition, an object of
$ \cM(X)^{\go} $ is obtained by iterating the above cohomological
functors and subquotients to
$ A^{\cM} $, and we get
$$
\cM(X)^{\go}\,\,\, \text{are stable by }\bD,\boxtimes
\leqno(7.1.1)
$$
by (2.4.2) (3.1.3) and (2.7) (3.1.9), where the external product is exact
for both factors by
$ (1.5.2) $.
By the next proposition,
$ \cM(X)^{\go} $ is an abelian full subcategory of
$ \cM(X) $ stable by
$ \Ker, \Coker $, and finite direct sum in
$ \cM(X) $, and
$ \cM(X)^{\go} $ form a theory of mixed sheaves on
$ \cV(k) $ so that the functors are compatible with the natural functor
$ \cM(X)^{\go} \to \cM(X) $.

\medskip\noindent{{\bf 7.2.~Proposition.}} {\it
The subcategory
$ \cM(X)^{\go} $ is stable by finite direct sum in
$ \cM(X) $.
For
$ M \in \cM(X) $, we have
$ M \in \cM(X)^{\go} $ if and only if, for any point
$ x $ of
$ X $, there exist an open subvariety
$ U $ of
$ X $ containing
$ x $, a closed embedding
$ i : U \to Z $, a quasi-projective morphism
$ \pi : Y \to Z $, a divisor
$ D $ on
$ Y $ with natural inclusions
$ j_{U} : U \to X $,
$ j : Y \setminus D \to Y $, such that
$ Y $ is smooth,
$ D $ is a reduced divisor with normal crossings, and
$ i_{*}(j_{U})^{*}M $ is isomorphic to a subquotient of
$$
H^{m}\pi_{*}j_{!}A_{Y\setminus D}^{\cM}(n)
\leqno(7.2.1)
$$
for
$ n, m \in \bZ $.
}

\medskip\noindent{\it Proof.}
We first note that we may assume
$ U, Z $ affine and
$$
Z\,\,\, \text{is an affine space }\bA^{r}
\leqno(7.2.2)
$$
by replacing
$ Z $ with an affine space containing
$ Z $.

Let
$ \cM(X)' $ denote the full subcategory of
$ \cM(X) $ defined by the above condition.
We first show
$$
\cM(X)' \subset \cM(X)^{\go}.
\leqno(7.2.3)
$$
We have an affine open covering
$ \{U_{k}\} $ of
$ X $ with
$ \pi_{k} : Y_{k} \to Z_{k} $,
$ j_{k} : Y_{k} \setminus D_{k}\to Y_{k} $,
$ i_{k} : U_{k} \to Z_{k} $ satisfying (7.2.1) as above.
By (4.2.1) (4.12) (5.7.7), we have
$$
\aligned
H^{0}(a_{\bP^{1}})_{*}A_{{P}^{1}}^{\cM}
&= A^{\cM},\,\,\,H^{2}(a_{\bP^{1}})_{*}A_{{P}^{1}}^{\cM}
=A^{\cM}(-1)
\\
H^{0}(a_{S^{*}})_{*}A_{{S}^{*}}^{\cM}
&= A^{\cM},\,\,\,H^{1}(a_{S^{*}})_{*}A_{{S}^{*}}^{\cM}
= A^{\cM}(-1).
\endaligned
\leqno(7.2.4)
$$
So we may assume
$ m, n $ in (7.2.1) is independent of
$ k $ by (2.7.3) (2.7.5), replacing
$ Y_{k} $ by its product with
$ (\bP^{1})^{b_{k}}\times (S^{*})^{c_{k}} $.
Let
$ X', Y', D', Z' $ be the disjoint union of
$ U_{k}, Y_{k}, D_{k}, Z_{k} $ respectively with natural morphisms
$ \pi ' : Y' \to Z' $,
$ j' : Y' \setminus D' \to Y' $,
$ i' : X' \to Z' $,
$ j'' : X' \to X $.
Then there is a subquotient
$ M' $ of
$ H^{m}\pi '_{*}j'_{!}A_{Y'\setminus D'}^{\cM}(n) $ such that
$ \supp M' \subset \Im \,i' $ and
$ M $ is isomorphic to a subobject of
$ j''_{*}H^{0}i'^{*}M' $ by (2.2).
So we get (7.2.3).

By a similar argument, we can check
$$
\cM(X)'\,\,\, \text{is stable by finite direct sum in }
\cM(X).
\leqno(7.2.5)
$$
In fact, for
$ M_{1}, M_{2} \in \cM(X)' $ and a point of
$ X $, there are
$ \pi_{a} : Y_{a} \to Z_{a} $,
$ j_{a} : Y_{a} \setminus D_{a}\to Y_{a} $,
$ i_{a} : U \to Z_{a} $ associated with
$ M_{a} $ for
$ a = 1, 2 $.
We have a section
$ i'_{a} $ of the projection of
$ Z_{1}\times Z_{2} $ onto
$ Z_{a} $ such that
$ i'_{1}i_{1} = i'_{2}i_{2} $, assuming (7.2.2).
So we may assume
$ Z_{1} = Z_{2}, i_{1} = i_{2} $, by replacing
$ Z_{a} $ with
$ Z_{1}\times Z_{2} $, and
$ i_{a} $ by
$ i'_{a}i_{a} $.
Then (7.2.5) is clear, because we may assume
$ m_{1} = m_{2}, n_{1} = n_{2} $ using (7.2.4).

For the converse of (7.2.3), it is enough to show that
$ \cM(X)' $ are stable by the cohomological functors in the definition
of
$ \cM(X)^{\go} $, because
$ \cM(X)' $ is clearly stable by subquotient.
Here we may assume
$ Y $ in (7.2.1) is pure dimensional, using the product of each connected
component with an affine space of appropriate dimension.
By (2.1.4) (4.3.1), we have
$$
\bD(H^{m}\pi_{*}j_{!}A_{Y\setminus D}^{\cM}(n)) = H^{d-m}\pi_{!}j_{*}
A_{Y\setminus D}^{\cM}(2d-n)
\leqno(7.2.6)
$$
where
$ d = \dim Y $.
Let
$ \overline{\pi} : \overline{Y} \to Z $ be a projective morphism with
an open embedding
$ \overline{j} : Y \to \overline{Y} $ such that
$ \overline{Y} $ is smooth and
$ (\overline{Y} \setminus Y) \cup D $ is a divisor with normal crossings
[17].
Let
$ Y' $ be the complement of the closure of
$ D $ in
$ \overline{Y}, D' = Y' \setminus (Y \setminus D) $ with natural morphisms
$ \pi ' : Y' \to Z, j' : Y' \setminus D' \to Y' $,
$ \overline{j}': Y' \to \overline{Y} $.
Then
$ Y' \setminus D' = Y \setminus D $, and
$$
\overline{j}_{!}j_{*}A_{Y\setminus D}^{\cM}[d] = \overline{j}'_{*}j'_{!}
A_{Y'\setminus D'}^{\cM}[d].
\leqno(7.2.7)
$$
In fact, their restrictions to
$ Y $ coincide, and we have the morphism (7.2.7) by
$ \overline{j}_{!}\overline{j}^{!} \to id $.
This is an isomorphism, considering their underlying perverse sheaves.
By (2.6) (4.12), we have
$$
H^{d-m}\pi_{!}j_{*}A_{Y\setminus D}^{\cM}(2d-n) = H^{d-m}\pi '_{*}j'_{!}
A_{Y'\setminus D'}^{\cM}(2d-n).
\leqno(7.2.8)
$$
So (7.2.1) is essentially self-dual, and it is enough to show the stability
by
$ H^{i}f_{*} $ and
$ H^{i}f^{!} $.

We first show the stability by
$ H^{i}j_{*}j^{*}, H^{i}j_{!}j^{!} $ for an open inclusion
$ j : X' \to X $.
Since the assertion is local, we may assume
$ X $ is affine,
$ C := X \setminus X' $ is an intersection of principal divisors, and
$ M $ is a subquotient of (7.2.1) with
$ U = X $.
By the same argument as (2.4), the assertion is reduced to the case where
$ C $ is a principal divisor using a
Cech complex associated with an affine open covering of
$ X $.
Then
$ C $ is extended to a principal divisor
$ C' $ on
$ Z $, because
$ Z $ is affine (cf. (7.2.2)).
So it is enough to show the assertion for the open inclusion
$ Z \setminus C' \to Z $,
because the functors
$ j_{*}j^{*}, j_{!}j^{!} $ are exact and commute with subquotients.
Then the assertion is clear by base change.

We now consider the direct image
$ H^{i}f_{*} $.
By the above arguments, the resolution in (2.4) is defined in
$ \cM(X)' $, and we may assume
$ X $ affine and
$ H^{i}f_{*}M = 0 $ for
$ i \ne 0 $, replacing
$ X $ by
$ U_{I} $ in (2.4).
Since the functor
$ H^{0}f_{*} $ is right exact, we may replace
$ M $ with
$ \widetilde{j}_{!}\widetilde{j}^{!}M $,
where
$ \widetilde{j} : X \setminus H \to X $ is a natural inclusion with
$ H $ a sufficiently general hyperplane section of
$ X $ (here
$ X $ is viewed as a quasiprojective variety).
Using again a sufficiently fine covering
$ \{U_{i}\} $ of
$ X $ and replacing
$ X $ with
$ U_{I} $,
we may assume that the condition (7.2.1) is satisfied with
$ X = U $ and
$ f $ is extended to a morphism of
$ Z $ to an affine space containing an open subvariety of
$ Y $,
where the hyperplane section
$ H $ is taken generically with respect to
$ U_{i} $ (cf. [2]).
So the assertion is reduced to the case
$ Z = X $.
Furthermore we may assume that
$ H^{m'}\pi_{*}j_{!}A_{Y\setminus D}^{\cM}(n) (m ' \in \bZ) $ and the
subobjects of
$ H^{m}\pi_{*}j_{!}A_{Y\setminus D}^{\cM}(n) $ defining
$ M $ are
$ f_{*} $-acyclic, replacing them with the image of the functor
$ \widetilde{j}_{!}\widetilde{j}^{!} $.
Then the assertion follows from the Leray spectral sequence.
The argument is similar for
$ H^{i}f_{!} $.

It remains to show the stability by pull-back.
Since the assertion is local, and is clear for a smooth projection, we may
assume that
$ f $ is a closed embedding.
By definition (3.3), the assertion follows from the stability by open direct
image and finite direct sum.

\medskip\noindent
{\it Remark.} We can show that the condition in (7.2) is satisfied for any
affine open subvariety
$ U $ of
$ X $, using a similar argument.

\medskip\noindent
{\bf 7.3.~Definition.} Let
$ \cM(X,w)^{\go}, \cM_{Z}(X,w)^{\go} $ be the full subcategories of
$ \cM(X)^{\go} $ defined as in (6.1), where
$ \cM_{Z}(Z,w)^{\go} = \cM_{Z}(X,w)^{\go} $ by definition and (1.4.6).
We say that
$ M \in \cM_{Z}(Z,w)^{\go} $ has {\it geometric level}
$ \le n $, if
$ M $ is isomorphic to a direct factor of
$ H^{m}\pi_{*}A_{Y}^{\cM}(r) $ for a projective morphism
$ \pi $ of a smooth variety
$ Y $ onto
$ Z $ such that
$ \dim Y \le n $.
We say that
$ M \in \cM(X)^{\go} $ has {\it geometric level}
$ \le n $ if so do the direct factors of
$ \Gr_{i}^{W}M $ with strict support, and
$ M \in \cM(X)^{\go} $ has geometric level
$ n $ if it has geometric level
$ \le n $ but not
$ \le n - 1 $.
We denote by
$ \cM(X{)}_{\gl \le n}^{\go} $ the full subcategory of
$ \cM(X)^{\go} $ consisting of objects with geometric level
$ \le n $, and let
$ \cM(X,w{)}_{\gl\le n}^{\go} = \cM(X,w)\cap\cM(X{)}_{\gl\le n}^{\go} $.

\medskip\noindent
{\it Remarks.} (i) The full subcategory
$ \cM(X{)}_{\gl \le n}^{\go} $ is stable by the pull-back under an open
embedding.
Conversely, for
$ M \in \cM_{Z}(Z,w)^{\go} $, we have
$ M \in \cM(Z{)}_{\gl \le n}^{\go} $ if
$ j^{*}M \in \cM(U{)}_{\gl \le n}^{\go} $ for some non empty open
subvariety
$ U $ of
$ Z $ with
$ j : U \to Z $ the natural inclusion.
In fact, a projective morphism onto an open subvariety of
$ Z $ is extended to a projective morphism onto
$ Z $ by Nagata-Hironaka (where we may assume
$ Y $ affine), and
$ M $ must be isomorphic to a direct factor of the direct image of the
constant sheaf by (6.4) (6.7).

\medskip

(ii) Let
$ n $ be the geometric level of
$ M \in \cM(X,w)^{\go} $.
Then
$$
n \equiv w\mod 2,
\leqno(7.3.1)
$$
because we may assume
$ Y $ pure dimensional and
$$
m = \dim Y,
\leqno(7.3.2)
$$
in the definition of geometric
$ \le n $, where
$ m - 2r = w $ by (6.7).
In fact, a connected component
$ Y_{i} $ of
$ Y $ may be replaced by its product with
$ \bP^{k} $ with
$ k = \dim Y - \dim Y_{i} $ using (2.7.3) (2.7.5).
If
$ m > \dim Y $, we may may replace
$ m $ by
$ 2 \dim Y - m $ using (7.4.1) below.
If
$ m < \dim Y $, we can replace
$ Y $ by its intersection with a generic hyperplane by the weak Lefschetz
property of a generic fiber of
$ f $, where we can restrict
$ Y $ to its open subvariety by (i) so that we may assume
$ f $ is smooth and
$ X $ is a closed subvariety of
$ \bP^{N}\times Y $.
So we get (7.3.2).

\medskip

(iii) In the definition of geometric
$ \le n $, we can replace
$ Y $ by a smooth variety
$ Y' $ with a birational proper morphism
$ \pi : Y' \to Y $.
In fact,
$ A_{Y}^{\cM} $ is a direct factor of
$ \pi_{*}A_{Y'}^{\cM} $ by the decomposition (6.10.1) for
$ \pi_{*}A_{Y'}^{\cM} $ and the strict support decomposition (6.2.2),
because the direct factor of
$ H^{0}\pi_{*}(A_{Y'}^{\cM}[\dim Y]) $ with strict support
$ Y $ is
$ A_{Y}^{\cM}[\dim Y] $ by (6.4).
Here we may assume
$ Y $ connected.

\medskip

(iv) By (7.3.2),
$ \cM(X{)}_{\gl \le n}^{\go} $ is stable by finite direct sum.
By semi-simplicity of pure objects,
$ \cM(X{)}_{\gl \le n}^{\go} $ is an abelian full subcategory of
$ \cM(X)^{\go} $ stable by subquotient and extension in
$ \cM(X)^{\go} $.
It is also stable by the dual functor
$ \bD $ using the duality (4.3) and (4.12).

\medskip\noindent{{\bf 7.4.~Proposition.}} {\it
Let
$ f : X \to Y $ be a projective morphism of varieties, and
$ l\in H^{2}(X(\bC),A(1)) $ the Chern class of a relatively ample line
bundle.
Then, for
$ M \in \cM(X,w)^{\go} $, we have the relative hard Lefschetz:
$$
l^{k} : H^{-k}f_{*}M \simto H^{k}f_{*}M(k)\quad \text{for}\quad
k > 0,
\leqno(7.4.1)
$$
whose underlying morphism is induced by the natural action of
$ l $ on
$ A $-complexes.
}

\medskip\noindent{\it Proof.}
It is enough to define the action of
$ l $ on
$ H^{\ssbull}f_{*}M $, because the isomorphism on the underlying
perverse sheaves follows from [27] [28] using (7.2).
So the assertion is local, we may assume
$ f $ is the projection of
$ \bP^{N}\times Y $ to
$ Y $.
Then the assertion is proved as in [3].
In fact, let
$ P^{\vee} $ be the dual projective space of
$ P := \bP^{N} $ with the universal hyperplane
$ H $ of
$ P\times P^{\vee} $.
Then the action of
$ l $ is defined on
$ (H^{\ssbull}f_{*}M)\,\boxtimes\,A_{{P}^{\vee}}^{\cM} $ by the
restriction and Gysin morphism using (5.12), and it induces the action on
$ H^{\ssbull}f_{*}M $.

\medskip\noindent{{\bf 7.5.~Proposition.}} {\it
Let
$ j : U \to X $ be an affine open embedding of varieties.
Then the direct images
$ j_{*}, j_{!} $ induce
}
$$
j_{*},\,\,\,j_{!} : \cM(U{)}_{\gl \le n}^{\go}\to
\cM(X{)}_{\gl \le n}^{\go}.
\leqno(7.5.1)
$$

\medskip\noindent
{\it Proof.} Since the functors are exact, it is enough to consider the
direct image of a pure object
$ M $ of
$ \cM(U{)}_{\gl \le n}^{\go} $ with strict support.
Let
$ \pi : Y \to U $ be as in (7.3) so that
$ M $ is a direct factor of
$ H^{m}\pi_{*}A_{Y}^{\cM}(r) $.
Let
$ \pi ' : Y' \to X $ be its extension such that
$ \pi ' $ is projective,
$ Y' $ is smooth and
$ D := Y' \setminus Y $ is a divisor with normal crossings whose
irreducible components are smooth.
By duality, the assertion is reduced to the following:

\medskip\noindent{{\bf 7.6.~Proposition.}} {\it
Let
$ X $ be a smooth variety with pure dimension
$ d $, and
$ D $ a reduced divisor with normal crossings whose irreducible
components
$ D_{i} $ are smooth.
Let
$ U = X \setminus D $ and
$ D_{I} = \bigcap_{i\in I} D_{i} $ with natural inclusions
$ j : U \to X $ and
$ a_{I} : D_{I} \to X $, where
$ D_{\emptyset} = X $.
Then we have an isomorphism
$$
\Gr_{k}^{W}(j_{!}A_{U}^{\cM}[d]) = \mopls_{|I|=d-k} (a_{I})_{*}
A_{{D}_{I}}^{\cM}[k]\quad \text{in}\quad
\cM(X),
\leqno(7.6.1)
$$
where
$ W $ is the weight filtration of
$ j_{!}A_{U}^{\cM}[d] \in \cM(X) $.
}

\medskip\noindent{\it Proof.}
For
$ J \supset I $, we have the restriction morphism
$$
(a_{I})_{*}A_{{D}_{I}}^{\cM} \to (a_{J})_{*}A_{{D}_{J}}^{\cM}
\leqno(7.6.2)
$$
compatible with composition by (3.9.2).
This induces
$$
d^{p} : \mopls_{|I|=p} (a_{I})_{*}A_{{D}_{I}}^{\cM} \to
\mopls_{|I|=p+1} (a_{I})_{*}A_{{D}_{I}}^{\cM}
\leqno(7.6.3)
$$
such that
$ d^{p+1}\scirc d^{p} = 0 $ by defining
$ d^{p} $ with a sign like a
Cech complex.
We have
$$
\Ext^{k}((a_{I})_{*}A_{{D}_{I}}^{\cM}, (a_{J})_{*}A_{{D}_{J}}^{\cM})
= 0\quad \text{for}\quad k < 0,
\leqno(7.6.4)
$$
by adjunction, because
$$
i^{!}A_{Y}^{\cM} = A_{X}^{\cM}(-r)[-2r]
\leqno(7.6.5)
$$
for a closed embedding of smooth varieties
$ i : X \to Y $ with codimension
$ r $ by (4.3).
We define inductively
$ M'_{p} \in D^{b}\cM(X) $ by
$$
M'_{p} = C(u_{p} : M'_{p-1}[-1] \to \mopls_{|I|=p} (a_{I})_{*}
A_{{D}_{I}}^{\cM}[d-p])
\leqno(7.6.6)
$$
where the morphism
$ u_{p+1} $ is induced by (7.6.3--4).
Then we get inductively
$ M'_{p} \in \cM(X) $, and they determine a filtration
$ G $ of
$ M'_{d} $ such that
$$
\Gr_{k}^{G}M'_{d} = \mopls_{|I|=d-k} (a_{I})_{*}A_{{D}_{I}}^{\cM}[k].
\leqno(7.6.7)
$$
Applying the same construction to the underlying
$ A $-complexes, we can check
$$
\For(M'_{d}) = j_{!}A_{U}[d]
\leqno(7.6.8)
$$
using (7.6.4).
Since
$ j^{!}M'_{d} = A_{U}^{\cM} $, we have a natural morphism
$$
M'_{d} \to j_{!}A_{U}^{\cM}
\leqno(7.6.9)
$$
which is an isomorphism by (7.6.8).
Since
$ \Gr_{k}^{G}M'_{d} $ is pure of weight
$ k $, we have
$ G = W $, and the assertion follows.

\medskip\noindent{{\bf 7.7.~Corollary}} {\it
With the above notation and assumption, assume further
$ D = g^{-1}(0)_{\red} $ for a function
$ g $.
Then
$$
P_{N}\Gr_{d-1+k}^{W}(\psi_{g,1}A_{X}^{\cM}[d]) = \mopls_{|I|=k+1}
(a_{I})_{*}A_{{D}_{I}}^{\cM}(-k)[d-k
-1]\,\,\, \text{for }k \ge 0
\leqno(7.7.1)
$$
where the left hand side is the primitive part of
$ \Gr_{d-1+k}^{W}(\psi_{g,1}A_{X}^{\cM}[d]) $ defined by the kernel of
$ N^{k+1} $.
}

\medskip\noindent{\it Proof.}
Let
$ i : D \to X $ denote the natural inclusion.
Since
$ A_D^{\cM} = i^{*}A_{X}^{\cM} $, and
$ \For(A_D^{\cM}[d-1]) $ is a perverse sheaf, we have an isomorphism
$$
A_D^{\cM}[d-1] = \Ker\,N \subset \psi_{g,1}A_{X}^{\cM}[d]
\leqno(7.7.2)
$$
by (5.9.1), because
$ \Var $ is injective in this case (cf. [27]).
Since
$ N $ is strictly compatible with
$ W $, we get
$$
\Gr_{k}^{W}(A_D^{\cM}[d-1]) = \Ker\,N \subset \Gr_{k}^{W}
\psi_{g,1}A_{X}^{\cM}[d].
\leqno(7.7.3)
$$
On the other hand, we have
$$
\Gr_{k}^{W}(A_D^{\cM}[d-1]) = \mopls_{|I|=d-k} (a_{I})_{*}
A_{{D}_{I}}^{\cM}[k]\quad \text{for}\quad
k \le d - 1
\leqno(7.7.4)
$$
by (7.6) using the short exact sequence in
$ \cM(X) $:
$$
0 \to A_{X}^{\cM}[d] \to j_{!}A_{U}^{\cM}[d] \to A_D^{\cM}
[d-1] \to 0,
\leqno(7.7.5)
$$
cf. (3.3.2).
Since
$ A_{X}^{\cM}[d] $ is pure of weight
$ d $ by (6.6), the weight filtration on
$ \psi_{g}A_{X}^{\cM}[d] $ is the monodromy filtration shifted by
$ d - 1 $ (cf. (6.11)).
So we have the primitive decomposition:
$$
\Gr_{k}^{W}\psi_{g}A_{X}^{\cM}[d] = \mopls_{i\ge 0} N^{i}(P_{N}
\Gr_{k+2i}^{W}\psi_{g}A_{X}^{\cM}[d])(i)
\leqno(7.7.6)
$$
where
$ P_{N}\Gr_{k}^{W} = 0 $ for
$ k > d - 1 $.
Then
$$
N^{k} : P_{N}\Gr_{d-1+k}^{W}\psi_{g}A_{X}^{\cM}[d] \simto \Ker\,N
\subset P_{N}\Gr_{d-1-k}^{W}
\psi_{g}A_{X}^{\cM}[d],
\leqno(7.7.7)
$$
and the assertion follows from (7.7.2--3).

\medskip\noindent{{\bf 7.8.~Proposition.}} {\it
Let
$ g $ be a function on
$ X $.
Then
}
$ \psi_{g}, \varphi_{g,1}, \xi_{g} $ induce
$$
\psi_{g}, \varphi_{g,1}, \xi_{g} : \cM(X{)}_{\gl \le n}^{\go} \to
\cM(X{)}_{\gl \le n}^{\go}.
\leqno(7.8.1)
$$

\medskip\noindent
{\it Proof.} Let
$ X^{*} = X \setminus g^{-1}(0) $.
We first reduce the assertion to:
$$
\psi_{g,1}M \in \cM(X{)}_{\gl \le n}^{\go}\,\,\, \text{for }M
\in \cM_{Z}(X^{*},w{)}_{\gl \le n}^{\go}.
\leqno(7.8.2)
$$
Since
$ \psi_{g}, \varphi_{g,1} $ are exact,we may assume
$ M $ pure, and then
$ M \in \cM_{Z}(X,w{)}_{\gl \le n}^{\go} $ by the strict support
decomposition.
Then
$ \can : \psi_{g,1}M \to \varphi_{g,1}M $ is surjective for
$ M \in \cM_{Z}(X,w{)}_{\gl \le n}^{\go} $ such that
$ Z \not\subset g^{-1}(0) $ (cf. [27, 5.1.4]), and the assertion on
$ \varphi_{g,1} $ is reduced to (7.8.2).
The assertion on
$ \xi_{g} $ follows from (5.5.1).
So it remains to reduce the assertion on
$ \psi_{g} $ to (7.8.2).

By definition, we may replace
$ g : X \to S $ by the projection
$ p : Y = X\times S \to S $ in (5.2), because the stability by
$ (i_{g})_{*} $ is clear.
Let
$ Y' \to S $ be the base change of
$ p $ by the
$ n $-fold ramified covering
$ \pi : S \to S $ as in the proof of (5.7), and
$ \pi ' : Y' \to Y $ the natural morphism.
Then we can check
$$
\pi '_{*}A_{Y'}^{\cM}\otimes M = \pi '_{*}\pi '^{*}M
$$
using the diagram
$$
\CD
Y' @= Y' @>>> Y
\\
@VVV @VVV @VVV
\\
Y'\times Y' @>>> Y'\times Y @>>> Y\times Y
\endCD
$$
and (3.10.2) (4.1.5).
Similarly, we have
$$
\pi '_{*}A_{Y'}^{\cM} = p^{*}\pi_{*}A_{S}^{\cM}.
$$
Since
$ \cM(Y{)}_{\gl \le n}^{\go} $ is stable by the pull-back under an \'etale
morphism, the assertion is reduced to that for
$ \psi_{p\pi ',1} $ by definition (5.7.6) and (5.10.1).

Now we show (7.8.2).
By definition (7.3),
$ M $ is a direct factor of
$ H^{m}\pi_{*}A_{Y}^{\cM}(r) $, where we may assume
$ (g\pi )^{-1}(0)_{\red} $ is a divisor with normal crossings by
Remark (iii) after (7.3) and [17].
Then the assertion follows from (7.7) and (5.10.1).

\medskip\noindent{{\bf 7.9.~Proposition.}} {\it
Let
$ i : X \to Y $ be a closed embedding of varieties.
Then the direct image induces equivalences of categories
$$
i_{*} : \cM(X{)}_{\gl \le n}^{\go} \simto \cM_{X}(Y{)}_{\gl \le n}^{\go}
\leqno(7.9.1)
$$
$$
i_{*} : D^{b}\cM(X{)}_{\gl \le n}^{\go} \simto {D}_{X}^{b}\cM(Y
{)}_{\gl \le n}^{\go}
\leqno(7.9.2)
$$
where
$ \cM_{X}(Y{)}_{\gl \le n}^{\go} $,
$ {D}_{X}^{b}\cM(Y{)}_{\gl \le n}^{\go} $ are respectively
the full subcategory of
$ \cM_{X}(Y{)}_{\gl \le n}^{\go} $,
$ D^{b}\cM(Y{)}_{\gl \le n}^{\go} $ consisting of
the objects whose \(cohomological\,\) support is contained in
$ i(X) $.
}

\medskip\noindent{\it Proof.}
The first assertion is clear by definition using (1.4.6) (3.7.3).
Then the second follows from the same argument as (5.6) using (7.5) (7.8).

\medskip\noindent{{\bf 7.10.~Proposition.}} {\it
Let
$ f : X \to Y $ be a morphism of varieties.
Then we have naturally the direct image functors
$$
f_{*}, f_{!} : D^{b}\cM(X{)}_{\gl \le n}^{\go} \to D^{b}\cM
(Y{)}_{\gl \le n}^{\go}
\leqno(7.10.1)
$$
compatible with the composition of morphisms of varieties.
Moreover, we have a canonical isomorphism
$ f_{!}\scirc \bD = \bD\scirc f_{*} $ and a canonical morphism
$ f_{!} \to f_{*} $.
These functors and morphisms are compatible with the direct images in
{\rm (2.4)} and the corresponding morphisms {\rm (2.4.2) (4.12.1)}
by the natural functors
$ D^{b}\cM(X{)}_{\gl \le n}^{\go} \to D^{b}\cM(X) $.
}

\medskip\noindent{\it Proof.}
For (7.10.1), it is enough to show
$$
H^{i}f_{*}M, H^{i}f_{!}M \in \cM(Y{)}_{\gl \le n}^{\go}\quad\text{for}\quad
M \in \cM(X{)}_{\gl \le n}^{\go}
\leqno(7.10.2)
$$
for an affine morphism
$ f $, using (7.5) and the same argument as in the proof of (2.4).
Using a compactification
$ f' : X' \to Y $ of
$ f $ such that
$ X = X' \setminus C $ for a locally principal divisor
$ C $ on
$ X' $, we may assume
$ f $ projective by (7.5).
We may also assume
$ M $ pure with strict support by the spectral sequence associated with
the weight filtration.
Then the assertion is clear, using the decomposition (6.10) for the direct
image in the definition of geometric level in (7.3).

The isomorphism
$ f_{!}\scirc \bD = \bD\scirc f_{*} $ is clear by definition.
For
$ f_{!} \to f_{*} $, it is enough to show
$ f_{!} = f_{*} $ for
$ f $ proper, because we have
$ j_{!} \to j_{*} $ for an open embedding
$ j : X \to X' $ such that
$ j(X) = X' \setminus C $ for a locally principal divisor
$ C $ on
$ X' $.
By the above construction,
$ M \in D^{b}\cM(X{)}_{\gl \le n}^{\go} $ is represented by a complex
whose components
$ M^{i} $ satisfy
$ H^{k}f_{*}M^{i} = 0 $ for
$ k \ne 0 $, and we have
$ H^{k}f_{!}M^{i} = H^{k}f_{*}M^{i} $ for any
$ k $ by (4.12).
So we get
$ f_{!}M = f_{*}M $, where we can check the independence of the
representative as in (2.4).

\medskip\noindent{{\bf 7.11.~Proposition.}} {\it
Let
$ i : X \to Y $ be a closed embedding of varieties, and
$ U = Y \setminus X $ with
$ j : U \to Y $.
Then the pull-back functors induce functors
$$
i^{!}, i^{*} : D^{b}\cM(Y{)}_{\gl \le n}^{\go} \to D^{b}\cM
(X{)}_{\gl \le n}^{\go}
\leqno(7.11.1)
$$
with the triangles of functors {\rm (3.3.2)}, such that they are
compatible with those in {\rm (3.7)} by the natural functors
$ D^{b}\cM(X{)}_{\gl \le n}^{\go} \to D^{b}\cM(X) $.
Moreover,
$ j_{!}, j_{*} $ are the adjoint functors of
$ j^{!} = j^{*} $.
}

\medskip\noindent{\it Proof.}
This follows from the same argument as (3.2) (3.7) using (7.5) (7.9).

As a corollary, we get

\medskip\noindent{{\bf 7.12.~Proposition.}} {\it
Let
$ \cA = \cM(X), \cM(X)^{\go} $ or
$ \cM(X{)}_{\gl \le n}^{\go} $.
Then we have a `classical'
$ t $-structure
$ ({}^{c}\cD^{\le 0},{}^{c}\cD^{\ge 0}) $ on
$ D^{b}\cA $ such that
$ M \in{}^{c}\cD^{\le 0} $ \(resp.
$ M \in{}^{c}\cD^{\ge 0}) $ if and only if
$ \cH^{j} \For(M) = 0 $ for
$ j > 0 $ \(resp.
$ j < 0) $, where
$ \cH^{j} $ is the natural \(i.e.
classical\,\) cohomology functor of
$ {D}_{c}^{b}(X(\bC),A) $.
Moreover, we have
$$
\For : \cC \to M(A_{X(\bC)})\,\,\,
\text{is exact and faithful,}
\leqno(7.12.1)
$$
where
$ \cC $ is the heart of the
$ t $-structure.
}

\medskip\noindent{\it Proof.}
Using the theory of gluing
$ t $-structure in [3], we can define
$ {}^{c}\cD^{\le 0} $ (resp.
$ {}^{c}\cD^{\ge 0}) $ by the following condition:

\medskip\noindent
For any closed embedding
$ i_{S} : S \to X $ of an irreducible variety
$ S $, there is a non empty open subvariety
$ U $ of
$ S $ with natural inclusion
$ j_{U} : U \to S $ such that
$ (j_{U})^{*}H^{k}{i}_{S}^{*}M = 0 $ for
$ k > \dim S $ (resp.
$ (j_{U})^{*}H^{k}{i}_{S}^{!}M = 0 $ for
$ k < \dim S) $.

\medskip\noindent
This condition depends only on the underlying
$ A $-complexes, and we can check the coincidence with the condition in
Proposition using the distinguished triangles in (3.3.2).
For the condition of
$ {}^{c}\cD^{\ge 0} $, we use also the left exactness of the direct
image
$ j_{*} $ with respect to the classical
$ t $-structure on
$ {D}_{c}^{b}(X(\bC),A) $ for an open embedding
$ j $.

Let
$ {}^{c}H^{k} : D^{b}\cA \to \cC $ be the cohomology functor
associated with the
$ t $-structure, so that
$ M \in D^{b}\cA $ belongs to
$ \cC $ if and only if
$ {}^{c}H^{k}(M) = 0 $ for
$ k \ne 0 $ (cf. [3]).
The functor
$ {}^{c}H^{k} $ corresponds to the natural cohomology functor
$ \cH^{k} : {D}_{c}^{b}(X(\bC),A) \to M(A_{X(\bC)}) $ by the functor
$ \For $, and
$ \For : \cC \to M(A_{X(\bC)}) $ is exact using the associated
long exact sequence.
It is also faithful, because
$ \Im $ commutes with
$ \For $, and
$ M \in \cC $ is
$ 0 $ if and only if
$ \For(M) = 0 $.
So we get the last assertion.

\medskip\noindent{{\bf 7.13.~Proposition.}} {\it
Let
$ X $ be a variety of dimension
$ \le n $.
Then
$ A_{X}^{\cM} \in D^{b}\cM(X) $ is naturally lifted to
$ D^{b}\cM(X{)}_{\gl \le n}^{\go} $.
More precisely, there exists uniquely an object of
$ D^{b}\cM(X{)}_{\gl \le n}^{\go} $, denoted also by
$ A_{X}^{\cM} $, such that
$ \For(A_{X}^{\cM}) = A_{X(\bC)} $ and the restriction to a smooth open
dense subvariety
$ U $ is isomorphic to
$ A_{U}^{\cM} $ in
$ D^{b}\cM(U) $.
Moreover, for
$ f : X \to Y $ a morphism of varieties of dimension
$ \le n $, the restriction and Gysin morphisms
$ A_{Y} \to f_{*}A_{X(\bC)} $ and
$ f_{!}\bD A_{X(\bC)} \to \bD A_{Y(\bC)} $ are uniquely lifted to
morphisms of
$ D^{b}\cM(Y{)}_{\gl \le n}^{\go} $, and they are compatible with the
composition of morphisms.
}

\medskip\noindent{\it Proof.}
The first assertion is clear if
$ X $ smooth.
We first show the uniqueness of
$ A_{X}^{\cM} $ in
$ D^{b}\cM(X{)}_{\gl \le n}^{\go} $.
Let
$ M, M' $ be objects of
$ D^{b}\cM(X{)}_{\gl \le n}^{\go} $ satisfying the condition of
$ A_{X}^{\cM} $.
Then the functor
$ \For $ induces an injective morphism
$$
\Hom_{gl\le n}(M, M') \to \Hom(A_{X(\bC)}, A_{X(\bC)})
\leqno(7.13.1)
$$
by (7.12.1), where
$ \Hom_{gl\le n} $ is taken in
$ D^{b}\cM(X{)}_{\gl \le n}^{\go} $.
Let
$ U $ be as in Proposition, and
$ j : U \to X $ the inclusion.
We have an injective morphism
$$
M' \to{}^{c}H^{0}j_{*}A_{U}^{\cM}\,\,\, \text{in}\quad \cC
\leqno(7.13.2)
$$
by the long exact sequence associated with a triangle in (3.3.2).
By (7.13.1), we have the injectivity of
$$
\Hom_{gl\le n}(M, M') \to \Hom_{gl\le n}(M,{}^{c}H^{0}
j_{*}A_{U}^{\cM})=\Hom_{gl\le n}(A_{U}^{\cM}, A_{U}^{\cM})
\leqno(7.13.3)
$$
because it is clear for the underlying
$ A $-complexes and
$ j^{*}M = j^{*}M' = A_{U}^{\cM} $.
Here the last isomorphism follows from the adjoint relation, because we may
replace
$ {}^{c}H^{0}j_{*}A_{U}^{\cM} $ with
$ j_{*}A_{U}^{\cM} $ using
$ {}^{c}H^{k}j_{*}A_{U}^{\cM} = 0 $ for
$ k < 0 $.
So it is enough to show that the canonical morphism
$ u : M \to{}^{c}H^{0}j_{*}A_{U}^{\cM} $ belongs to the image of the first
morphism of (7.13.3), i.e., the image of
$ u $ in
$ \Hom_{gl\le n}(M,{}^{c}H^{0}j_{*}A_{U}^{\cM}/M') $ is zero.
But it is true for the underlying
$ A $-complexes, and the assertion follows from (7.12.1).

For existence, we extend
$ A_{U}^{\cM} $ inductively to larger open subvarieties of
$ X $ using noetherian induction.
Assume the assertion is proved for an open subvariety
$ U $ which is not necessary smooth.
Applying (7.12) to
$ \cA = \cM(X) $, we have a natural morphism
$$
A_{X}^{\cM} \to{}^{c}H^{0}j_{*}A_{U}^{\cM}\,\,\, \text{in }D^{b}\cM
(X)
\leqno(7.13.4)
$$
by the same argument as (7.13.2).
Here the cohomological functor
$ {}^{c}H^{j} $ is compatible with the natural functor
$ \cM(X{)}_{\gl \le n}^{\go} \to \cM(X)^{\go} \to
\cM(X) $, because so is the
$ t $-structure by definition.
The morphism (7.13.4) is injective in the heart of the
$ t $-structure of
$ D^{b}\cM(X) $, and its cokernel is supported in
$ Z := X \setminus U $.
Let
$ i : Z \to X $ be a natural morphism, and
$ N \in D^{b}\cM(Z) $ such that
$ i_{*}N $ is isomorphic to the cokernel, i.e.,
$$
i_{*}N = C(A_{X}^{\cM} \to{}^{c}H^{0}j_{*}A_{U}^{\cM})\,\,\,
\text{in }D^{b}\cM(X).
$$
Restricting
$ X $ to an open subvariety containing
$ U $, we may assume
$ Z $ is smooth and pure dimensional, and
$ \For(N), i^{*} \For({}^{c}H^{0}j_{*}A_{U}^{\cM}) $ are local systems on
$ Z $.
In particular,
$ N[-d], i^{*}({}^{c}H^{0}j_{*}A_{U}^{\cM})[-d] \in \cM(Z) $, where
$ d = \dim Z $.
Then, it is enough to show that
$ N[-d] \in \cM(Z{)}_{\gl \le n}^{\go} $ and the natural morphism
$ {}^{c}H^{0}j_{*}A_{U}^{\cM} \to i_{*}N $ is naturally lifted to a
morphism in
$ D^{b}\cM(X{)}_{\gl \le n}^{\go} $.
By adjunction, it is equivalent to consider a morphism
$$
i^{*}({}^{c}H^{0}j_{*}A_{U}^{\cM}) \to N.
$$
By hypothesis, it is a surjective morphism of
$ \cM(Z) $, and hence of
$ \cM(Z{)}_{\gl \le n}^{\go} $,
up to shift of complexes by
$ -d $, because
$ \cM(Z{)}_{\gl \le n}^{\go} $ is stable by subquotients in
$ \cM(Z) $.
So we get the assertion.

For the last assertion, it is enough to consider the restriction morphism
by duality.
We have the injectivity of
$$
\Hom_{gl\le n}(A_{Y}^{\cM},{}^{c}H^{0}f_{*}A_{X}^{\cM}) \to
\Hom(A_{Y(\bC)},\cH^{0}f_{*}A_{X(\bC)})
\leqno(7.13.5)
$$
by (7.12.1).
Since
$ \cH^{k}f_{*}A_{X(\bC)} = 0 $ for
$ k < 0 $,
we may replace
$ {}^{c}H^{0}f_{*}A_{X}^{\cM}, \cH^{0}f_{*}A_{X(\bC)} $ by
$ f_{*}A_{X}^{\cM}, f_{*}A_{X(\bC)} $.
So it is enough to show that the natural restriction morphism
$ v : A_{Y(\bC)} \to \cH^{0}f_{*}A_{X(\bC)} $ is lifted to a morphism of
$ D^{b}\cM(Y{)}_{\gl \le n}^{\go} $.
Here we may assume
$ f(X) $ is dense in
$ Y $ by replacing
$ Y $ with the closure of
$ f(X) $, because the assertion is clear in the closed embedding case
by (7.11).
Note that, if
$ Y $ is smooth and
$ \cH^{0}f_{*}A_{X(\bC)} $ is a local system, we have the assertion,
because
$ v $ is lifted to a morphism of
$ D^{b}\cM(Y) $ by (3.9.2), and
$ \cM(Y{)}_{\gl \le n}^{\go} $ is a full subcategory of
$ \cM(Y) $.

Let
$ U $ be a smooth open dense subvariety of
$ Y $ such that the restriction of
$ \cH^{0}f_{*}A_{X(\bC)} $ to
$ U $ is a local system.
By adjunction, we have a canonical morphism
$$
v' : A_{Y}^{\cM} \to{}^{c}\cH^{0}j_{*}j^{*}(^{c}\cH^{0}f_{*}
A_{X}^{\cM}),
$$
corresponding to the restriction morphism for the restriction of
$ f $ over
$ U $.
We have also a natural injective morphism in
$ \cC $ on
$ Y $:
$$
{}^{c}\cH^{0}f_{*}A_{X}^{\cM} \to{}^{c}\cH^{0}j_{*}j^{*}(^{c}
\cH^{0}f_{*}A_{X}^{\cM}),
$$
because the underlying
$ A $-complex of its kernel is zero.
Let
$ Q $ denote the cokernel of the injection.
By the same argument as above, it is enough to show that the composition of
$ v' $ with the projection to
$ Q $ is zero, and the assertion follows from the compatibility with the
functor
$ \For $.

\medskip\noindent
{\it Remark.} If
$ X_{\bC} $ is connected, the injectivity of (7.13.1) implies
$$
\End(A_{X}^{\cM}) = A
\leqno(7.13.6)
$$
But this is not true in general (depending on
$ \cM) $ even if
$ X $ is connected.
Here we may assume
$ X $ smooth and connected by (7.13.3), restricting
$ X $ to its open dense subvariety.
Then (7.13.6) is equivalent to the
$ \cM $-Hodge conjecture for
$ p = 0 $ (cf. (8.5)).

\medskip\noindent{{\bf 7.14.~Proposition.}} {\it
Let
$ X \in \cV(k) $, and
$ X' $ be the union of the irreducible components
$ X_{i} $ of
$ X $ such that
$ \dim X_{i} = d := \dim X $.
Let
$ \IC_{X'}A^{\cM} = \mopls_{i} \IC_{X_{i}}A^{\cM} $ with
$ \IC_{X_{i}}A^{\cM} = (j_{i})_{!*}A_{{U}_{i}}^{\cM}[d] \in
\cM(X{)}_{\gl \le n}^{\go} $ where
$ j_{i} : U_{i} \to X $ is the inclusion of a dense smooth open subvariety
$ U_{i} $ of
$ X' $ into
$ X $ \(cf. {\rm (6.3))}.
Then
$ \IC_{X'}A^{\cM}, \IC_{X_{i}}A^{\cM} $ are self-dual up to Tate twist,
i.e.,
$$
\bD(\IC_{X'}A^{\cM}) = \IC_{X'}A^{\cM}(d)\,\,\, \text{(same for }
\IC_{X_{i}}A^{\cM}).
\leqno(7.14.1)
$$
Moreover, we have a canonical morphism
$$
A_{X}^{\cM}[d] \to \IC_{X'}A^{\cM}\,\,\, \text{in }D^{b}\cM
(X{)}_{\gl \le n}^{\go}
\leqno(7.14.2)
$$
which induces isomorphisms
$$
\aligned
\Hom(A_{X}^{\cM}[d], (\bD A_{X}^{\cM})(-d)[-d])
&\buildrel{\sim}\over\leftarrow \Hom(\IC_{X'}A^{\cM},
(\bD A_{X}^{\cM})(-d)[-d])
\\
\buildrel{\sim}\over\leftarrow \Hom(\IC_{X'}A^{\cM}, \IC_{X'}A^{\cM})
&= \mopls_{i} \Hom(\IC_{X_{i}}A^{\cM}, \IC_{X_{i}}A^{\cM})
\endaligned
\leqno(7.14.3)
$$
where
$ \Hom $ can be taken in any of
$ D^{b}\cM(X{)}_{\gl \le n}^{\go}, D^{b}\cM(X)^{\go} $ and
$ D^{b}\cM(X) $.
}

\medskip\noindent{\it Proof.}
The first assertion is clear by definition of intermediate direct image
and (4.3).
The remaining assertions follow from:

\medskip\noindent{{\bf 7.15.~Lemma.}} {\it
With the above notation, we have
}
$$
H^{k}A_{X}^{\cM} = 0\,\,\, \text{for }k > d
\leqno(7.15.1)
$$
$$
\Gr_{k}^{W}H^{d}A_{X}^{\cM} = 0\,\,\, \text{for}\quad k > d
\leqno(7.15.2)
$$
$$
\Gr_{d}^{W}H^{d}A_{X}^{\cM} = \IC_{X'}A^{\cM}.
\leqno(7.15.3)
$$

\medskip\noindent
{\it Proof.} For (7.15.1), it is enough to show
$^{p}\cH^{k}A_{X(\bC)} = 0 $ for
$ k > d $.
But it follows from the definition of perverse sheaf (cf. [3]).
The second assertion is clear by (6.7).
The last assertion is reduced to
$$
\Hom(A_{X}^{\cM}[d], M) = 0
\leqno(7.15.4)
$$
for
$ M \in \cM(X) $ such that
$ \dim \supp M < d $, because
$ \Gr_{d}^{W}H^{d}A_{X}^{\cM} $ is semisimple and the restriction of
(7.15.3) to
$ U $ is clear.
Let
$ Z = \supp M $.
By adjunction we may replace
$ A_{X}^{\cM}[d] $ by
$ A_{Z}^{\cM}[d] $.
Then the assertion follows from (7.15.1) for
$ Z $.

\medskip\noindent
{\bf 7.16.}~{\it Remark.} Let
$ X \in \cV(k) $, and
$ n = \dim X $.
By (7.14.2) and its dual, we have a canonical morphism in
$ D^{b}\cM(X{)}_{\gl \le n}^{\go} $ (and in
$ D^{b}\cM(X)) $:
$$
A_{X}^{\cM}(n)[2n] \to \bD A_{X}^{\cM}
\leqno(7.16.1)
$$
which corresponds to the identity on
$ \IC_{X'}A^{\cM} $ by (7.14.3).
Let
$ \pi : Y \to X $ be a proper morphism such that
$ Y $ is smooth of pure dimension
$ n $ and
$ \pi (Y) = X' $ in (7.14).
Then (7.16.1) coincides with the composition of the restriction and Gysin
morphisms:
$$
A_{X}^{\cM}(n)[2n] \to \pi_{*}A_{Y}^{\cM}(n)[2n] \to \bD
A_{X}^{\cM},
\leqno(7.16.2)
$$
where (4.3.1) is also used.
In fact, the morphisms of (7.14.3) are uniquely determined by their
restriction to
$ U $ in (7.4), because of (6.3.4).

\bigskip\bigskip
\centerline{\bf \S 8. Cycle Class}

\bigskip

From now on, we assume
$ \cM(X) = \cM(X)^{\go} $, by replacing
$ \cM(X) $ by its full subcategory
$ \cM(X)^{\go} $ (cf. (7.1)).

\medskip\noindent
{\bf 8.1.~Definition.} For
$ X \in \cV(k) $ with natural morphism
$ a_{X} : X \to \Spec k $, let
$$
\aligned
{H}_{i}^{\cM}(X,A(j))
&= \Hom(A_{X}^{\cM}, \bD A_{X}^{\cM}(-j)[-i])
\\
&=\Hom(A^{\cM}, (a_{X})_{*}\bD A_{X}^{\cM}(-j)[-i])
\endaligned
\leqno(8.1.1)
$$
$$
\aligned
{H}_{M}^{i}(X,A(j))
&= \Hom(A_{X}^{\cM}, A_{X}^{\cM}(j)[i])
\\
&=\Hom(A^{\cM}, (a_{X})_{*}A_{X}^{\cM}(j)[i])
\endaligned
\leqno(8.1.2)
$$
where
$ \Hom $ is taken in
$ D^{b}\cM(X) $ or
$ D^{b}\cM(\Spec k) $, and the last isomorphisms of (8.1.1--2) are
induced by adjunction.

For a proper morphism
$ f : X \to Y $, the direct image
$$
f_{\#} : {H}_{i}^{\cM}(X,A(j)) \to {H}_{i}^{\cM}(Y,A(j))
\leqno(8.1.3)
$$
is defined by the composition with the restriction and Gysin morphisms
(3.9.2), using the first expression of (8.1.1).

For a morphism of smooth varieties
$ f : X \to Y $, the pull-back
$$
f^{\#} : {H}_{M}^{i}(Y,A(j)) \to {H}_{M}^{i}(X,A(j))
\leqno(8.1.4)
$$
is defined by applying the functor
$ f^{*} $ to the first expression of (8.1.2).

\medskip\noindent
{\it Remarks.} (i) By (7.15.1), we have
$$
{H}_{i}^{\cM}(X,A(j)) = 0\,\,\, \text{for }i > 2 \dim X.
\leqno(8.1.5)
$$

\medskip

(ii) If
$ X $ is smooth and pure dimensional, we have by (4.3.1):
$$
{H}_{i}^{\cM}(X,A(j)) = {H}_{M}^{2dimX-i}(X, A(\dim X-j)).
\leqno(8.1.6)
$$

\medskip

(iii) Assume
$ X $ smooth proper.
Then we have a decomposition
$$
(a_{X})_{*}A_{X}^{\cM} \simeq \mopls_{i} (H^{i}(a_{X})_{*}A_{X}^{\cM})
[-i]\quad \text{in }D^{b}\cM(\Spec k)
\leqno(8.1.7)
$$
by (6.10.1) applied to
$ \cM(\Spec k) $.
So the canonical filtration
$ \tau (cf. [9]) $ on
$ (a_{X})_{*}A_{X}^{\cM} $ induces a decreasing filtration
$ L $ on
$ {H}_{M}^{i}(X,A(j)) $ such that
$$
\Gr_{L}^{k}{H}_{M}^{i}(X,A(j)) = \Ext^{k}(A^{\cM}, H^{i-k}(a_{X})_{*}
A_{X}^{\cM}(j)).
\leqno(8.1.8)
$$
where
$ \Ext^{k} $ is taken in
$ D^{b}\cM(\Spec k) $.

(iv) Using the second expression of (8.1.1--2), the direct image and
pull-back are defined by the composition with the Gysin morphism and the
restriction morphism respectively.

(v) For a closed subvariety
$ Z $ of
$ X $ with natural inclusion
$ i : Z \to X $, define
$$
{H}_{Z,M}^{i}(X,A(j)) = \Hom(A^{\cM}, (a_{Z})_{*}i^{*}A_{X}^{\cM}(j)[i]).
\leqno(8.1.9)
$$
We can show that
$ {H}_{i}^{\cM}(X,A(j)) $ and
$ {H}_{Z,M}^{i}(X,A(j)) $ form a twisted Poincar\'e duality theory [7].

\medskip\noindent{{\bf 8.2.~Theorem.}} {\it
Let
$ X \in \cV(k) $, and
$ \CH_{d}(X, A) $ the Chow group of
$ X $ with
$ A $-coefficients.
Then we have naturally a cycle class map:
$$
cl^{\cM} : \CH_{d}(X, A) \to {H}_{2d}^{\cM}(X,A(d))
\leqno(8.2.1)
$$
compatible with the direct image by a proper morphism and the pull-back by
a morphism of smooth varieties.
}

\medskip\noindent{\it Proof.}
Let
$ Z $ be a reduced irreducible closed subvariety of
$ X $ with dimension
$ d $.
Composing (7.16.1) for
$ Z $ with the restriction and Gysin morphisms (3.9.2), we get
$$
A_{X}^{\cM}(d)[2d] \to \bD A_{X}^{\cM},
\leqno(8.2.2)
$$
which defines
$ cl^{\cM}([Z]) \in {H}_{2d}^{\cM}(X,A(d)) $.

For well-definedness of
$ cl^{\cM} $, consider a reduced irreducible closed subvariety
$ Z' $ of
$ Y := X\times S $ with dimension
$ d + 1 $, where
$ S = \Spec k[t] $.
Then
$ Z' $ defines
$$
A_{Y}^{\cM}(d+1)[2d+2] \to \bD A_{Y}^{\cM}
\leqno(8.2.3)
$$
by (8.2.2).
Let
$ P_{0}, P_{1} $ be closed points of
$ S $ defined by
$ t = 0 $ and
$ t = 1 $ respectively.
Let
$ i_{a} : X\times \{P_{a}\} \to Y, \pi : Y \to X $ be natural
morphisms.
Then we have functorial morphisms
$ \pi_{*} \to (i_{a})* $ by (3.3.2), which induce isomorphisms on
$ A_{Y}^{\cM}, \bD A_{Y}^{\cM} $.
So it is enough to show that
$ (i_{0})* $ of (8.1.4) is induced by
$ \zeta := [Z']^{\ssbull}[X\times \{P_{0}\}] $ (the argument is same for
$ a = 1) $.
Let
$ Z = Z' \cap X\times \{P_{0}\} $.
Then the assertion is reduced to the commutative diagram:
$$
\CD
A_{Z'}^{\cM} @>>> \bD A_{Z'}^{\cM}(-d-1)[-2d-2] @>>>
\bD A_{Y}^{\cM}(-d-1)[-2d-2]
\\
@VVV @. @VVV
\\
A_{Z}^{\cM} @>{\alpha}>> \bD A_{Z}^{\cM}(-d)[-2d] @>>>
\bD A_{X}^{\cM}(-d)[-2d].
\endCD
\leqno(8.2.4)
$$
where
$ X $ is identified with
$ X\times \{P_{0}\} $, the vertical morphisms are induced by the restriction
morphisms, and
$ \alpha $ is induced by the coefficient of
$ \zeta $ using (7.14.3).
By adjunction,
$ \alpha $ is uniquely determined by the diagram, and the assertion is
reduced to that for the underlying
$ A $-complexes, using (7.14.3).
So we may assume
$ k = \bC $.
Let
$ p : Y \to S $ denote a natural projection.
Applying the functorial morphism
$ (i_{0})* \to \psi_{p,1} $ to the composition of the morphisms
of the first row of (8.2.4), we get a commutative diagram:
$$
\CD
(i_{0})^{*}A_{Z'} @>{\beta}>> (i_{0})^{*}\bD A_{Y}(-d-1)[-2d-2]
\\
@VVV @VV{\simeq}V
\\
\psi_{p,1}A_{Z'} @>>> \psi_{p,1}\bD A_{Y}(-d-1)[-2d-2]
\endCD
\leqno(8.2.5)
$$
Here
$ \beta $ coincides the composition of
$ \alpha $ with the Gysin morphism.
Let
$ m_{i} $ be the coefficient of
$ \alpha $ at a general point of an irreducible component
$ Z_{i} $ of
$ Z $, (i.e.,
$ \alpha = m_{i}id) $.
Then
$ m_{i} $ coincides with the intersection multiplicity of
$ Z' $ and
$ X\times \{P_{0}\} $, because they coincide with the number of connected
components of the Milnor fiber of
$ p : Z' \to S $.
So we get the assertion.

The compatibility with direct image is clear by definition (8.1.3)
(cf. [30, II]).
The compatibility with pull-back is checked as in [loc. cit.] using
specialization, and reducing to the case of closed embedding of
codimension one as above.

\medskip\noindent
{\it Remark.} If
$ X $ is smooth, (8.2.1) becomes
$$
cl^{\cM} : \CH^{p}(X, A) \to {H}_{M}^{2p}(X,A(p))
\leqno(8.2.6)
$$
by (8.1.6).
For
$ X $ smooth proper, let
$ L $ be the filtration on
$ \CH^{p}(X) $ induced by the filtration
$ L $ on
$ {H}_{M}^{2p}(X,A(p)) $ in (8.1.8) via the cycle map
$ cl^{\cM} $ in (8.2.6).
Then the cycle map induces an injective morphism
$$
\Gr_{L}^{k}cl^{\cM} : \Gr_{L}^{k}\CH^{p}(X) \to \Ext^{k}(A^{\cM},
H^{2p-k}(a_{X})_{*}A_{X}^{\cM}(p))
\leqno(8.2.7)
$$
by (8.1.8).
The existence of such a filtration was suggested by Bloch [5].
We can show
$$
cl^{\cM}(\CH^{p}(X, A)) \cap L^{p+1}{H}_{M}^{2p}(X,A(p)) = 0
\leqno(8.2.8)
$$
by reducing to the smooth projective case, and using the weak Lefschetz
theorem (cf. [30, II]).
This implies
$$
\Gr_{L}^{k}\CH^{p}(X) = 0\quad \text{for}\quad k > p,
\leqno(8.2.9)
$$
and the separatedness of
$ L $ is equivalent to the injectivity of the cycle map.

\medskip\noindent{{\bf 8.3.~Proposition.}} {\it
Let
$ \CH_{d}(X,n) = \CH^{p}(X,n) $ be Bloch's higher Chow group {\rm [6]}
for
$ d = \dim X - p $, and
$ \CH_{d}(X,n)_{A} $ its scalar extension.
Then we have a cycle class map
}
$$
cl^{\cM} : \CH_{d}(X,n)_{A} \to {H}_{2d+n}^{\cM}(X,A(d)).
\leqno(8.3.1)
$$

\medskip\noindent
{\it Remark.} If
$ X $ is smooth and pure dimensional, (8.3.1) becomes
$$
cl^{\cM} : \CH^{p}(X,n)_{A} \to {H}_{M}^{2p-n}(X,A(p))
\leqno(8.3.2)
$$
by (8.1.6), where
$ p = \dim X - d $.

\medskip\noindent{\it Proof of (8.3).}
We can construct a cycle map of
$ \bigcap_{i} \Ker\,\partial_{i} $ in [loc. cit.], applying the next Lemma
to
$ X\times \Delta^{n}, X\times S^{n-1} $ in the notation of [loc. cit.],
because (7.6.1) implies:
$$
(a_{\Delta^{n}})_{*}j_{!}A_{U}^{\cM} = A^{\cM}[-n],
\leqno(8.3.3)
$$
where
$ j : U := \Delta^{n} \setminus S^{n-1} \to \Delta^{n} $.
Then well-definedness of the map is reduced to the invariance of the map by
a deformation parametrized by
$ \Spec k[t] $, and is checked by the same argument as in the proof of
(8.2).

\medskip\noindent{{\bf 8.4.~Lemma.}} {\it
Let
$ X = X'\times S, Y = X'\times D $ for
$ X', S, D \in \cV(k) $ such that
$ S $ is smooth, and
$ D $ is a divisor with normal crossings on
$ S $ whose irreducible components
$ D_{i} $ are smooth.
Let
$ Y_{i} = X'\times D_{i}, Y_{I} = \bigcap_{i\in I} Y_{i} $, and
$ i : Y \to X, j : U := X \setminus Y \to X $ denote natural
inclusions.
Let
$ \zeta = \sum_{k} n_{k}[Z_{k}] $ a cycle of dimension
$ d $ on
$ X $, where
$ Z_{k} $ are reduced irreducible closed subvarieties of
$ X $, which intersect properly with
$ Y_{I} $ for any
$ I $.
Let
$ Z = \bigcup_{k}Z_{k} $.
Then a cycle
$ \zeta $ induces a morphism
$$
A_{Z}^{\cM} \to \bD A_{X}^{\cM}(-d)[-2d]
\leqno(8.4.1)
$$
as in {\rm (8.2.2)} using {\rm (7.14.3)} and the Gysin morphism,
where the direct image by closed embedding is omitted
\(cf. Remark after {\rm (5.6))}.
This morphism is uniquely lifted to a morphism
$$
A_{Z}^{\cM} \to j_{!}\bD A_{U}^{\cM}(-d)[-2d],
\leqno(8.4.2)
$$
if
$ \zeta^{\ssbull}[Y_{i}] = 0 $ in
$ Z \cap Y_{i} $ for any
$ i $.
}

\medskip\noindent{\it Proof.}
It is enough to show that the composition of (8.4.1) with the natural
morphism
$$
\bD A_{X}^{\cM}(-d)[-2d] \to i^{*}\bD A_{X}^{\cM}(-d)[-2d]
\leqno(8.4.3)
$$
is zero, and
$$
\Hom(A_{Z}^{\cM}, i^{*}\bD A_{X}^{\cM}(-d)[-2d-1]) = 0.
\leqno(8.4.4)
$$
Let
$ n = \dim S $ (where we may assume
$ S $ is connected).
By (3.1.9) (4.3.1), we have
$$
i^{*}\bD A_{X}^{\cM} = \bD A_{X'}^{\cM}\,\boxtimes\,A_D^{\cM}(n)[2n].
\leqno(8.4.5)
$$
By (7.6.1), we have
$$
\Gr_{k}^{W}(A_D^{\cM}[n-1]) = \mopls_{|I|=n-k} A_{{D}_{I}}^{\cM}[k]
\quad \text{for }0 \le k \le n - 1,
\leqno(8.4.6)
$$
and it is zero otherwise.
Since
$ A_{{D}_{I}}^{\cM}[k] = \bD A_{{D}_{I}}^{\cM}(-k)[-k] $, we get a
decreasing filtration
$ G $ on
$ i^{*}\bD A_{X}^{\cM} $ such that
$$
\Gr_{G}^{p}(i^{*}\bD A_{X}^{\cM}) = \mopls_{|I|=p}
\bD A_{{Y}_{I}}^{\cM}(p)[p+1]\quad \text{for }1 \le p \le n.
\leqno(8.4.7)
$$
where
$ p = n - k $ is the codimension of
$ Y_{I} $.
The composition of (8.4.1) (8.4.3) and the projection
$$
i^{*}\bD A_{X}^{\cM}(-d)[-2d] \to \Gr_{G}^{1}(i^{*}\bD A_{X}^{\cM})
(-d)[-2d]
\leqno(8.4.8)
$$
is zero by
$ \zeta^{\ssbull}[Y_{i}] = 0 $ using the compatibility the cycle map with
pull-back (cf. (8.2)).
Using (8.4.7) (7.15.1) and adjunction, we can check
$$
\Hom(A_{Z}^{\cM}, \Gr_{G}^{p}(i^{*}\bD A_{X}^{\cM})(-d)[-2d]) = 0\quad
\text{for}\quad p > 1
\leqno(8.4.9)
$$
and (8.4.4), because
$ Z $ intersects properly with
$ Y_{I} $.

\medskip\noindent
{\bf 8.5.~Definition.} Let
$ X $ be a smooth proper variety in
$ \cV(k) $.
An element of
$$
\Hom(A^{\cM}, H^{2p}(a_{X})_{*}A_{X}^{\cM}(p))
\leqno(8.5.1)
$$
is called an
$ \cM $-{\it Hodge cycle.} We say that the
$ \cM $-{\it Hodge type conjecture} is true if
$ \Gr_{L}^{0}cl^{\cM} $ in (8.2.7) is surjective.

\medskip\noindent
{\it Remarks.} (i) Let
$ X $ be as above.
We have a birational morphism
$ \pi : X' \to X $ such that
$ X' $ is smooth projective.
Then the
$ \cM $-Hodge type conjecture is true for
$ X $ if it is true for
$ X' $.
In fact, the composition of the restriction morphism
$ A_{X}^{\cM} \to \pi_{*}A_{X'}^{\cM} $ and the restriction morphism
$ \pi_{*}A_{X'}^{\cM} \to A_{X}^{\cM} $ is the identity by the
faithfulness of
$ \For : \cM(X) \to \Perv(X(\bC),A) $, and the assertion follows from the
compatibility of the cycle map with direct image.

\medskip

(ii) If
$ k = \bC, A = \bQ $ and
$ \cM(X) = \MHM(X,A) $ as in the example (i) of (1.8), then
$ \cM $-Hodge cycle is same as Hodge cycle, and
$ \cM $-Hodge type conjecture as Hodge conjecture.

\medskip

(iii) If
$ A = \bQ $ and
$ \cM(X) $ is as in the example (iii) (iv) or
$ (v) $ of (1.8), then the
$ \cM $-Hodge type conjecture is reduced to the Hodge conjecture for
$ X_{\bC} $, and is true for
$ p \le 1 $.
In fact, it is enough to consider the example (iii).
If the Hodge conjecture for
$ X_{\bC} $ is true, the image of an
$ \cM $-Hodge cycle in
$ H^{2p}(X(\bC), \bQ(p)) $ belongs to the image of an algebraic
cycle which is defined over a field
$ K $ finitely generated over
$ k $.
Then we can reduce to the case
$ K $ is algebraic over
$ k $, taking a finitely generated
$ k $-algebra
$ R $ in
$ K $, which generates
$ K $, and using the reduction at a sufficiently general closed point
of
$ \Spec R $.
In this case, we may assume
$ K $ is a Galois extension.
Then the assertion is proved by taking the average of the images of the
algebraic cycle by the action of Galois group.
In fact, the algebraic cycle defines the cycle classes in
$ H^{2p}(X(\bC), \bQ(p)) $ and
$ {H}_{\acute et}^{2p}(\overline{X}, \bQ_{l}(p)) $ compatible with the
comparison isomorphism, and its \'etale part is
$ G $-invariant, because it coincides with the original
$ \cM $-Hodge cycle by its coincidence on
$ H^{2p}(X(\bC), \bQ(p)) $.
Compare to [20].

\medskip

(iv) By a similar argument, we see that the natural morphism
$$
\CH_{d}(X,A) \to \CH_{d}(X_{K},A)
\leqno(8.5.2)
$$
is injective for a field extension
$ k \to K $, where
$ X_{K} = X\otimes_{k}K $.

\medskip\noindent
{\bf 8.6.~Definition.} Let
$ X, Y $ be smooth proper varieties.
The group of correspondences of
$ X $ to
$ Y $ with
$ A $-coefficients is defined by
$$
C^{i}(X,Y) = \CH^{p}(X\times Y,A)\quad \text{for}\quad p = \dim X + i.
\leqno(8.6.1)
$$
The composition is defined by
$$
\xi\scirc \zeta = (p_{13})_{*}((p_{12})^{*}\xi^{\ssbull}
(p_{23})^{*}\zeta )
\leqno(8.6.2)
$$
for
$ \xi \in C^{i}(X,Y), \zeta \in C^{j}(Y,Z) $, where
$ p_{ij} $ is the projection of
$ X\times Y\times Z $ to the product of the
$ i^{th} $ and
$ j^{th} $ components (cf. [23] [25]).

\medskip\noindent
{\it Remark.} We have canonical isomorphisms
$$
\aligned
{H}_{M}^{2p}(X\times Y, A(p))
&= \Ext^{2i}(A^{\cM}, (a_{X})_{*}\bD A_{X}^{\cM} \otimes
(a_{Y})_{*}A_{Y}^{\cM}(i))
\\
&=\Ext^{2i}((a_{X})_{*}A_{X}^{\cM}, (a_{Y})_{*}A_{Y}^{\cM}(i)).
\endaligned
\leqno(8.6.3)
$$
where the first isomorphism follows form (2.7.5) (4.3.1), and the second
from the dual of (4.11.1) on
$ \Spec k $.
This isomorphism is compatible with the filtration
$ L $ in (8.1.8), where
$ L $ on the last term is defined by using the decomposition (8.1.7) so
that
$$
\Gr_{L}^{k} = \mopls_{j} \Ext^{k}(H^{j}(a_{X})_{*}A_{X}^{\cM},
H^{j+2i-k}(a_{Y})_{*}A_{Y}^{\cM}(i) ,
\leqno(8.6.4)
$$
cf. [30, II, (4.1)].
Combined with (8.2), we get
$$
\lambda : C^{i}(X,Y) \to \Ext^{2i}((a_{X})_{*}A_{X}^{\cM}, (a_{Y})_{*}
A_{Y}^{\cM}(i) .
\leqno(8.6.5)
$$

By the same argument as in [30, II], we can show the following two
propositions:

\medskip\noindent{{\bf 8.7.~Proposition.}} {\it
The morphism {\rm (8.6.5)} is compatible with the composition of cycle,
where the composition on the target is defined by the composition of
morphisms.
}

\medskip\noindent{{\bf 8.8.~Proposition.}} {\it
Let
$ Z $ be a reduced irreducible closed subvariety of
$ X\times Y $.
Then the image of
$ [Z] $ by {\rm (8.6.5)} coincides with the composition of {\rm (7.16.1)}
\(applied to
$ Z) $ with the restriction and Gysin morphisms for the natural morphism
of
$ Z $ to
$ X $ and
$ Y $ respectively.
It coincides also with the composition of restriction and Gysin morphisms
for the morphism of
$ Z' $ to
$ X $ and
$ Y $ respectively, if
$ Z' \to Z $ is a resolution of singularity,
}

\medskip\noindent
Using these we can show the following (cf. [loc. cit.]):

\medskip\noindent{{\bf 8.9.~Theorem.}} {\it
The cycle map
$ cl^{\cM} $ in {\rm (8.2.1)} is surjective if the
$ \cM $-Hodge type conjecture is true for any smooth projective varieties
over
$ k $.
}

\medskip\noindent{{\bf 8.10.~Corollary}} {\it
The morphism
$ \lambda $ in {\rm (8.6.5)} is surjective if the
$ \cM $-Hodge type conjecture is true for any smooth projective varieties
over
$ k $.
}

\medskip\noindent
{\it Remarks.} (i) By definition of geometric origin, we need the
$ \cM $-Hodge type conjecture for {\it any} smooth projective varieties.

\medskip

(ii) We can consider the relative version of (8.6.5), and its surjectivity
is also reduced to the
$ \cM $-Hodge type conjecture (cf. [31]).
This implies the surjectivity of the natural morphisms
$$
\Ext_{\gl \le n}^{i}(M,N) \to \Ext_{\gl \le m}^{i}(M,N) \to
\Ext^{i}(M,N)
\leqno(8.10.1)
$$
for
$ m > n $ and
$ M \in \cM(X,w+i{)}_{\gl \le n}^{\go}, N \in \cM(X,w{)}_{\gl \le n}^{\go} $,
where
$ \Ext_{\gl \le n}^{i}, \Ext^{i} $ are taken in
$ D^{b}\cM(X{)}_{\gl \le n}^{\go}, D^{b}\cM(X)^{\go} = D^{b}\cM
(X) $ respectively.

\medskip\noindent{{\bf 8.11.~Proposition.}} {\it
For
$ X \in \cV(k) $ with
$ \dim X \le n $, we have naturally a cycle map
$$
cl^{\cM} : \CH_{d}(X, A) \to \Hom_{gl\le n}(A_{X}^{\cM}, \bD A_{X}^{\cM}
(-d)[-2d])
\leqno(8.11.1)
$$
whose composition with the natural morphism
$$
\Hom_{gl\le n}(A_{X}^{\cM}, \bD A_{X}^{\cM}(-d)[-2d]) \to
\Hom(A_{X}^{\cM},\bD A_{X}^{\cM}(-d)[-2d])
\leqno(8.11.2)
$$
coincides with {\rm (8.2.1)}.
}

\medskip\noindent{\it Proof.}
This follows from (7.13) (7.14.3).

\medskip
As to the injectivity of cycle map, we have the following by the same
argument as in [31], where the injectivity of (8.2.1) in the divisor
case is reduced to the case
$ X = \bC $ using the hypothesis (cf. Remark (iii) below).

\medskip\noindent{{\bf 8.12.~Theorem.}} {\it
Assume A
$ = \bQ $ and the functor
$ \For : \cM(X) \to \Perv(X(\bC),A) $ is factorized by
$ \MHM(X_{\bC},A) $ in a compatible way with the functors and morphisms
in {\rm \S 1}.
Let
$ X $ be a pure dimensional locally complete intersection quasiprojective
variety in
$ \cV(k) $.
Then {\rm (8.11.1)} is injective for
$ n = \dim X $, if {\rm (8.2.1)} for
$ d $ is injective for any locally complete intersection closed subvariety
$ Y $ of
$ X $ with pure dimension
$ \dim X - 1 $ and if the cycle map {\rm (8.2.1)} for
$ d + 1 $ is surjective for any smooth projective variety whose dimension
coincides with
$ \dim X $.
}

\medskip\noindent
{\it Remarks.} (i) We can replace the condition
$ X, Y $ locally complete intersection by
$ A_{X(\bC)}[\dim X], A_{Y(\bC)}[\dim Y] $ perverse sheaves
(cf. [loc. cit.]).

\medskip

(ii) By induction on dimension, the condition of (8.12) is replaced by
the following: (8.11.2) for
$ n = \dim Y $ is injective for any locally complete intersection closed
subvariety
$ Y $ of
$ X $ with pure dimension
$ < \dim X $, and the cycle map (8.2.1) for
$ d + 1 $ is surjective for any smooth projective variety of dimension
$ \le \dim X $.
If the
$ \cM $-Hodge type conjecture is true for any smooth projective variety
over
$ k $, then the second condition would be satisfied by (8.9).

\medskip

(iii) If
$ A = \bQ $ and the functor
$ \For $ is factorized by
$ \MHM(X_{\bC},A) $ as above (e.g. the examples of (1.8)), then the
injectivity of the cycle map (8.2.1) is reduced to the injectivity of the
cycle map on
$ X_{\bC} $ by the injectivity of (8.5.2).
In this case, the cycle map (8.11.1) is injective for
$ n = \dim X $ and
$ d \ge \dim X - 2 $ by [31, (0.6)], where
$ X $ is as in (8.12).

\medskip

As to the injectivity of (8.11.2), we have the following by the same
argument as [31, (2.14)]:

\medskip\noindent{{\bf 8.13.~Proposition.}} {\it
The morphism {\rm (8.11.2)} is injective for a variety of dimension
$ \le n $, if the natural morphisms
$$
\Hom_{gl\le n}(A_{Y}^{\cM}, A_{Y}^{\cM}(p)[2p-1]) \to
\Hom(A_{Y}^{\cM}, A_{Y}^{\cM}(p)[2p-1])
\leqno(8.13.1)
$$
$$
\Hom_{gl\le \dim Y}(A_{Y}^{\cM}, A_{Y}^{\cM}(p)[2p]) \to
\Hom_{gl\le n}(A_{Y}^{\cM},
A_{Y}^{\cM}(p)[2p])
\leqno(8.13.2)
$$
are surjective for any connected smooth \(locally closed\,\) subvariety
$ Y $ of
$ X $, where
$ p = \dim Y - d $.
}

\medskip\noindent
{\it Remarks.} (i) For
$ p \le 1, (8.13.1) $ is always surjective, and so is (8.13.2) for
$ p \le 0 $.
In fact, it is clear for
$ p \le 0 $, and (8.13.1) for
$ p = 1 $ follows from the definition of geometric level.
In particular, (8.11.2) is injective for
$ d = n - 1 $, because (8.13.2) is the identity if
$ p = 1 $ and
$ n = \dim X = \dim Y $.

\medskip

(ii) If the
$ \cM $-Hodge type conjecture is true for any smooth projective variety
over
$ k $, then (8.13.2) would be surjective by the surjectivity of the
morphisms of (8.10.1).
The surjectivity of (8.13.1) would be reduced to the surjectivity of
the cycle map (8.3.2) for
$ n = 1 $, because an element of
$ \CH^{p}(X,n)_{A} $ is represented by a cycle of dimension
$ \dim X + n - p $.

\medskip\noindent
{\bf 8.14.}~{\it Remark.} Let
$ X $ be a reduced irreducible variety, and
$ k(X) $ the field of rational functions on
$ X $.
We define
$$
{H}_{M}^{i}(k(X), A(j)) = \varinjlim_{U}{H}_{M}^{i}(U, A(j)),
\leqno(8.14.1)
$$
where
$ U $ runs over (smooth) nonempty open subvarieties of
$ X (cf. [7]) $.
If we have the surjectivity of the cycle map onto
$ {H}_{MH}^{i}(U, A(j)) $, then
$$
{H}_{MH}^{i}(k(X), A(j)) = 0\,\,\, \text{for }i > j,
\leqno(8.14.2)
$$
and
$ {H}_{MH}^{i}(k(X), A(i)) $ would be expressed in terms of Tate mixed
sheaves, and would be related with Milnor
$ K $-theory.

In the case
$ i = 2p, j = p > 0, (8.14.2) $ for irreducible subvarieties of a smooth
variety
$ X $ would imply the
$ \cM $-Hodge type conjecture for
$ X $ by induction on dimension, if the
$ \cM $-Hodge type conjecture is true for
$ p = 0 $.
Note that the last condition is satisfied for the examples (iii) (iv) or
$ (v) $ in (1.8).

In the case
$ i = 2p - 1, j = p > 1, (8.14.2) $ for irreducible subvarieties of a smooth
variety
$ Y $ would imply the surjectivity of (8.13.1), if we have an isomorphism
$$
\Gamma (X,{\cO}_{X}^{*}) \otimes_{\bZ} A \simto \Ext^{1}(A_{X}^{\cM},
A_{X}^{\cM}(1))
\leqno(8.14.3)
$$
for smooth varieties
$ X $ (cf. [31, (3.12)]) .
It is hoped that (8.14.3) would hold in the example
$ (v) $ of (1.8).

It is also interesting whether the cycle map
$$
L^{1}\CH^{1}(X,A) \to \Ext^{1}(A^{\cM}, H^{1}(a_{X})_{*}A_{X}^{\cM}(1))
\leqno(8.14.4)
$$
is surjective in the example
$ (v) $ of (1.8).
By (8.9), it is surjective if the absolute Hodge cycles are algebraic for
any smooth projective varieties over
$ k $.
Here we have the injectivity, because it can be reduced to the case
$ k = \bC $.

\end{document}